\documentclass[11pt]{article}
\usepackage[french]{babel}
\usepackage[ansinew]{inputenc}
\usepackage{enumerate}
\usepackage{amsfonts}
\usepackage{amssymb}
\usepackage{amsmath}
\usepackage{amscd}
\usepackage[all]{xy}
\newtheorem{theo}{Th\'eor\`eme}[section]

\newtheorem{lem}[theo]{Lemme}
\newtheorem{slem}[theo]{Sous-Lemme}
\newtheorem{cor}[theo]{Corollaire}
\newtheorem{rem}[theo]{Remarque}
\newtheorem{ex}[theo]{Exemple}
\newtheorem{asser}[theo]{Assertion}
\newtheorem{defi}[theo]{D\'efinition}
\def \endproof{\hfill$\square$ \bigskip}
\def\id{\textrm{\scriptsize id}}
\def\o{\textrm{O}}
\def\K{\mathbb{K}}
\def\D1{
\begin{equation*}
\raisebox{1.6cm}
{\xymatrix{ 
\ar[r]^-\id         & (\K,\o) \ar[r]^-\id \ar[d]^{\Gamma_i} & (\K,\o) \ar[r]^-\id \ar[d]^{\Gamma_{i-1}} & \cdots \ar[r]^-\id     & (\K,\o) \ar[r]^-\id \ar[d]^{\Gamma_1} & (\K,\o) \ar[d]^{\Gamma_0} \\
\ar[r]^-{\Pi_{i+1}} & (E_i,\o_i) \ar[r]^-{\Pi_i}            & (E_{i-1},\o_{i-1}) \ar[r]^-{\Pi_{i-1}}    & \cdots \ar[r]^-{\Pi_2} & (E_1,\o_1) \ar[r]^-{\Pi_1}            & (E_0,\o_0)                \\
\ar[r]^-{\omega_{i+1}}   & (Z_i,\o_i) \ar@{^{(}->}[u]_{J_i} \ar[r]^-{\omega_i} & (Z_{i-1},\o_{i-1}) \ar@{^{(}->}[u]_{J_{i-1}} \ar[r]^-{\omega_{i-1}}      & \cdots \ar[r]^-{\omega_2}   & (Z_1,\o_1) \ar@{^{(}->}[u]_{J_1} \ar[r]^-{\omega_1}              & (Z_0,\o_0) \ar@{^{(}->}[u]_{J_0}
}}
\end{equation*}
}

\begin{document}

\title{Sur quelques aspects de la g\'eom\'etrie de l'espace des arcs trac\'es sur un espace analytique}
\author{M. Hickel}
\date{ }\maketitle
\renewcommand{\abstractname}{Abstract}
\begin{abstract} Let $(X,x)$ be a germ of real or complex analytic space and $\mathcal{A}_{(X,x)}$ the space of germs of arcs on $(X,x)$. Let us consider $F_{x}: (X,x) \rightarrow (Y,y)$ a germ of a morphism and denote by $\mathcal{F}_{x}: \mathcal{A}_{(X,x)} \rightarrow \mathcal{A}_{(Y,y)}$ the induced morphism at the level of arcs. In this paper, we try to emphasize the analogies between the metric or local topological properties of $F_{x}$ and those of $\mathcal{F}_{x}$. We then define the notions of Nash sequence of multiplicities, Nash sequence of Hilbert-Samuel functions and Nash sequence of diagram of initial exponents  of $X$ along an arc $\varphi$, and study some of their basic properties. Some elementary connections between these notions and  motivic integration theory are also provided.
\end{abstract}
\section{Introduction}

Soient $\mathbb{K}$=$\mathbb{R}$ ou $\mathbb{C}$ et $(X,x)$ un germe de $\mathbb{K}$-espace analytique, d'alg\`ebre locale $\mathcal{O}_{X,x}$. Nous noterons par $\mathcal{A}_{(X,x)}$ l'espace des arcs trac\'es sur $(X,x)$, c'est \`a dire l'ensemble des germes de morphismes analytiques $\varphi$ de $(\mathbb{K},0)$ dans $(X,x)$, ou bien de fa\c{c}on \'equivalente l'ensemble des morphismes $\varphi^*$ de $\mathbb{K}$-alg\`ebre locale de $\mathcal{O}_{X,x}$ dans $\mathcal{O}_{1}=\mathbb{K} \{t\}$. De mani\`ere plus g\'en\'erale si $X$ est un espace analytique et $A$ un sous-ensemble de $\vert X \vert$, nous d\'esignerons par $\mathcal{A}_{(X,A)}$ (resp. $\mathcal{A}_{X}$) l'union $\coprod_{x \in A} \mathcal{A}_{(X,x)}$ (resp. $ \coprod_{x \in X} \mathcal{A}_{(X,x)}$). $\mathcal{A}_{X}$ est naturellement munie d'une distance ultram\'etrique not\'ee $\mathcal{D}_{X}$, d\'efinie de la mani\`ere suivante:

\vspace{5pt}

\begin{itemize}
\item[-] $\forall \varphi \in \mathcal{A}_{(X,x)}$, $\forall \phi \in \mathcal{A}_{(X,x')}$,  $x\not=x'$, on pose $\mathcal{D}_{X}(\varphi,\phi)=1$;
\item[-]$\forall \varphi, \phi \in \mathcal{A}_{(X,x)}$, on pose $\mathcal{D}_{X}(\varphi, \phi)=e^{-ord(\varphi - \phi)}$, o\`u $ord(\varphi - \phi)$ d\'esigne la valuation ou l'ordre de l'id\'eal de $\mathcal{O}_{1}$, $(\varphi^{*} -\phi^{*})(\mathcal{M}_{X,x})$, et $\mathcal{M}_{X,x}$ l'id\'eal maximal de $\mathcal{O}_{X,x}$.
\end{itemize} 

\vspace{5pt}

La topologie induite par $\mathcal{D}_{X}$ sur $\mathcal{A}_{X}$ sera dite la topologie de Krull sur $\mathcal{A}_{X}$.
Soit maintenant un germe de morphisme analytique $F_{x}: (X,x) \rightarrow (Y,y)$. Celui-ci induit un morphisme $\mathcal{F}_{x} : \mathcal{A}_{(X,x)} \rightarrow \mathcal{A}_{(Y,y)}$, o\`u pour $\varphi \in \mathcal{A}_{(X,x)}$, $\mathcal{F}_{x}(\varphi)$ est l'arc $(\mathbb{K}, 0) \rightarrow (X,x) \rightarrow (Y,y)$ obtenu par composition de $\varphi$ avec $F_{x}$. De fa\c{c}on plus g\'en\'erale si $F : X \rightarrow Y$ est un morphisme, nous noterons par $\mathcal{F} : \mathcal{A}_{X} \rightarrow \mathcal{A}_{Y}$ le morphisme induit au niveau de l'espace des arcs.

Dans cet article nous abordons deux th\`emes. Le premier consiste \`a tenter de souligner les analogies entre les propri\'et\'es topologiques locales d'un morphisme $F : X \rightarrow Y$, et celles du morphisme induit $\mathcal{F} : \mathcal{A}_{X} \rightarrow \mathcal{A}_{Y}$.

Nous donnons deux illustrations de ce jeu de miroirs. La premi\`ere con\-cer\-ne un r\'esultat du \`a J.C. Tougeron. Dans [T3], l'auteur pr\'ecit\'e donne des in\'eg\-a\-li\-t\'es de Lojasiewicz, avec exposant uniforme, par rapport aux fibres d'un morphisme analytique. Plus pr\'ecis\'ement:

\begin{theo}[J.C. Tougeron ]
Soit $F$ une application analytique d'un ouvert $\Omega_{n}$ de $\mathbb{R}^{n}$ dans un ouvert $\Omega_{p}$ de $\mathbb{R}^{p}$. Pour toute paire de compacts $K \subset \Omega_{n}$ et $K' \subset \Omega_{p}$, il existe une constante $\alpha \geq 0$ satisfaisant:
$$\forall y \in K', \hspace{3pt}\exists C_{y} > 0 \vert \forall x \in K,  \hspace{3pt}d_{p}(F(x),y) \geq C_{y}d_{n}(x,F^{-1}(y))^{\alpha}$$
o\`u $d_{p}$ (resp. $d_{n}$) d\'esigne la distance euclidienne sur $\mathbb{R}^{p}$ (resp. sur $\mathbb{R}^{n}$).
\end{theo}
Le th\'eor\`eme r\'esulte trivialement de l'assertion suivante:

\begin{asser}
Pour toute paire de compacts $K \subset \Omega_{n}$, $K' \subset \Omega_{p}$, il existe une constante $\alpha \geq 0$ satisfaisant:

$$\forall y \in K', \hspace{3pt}\forall x_{0} \in K,\hspace{3pt} \exists C_{(x_{0},y)}>0 \textit{ et } V_{x_{0}} \textit{ voisinage de } x_{0} \textit{ tel que:}$$
$$ \forall x \in V_{x_{0}}, \hspace{3pt}d_{p}(F(x),y)\geq C_{(x_{0},y)}d_{n}(x,F^{-1}(y))^{\alpha}$$
\end{asser}
Nous allons voir que l'on a un r\'esultat analogue, au niveau du morphisme induit au niveau des espaces des arcs, en rempla\c{c}ant la distance euclidienne par la distance ultram\'etrique d\'efinie plus haut. Nous d\'emontrerons en effet:
\begin{theo}
Soit $F: X \rightarrow Y$ un morphisme de $\mathbb{K}$-espace analytique et $\mathcal{F}: \mathcal{A}_{X} \rightarrow \mathcal{A}_{Y}$ le morphisme induit au niveau de l'espace des arcs. Alors pour toute paire de compacts $K \subset X$, $K' \subset Y$, il existe une constante $\alpha$ satisfaisant:
$$\forall \varphi \in \mathcal{A}_{(Y,K')}, \forall \phi_{0} \in \mathcal{A}_{(X,K)}, \hspace{3pt}\exists C_{(\phi_{0},\varphi)}>0 \textit{ et } \mathcal{V}{\phi_{0}} \textit{ voisinage de } \phi_{0} \textit{ tel que:}$$
$$ \forall \phi \in \mathcal{V}_{\phi_{0}}, \hspace{3pt}\mathcal{D}_{X}(\mathcal{F}(\phi),\varphi)\geq C_{(\phi_{0},\varphi)}\mathcal{D}_{Y}(\phi,\mathcal{F}^{-1}(\varphi))^{\alpha}.$$
\end{theo}
De la m\^eme fa\c{c}on que l'existence d'une majoration affine pour la fonction d'Artin-Greenberg (cf. [Gr]) \'equivaut \`a des in\'egalit\'es de Lojasiewicz au niveau de l'espace des arcs (cf. [H2]), le r\'esultat ci-dessus s'interpr\`ete comme une version partiellement uniforme du th\'eor\`eme de Greenberg, l'espace des param\`etres \'etant $ \mathcal{A}_{Y}$ (cf. section 3).

La seconde illustration des analogies pr\'ec\'edemment \'evoqu\'ees concerne la condition de rang de A.M. Gabrielov. Pour le r\^ole jou\'e par cette condition, nous renvoyons le lecteur \`a [B-M3], [Ga], [I1], [I2], [I3], [I4], [P] et plus particuli\`erement \`a [I4] et \`a sa bibliographie pour une historique.
Si $F_{x} : (X,x) \rightarrow (Y,y)$ est un germe de morphisme entre germes d'espaces irr\'eductibles, on d\'efinit le rang g\'en\'erique de $F_{x}$ par:
$$grk(F_{x})=\inf_{U}(\sup_{x'} \textit{ rang}( F_{x'}),\hspace{5pt} x' \in Reg(X\cap U))$$
o\`u $U$ parcourt une base de voisinage de $x$, $Reg(X\cap U)$ d\'esigne l'ensemble des points r\'eguliers de $X\cap U$, et $rang(F_{x'})$ est le rang de la matrice jacobienne de $F_{x'}$ en $x'$.
On dit que $F_{x}$ satisfait la condition de rang de A.M. Gabrielov si et seulement si:
$$grk(F_{x})= dim_{y}Y.$$
Pour $\mathbb{K}=\mathbb{C}$, il est relativement \'el\'ementaire de v\'erifier que cette condition peut se caract\'eriser topologiquement par l'\'equivalence des conditions suivantes:
\begin{itemize}
\item[-]$grk(F_{x})= dim_{y}Y$;
\item[-]$\forall U$ voisinage de $x$ assez petit, $\forall V$ ouvert de $Reg(X\cap U)$, $F_{x}(V)$ contient un ouvert de $Reg(Y)$;
\item[-]$\forall U$ voisinage de $x$ assez petit, $F_{x}(Reg(X\cap U))$ contient un ouvert de $Reg(Y)$.
\end{itemize}

\vspace{10pt}

D\'esignons maintenant par $(Sing(X),x)$ le germe en $x$ de l'ensemble des points singuliers de$(X,x)$. Soit alors $\mathcal{R}_{(X,x)}=\mathcal{A}_{(X,x)} - \mathcal{A}_{(sing(X),x)}$ l'ensemble des arcs trac\'es dans $(X,x)$, mais dont l'image n'est pas enti\`erement incluse dans $(Sing(X)$ $,\hspace{3pt}x)$. $\mathcal{R}_{(X,x)}$ joue comme nous le verrons plus loin le r\^ole de l'ensemble des points r\'eguliers de $\mathcal{A}_{(X,x)}$. Etant donn\'e un morphisme $F_{x}: (X,x) \rightarrow (Y,y)$, il semble naturel de se demander quand $\mathcal{F}_{x}(\mathcal{A}_{(X,x)})$ contient un ensemble repr\'esentatif de $\mathcal{A}_{(Y,y)}$, par exemple un ouvert de $\mathcal{R}_{(Y,y)}$.
L\`a encore le ph\'enom\`ene est similaire au niveau du morphisme $F_{x}$ et \`a celui de $\mathcal{F}_{x}$. On a en effet:

\begin{theo}
Pour $\mathbb{K}=\mathbb{C}$, soient $F_{x} : (X,x) \rightarrow (Y,y)$ un morphisme entre deux germes irr\'eductibles et $\mathcal{F}_{x} : \mathcal{A}_{(X,x)} \rightarrow \mathcal{A}_{(Y,y)}$ le morphisme induit.  Les propri\'et\'es suivantes sont \'equivalentes:
\begin{itemize}
\item[1)] $grk(F_{x})=dim_{y}Y$;
\item[2)] Pour tout ouvert $\mathcal{V}$ de $\mathcal{R}_{(X,x)}$, $\mathcal{F}_{x}(\mathcal{V})$ contient un ouvert de $\mathcal{R}_{(Y,y)}$;
\item[3)] $\mathcal{F}_{x}(\mathcal{R}_{(X,x)})$ contient un ouvert de $\mathcal{R}_{(Y,y)}$.
\end{itemize}
\end{theo}
Notre d\'emonstration est \'el\'ementaire dans le sens o\`u elle n'utilise pas l'existence de d\'esingularisation, ni aucun des r\'esultats concernant la condition de rang de A.M. Gabrielov  mentionn\'es ci-dessus. Par contre on peut retrouver partiellement \`a partir de $1.4$ certains de ces r\'esultats (cf. section 6).

\vspace{10pt}

Le deuxi\`eme th\`eme abord\'e dans cet article concerne la notion de suite des multiplicit\'es de Nash le long d'un arc. En un certain sens, nous cherchons des notions convenables de multiplicit\'e, de fonction de Hilbert-Samuel, et de diagramme des exposants initiaux de $\mathcal{A}_{(X,x)}$ en un point $\varphi$.
Dans [L-J], M. Lejeune-Jalabert introduisait de mani\`ere algorithmique, pour un germe d'hypersurface $(H,x)$, la notion de suite des multiplicit\'es de Nash de $(H,x)$ le long d'un arc. Dans notre travail [H1] \'egalement consacr\'e aux hypersurfaces, cette notion apparaissait sous un point de vue plus g\'eom\'etrique. Nous g\'en\'eralisons ici ce point de vue, pour d\'efinir les notions de suite des multiplicit\'es de Nash,  suite de Nash des fonctions de Hilbert-Samuel et suite de Nash des diagrammes des exposants initiaux de $(X,x)$ le long d'un arc $\varphi$, pour tout germe $(X,x)$.

Les d\'efinitions sont les suivantes. Consid\'erons un plongement $i_0$ de $(X,x)$ dans un $(\mathbb{K}^{n},0)$. Ceci induit un plongement de $\mathcal{A}_{(X,x)}$ dans $\mathcal{A}_{(\mathbb{K}^{n},0)}$. Soit $\varphi \in \mathcal{A}_{(\mathbb{K}^{n},0)}$. On pose:
\begin{itemize}
\item[-] $(E_{0},O_{0})=(\mathbb{K}\hspace{3pt},0)\times(\mathbb{K}^n,0)=(\mathbb{K}^{n+1},0)$;
\item[-] $(Z_{0},O_{0})=(\mathbb{K}\hspace{3pt},0)\times(X,x)$;
\item[-] $\Gamma_{0}: (\mathbb{K}\hspace{3pt},0) \rightarrow (E_{0},O_{0})$ le graphe de $\varphi$.
\end{itemize}

\vspace{10pt}

On consid\`ere alors le diagramme commutatif $(D_{1})$ suivant:
\D1
o\`u pour tout $i \geq 1$:
\begin{itemize}
\item[-] $\Pi_{i}$ est l'\'eclatement de $E_{i-1}$ de centre $O_{i-1}$;
\item[-] $\Gamma_{i}$ est le rel\`evement de $\Gamma_{i-1}$ \`a travers $\Pi_{i}$ et $O_{i}=\Gamma_{i}(0)$;
\item[-] $Z_{i}$ est le transform\'e strict de $Z_{i-1}$ par $\Pi_{i}$ et $\omega_{i}=\Pi_{i}\vert Z_{i}$, $J_i$ est le plonge\-ment induit (avec $J_0=(id,i_0)$).
\end{itemize}

\vspace{10pt}

Soient pour $i \in \mathbb{N}$:
\begin{itemize}
\item[-] $m_{i,\varphi}$ la multiplicit\'e de $Z_{i}$ au point $O_{i}$ (bien entendu si $O_{i}\not\in Z_{i}$, $m_{i,\varphi}=0)$;
\item[-] $H_{i,\varphi}$ la fonction de Hilbert-Samuel du germe $(Z_{i},O_{i})$ (bien entendu si $O_{i}\not\in Z_{i}$, $H_{i,\varphi}=0)$.
\end{itemize}

Soit maintenant $X=(X_{1},\ldots,X_{n})$ un syst\`eme de coordonn\'ees \`a l'origine de $\mathbb{K}^{n}$. Ecrivons:
$$\varphi (t) = \sum _{k=1}^{+ \infty} A_{k}.t^{k},\hspace{5pt} A_{k} \in \mathbb{K}^{n}, \hspace{10pt} A_{k}=(a_{k}^{1},\ldots,a_{k}^{n}).$$
L'\'ecriture en coordonn\'ees locales des transformations quadratiques, comme on le fait usuellement, fournit une suite d'isomorphismes locaux $\theta _{i}: (\mathbb{K}\hspace{3pt},0)\times (\mathbb{K}^{n},0) \rightarrow (E_{i},O_{i})$ telle que le diagramme  $(D_{2})$ suivant soit commutatif:
\begin{equation*} \footnotesize
\raisebox{2cm}
{\xymatrix{ 
(\K,\o) \ar[r]^-\id \ar[d]^{\Gamma'_i} & (\K,\o) \ar[r]^-\id \ar[d]^{\Gamma'_{i-1}} & \cdots \ar[r]^-\id     & (\K,\o) \ar[r]^-\id \ar[d]^{\Gamma'_1} & (\K,\o) \ar[d]^{\Gamma'_0} \\
(\K,\o) \times (\K^n,\o) \ar[r]^-{\Pi'_i} \ar[d]^{\theta_i} & (\K,\o) \times (\K^n,\o) \ar[r]^-{\Pi'_{i-1}} \ar[d]^{\theta_{i-1}} & \cdots \ar[r]^-{\Pi'_2}     & (\K,\o) \times (\K^n,\o) \ar[r]^-{\Pi'_1} \ar[d]^{\theta_1} & (\K,\o) \times (\K^n,\o) \ar[d]^{\theta_0=\id} \\
(E_i,\o_i) \ar[r]^-{\Pi_i}            & (E_{i-1},\o_{i-1}) \ar[r]^-{\Pi_{i-1}}    & \cdots \ar[r]^-{\Pi_2} & (E_1,\o_1) \ar[r]^-{\Pi_1}            & (E_0,\o_0)                \\
(Z_i,\o_i) \ar@{^{(}->}[u]_{J_i} \ar[r]^-{\omega_i} & (Z_{i-1},\o_{i-1}) \ar@{^{(}->}[u]_{J_{i-1}} \ar[r]^-{\omega_{i-1}}      & \cdots \ar[r]^-{\omega_2}   & (Z_1,\o_1) \ar@{^{(}->}[u]_{J_1} \ar[r]^-{\omega_1}              & (Z_0,\o_0) \ar@{^{(}->}[u]_{J_0}
}}
\end{equation*}

Via ces isomorphismes on a:
\begin{itemize}
\item[-] $\forall i\geq 1$, $\Pi '_{i}(t,X)=(t,t(A_{i}+X))$ o\`u $X=(X_{1},\ldots,X_{n})$;
\item[-] $\Gamma _{i}=\theta _{i}\circ \Gamma _{i}'$, $\Gamma _{i}'(t)=(t,\varphi _{i}(t))$ o\`u $\varphi _{i}(t)=\sum _{k>i} A_{k}.t^{k-i} \in \mathcal{A}_{(\mathbb{K}^{n},0)}.$
\end{itemize}

\vspace{8pt}

On notera de plus par:
\begin{itemize}
\item[-] $(Z'_{i},0)$ le germe image de $(Z_{i},O_{i})$ par $\theta ^{-1}_{i}\circ J_{i}$;
\item[-] $I_{i,A^{i}} \subset \mathbb{K}\{t,X_{1},\ldots,X_{n}\}$ l'id\'eal d\'efinissant $(Z'_{i},0)$, o\`u $A^{i}=(A_{1},\ldots,A_{i})$ $ \in \mathbb{K}^{ni}$.
\end{itemize}

\vspace{5pt}

Pour $i \in \mathbb{N}$, on note:

- $N_{i,\varphi}=N(I_{i,A^{i}})\subset\mathbb{N}^{n+1}$ le diagramme des exposants initiaux de $I_{i,A^{i}}$ pour le syst\`eme de coordonn\'ees $(t,X)$. (Bien entendu si $O_{i} \not\in Z_{i}$ i.e. $I_{i,A^{i}}=\mathbb{K}\{t,X\}$, on pose $N_{i,\varphi}=\mathbb{N}^{n+1})$.

On fait alors la d\'efinition suivante:

\begin{defi}
Soient $(X,x)$ un germe de $\mathbb{K}$-espace analytique et un plongement $(X,x) \hookrightarrow (\mathbb{K}^n,0)$. Pour $\varphi \in \mathcal{A}_{(\mathbb{K}^{n},0)}$, on appelle:
\begin{itemize}
\item[ 1)] Suite des multiplicit\'es de Nash de $(X,x)$ le long de $\varphi$ la suite:
$$\mathcal{M}_{X,\varphi}=(m_{0,\varphi},m_{1,\varphi},\ldots ,m_{i,\varphi},\ldots )$$
\item[ 2)] Suite de Nash des fonctions de Hilbert-Samuel de $(X,x)$ le long de $\varphi$ la suite:
$$\mathcal{H}_{X,\varphi}=(H_{0,\varphi},H_{1,\varphi},\ldots ,H_{i,\varphi},\ldots )$$  
\item[ 3)] Suite de Nash des diagrammes des exposants initiaux  de $(X,x)$ le long de $\varphi$ la suite:
$$\mathcal{N}_{X,\varphi}=(N_{0,\varphi},N_{1,\varphi},\ldots ,N_{i,\varphi},\ldots )$$
\end{itemize}
\end{defi} 
Pour le fait que le premier point de la d\'efinition pr\'ec\'edente g\'en\'eralise la notion correspondante pour les hypersurfaces, nous renvoyons le lecteur \`a [H1] et [G-S-L-J]  prop. 4.3 p. 18 et commentaires.

Ces d\'efinitions appellent quelques remarques. Tout d'abord pour $\varphi \in \mathcal{A}_{(X,x)}$, les deux premi\`eres notions pr\'ec\'edentes ne d\'ependent en rien du plongement. En effet, dans ce cas pour les d\'efinir, il suffit de consid\'erer les $\Gamma_{i}$ comme des \'el\'ements de $\mathcal{A}_{(Z_{i},O_{i})}$ et de dire que les $\omega_{i}$ sont les \'eclatements de $Z_{i-1}$ de centre $O_{i-1}$. D'autre part, pour un $i$ donn\'e, les tronqu\'es jusqu'\`a l'ordre $i$, $\mathcal{M}^{i}_{X,\varphi^{i}}$, $\mathcal{H}^{i}_{X,\varphi^{i}}$, $\mathcal{N}^{i}_{X,\varphi^{i}}$, des suites pr\'ec\'edentes ne d\'ependent, que du tronqu\'e $\varphi^{i}$ \`a l'ordre $i$ de $\varphi$. C'est \`a dire du morphisme de $\mathbb{K}$-alg\`ebre locale $\varphi^{i*}: \mathbb{K}\{X_{1},\ldots,X_{n}\} \rightarrow \mathbb{K}\{t\}/(t)^{i+1}$ obtenu par passage au quotient. Finalement la motivation de la consid\'eration de tels diagrammes est la suivante. Soit $i \in \mathbb{N}^{*}$ et $\varphi \in \mathcal{A}_{(\mathbb{K}^{n},0)}$. Ecrivons:
$$\varphi(t)=\sum_{k=1}^{i} A_{k}t^{k} +t^{i}\varphi_{i}(t) \textit{ avec } A_{k} \in \mathbb{K}^{n}, \hspace{3pt} \varphi_{i}\in \mathcal{A}_{(\mathbb{K}^{n},0)}$$
Alors $\varphi \in \mathcal{A}_{(X,x)}$ si et seulement si l'arc $t \rightarrow (t,\varphi_{i}(t))$ est dans $\mathcal{A}_{(Z'_{i},0)}$ $\simeq\mathcal{A}_{(Z_{i},O_{i})}$. Aussi nous semble-t-il naturel d'\'etudier de tels objets, puisqu'ils sont les premi\`eres mesures de la singularit\'e des $(Z'_{i},0)$.

 Nous donnons deux r\'esultats dans cette direction. Le premier g\'en\'eralise en toute codimension le th\'eor\`eme 5 p. 154 de [L-J].

\begin{theo}
Soient $(X,x) \hookrightarrow (\mathbb{K}^{n},0)$ et $\varphi \in \mathcal{A}_{(\mathbb{K}^{n},O)}$.
\begin{itemize}
\item[1)] La suite $\mathcal{M}_{X,\varphi}$ est d\'ecroissante et se stabilise \`a la valeur g\'en\'erique de la multiplicit\'e de $X$, en un point $\varphi(t)$, $t\neq 0$, assez petit.
\item[2)] La suite $\mathcal{H}_{X,\varphi}$ est d\'ecroissante (i.e. $\forall i,k \in \mathbb{N}, \hspace{3pt}H_{i+1,\varphi}(k) \leq H_{i,\varphi}(k)$) et se stabilise \`a la valeur g\'en\'erique de la fonction de Hilbert-Samuel de $Z_{0}$, en un point $(t,\varphi(t))$, $t\neq 0$, assez petit.
\end{itemize}
\end{theo}
Notons que la fonction de Hilbert-Samuel de $Z_{0}$ en un point $(t,y)$ se calcule directement \`a partir de celle de $X$ en $y$. Nous essayerons au cours de la preuve de cet \'enonc\'e de pr\'eciser l'ordre auquel ces suites se stabilisent.

Les propri\'et\'es de g\'en\'ericit\'e et de semi-continuit\'e  des applications $\varphi \rightarrow \mathcal{M}_{X,\varphi}$, $\varphi \rightarrow \mathcal{H}_{X,\varphi}$, et $\varphi \rightarrow \mathcal{N}_{X,\varphi}$ sont ensuite \'etudi\'ees.

 Pour cela, soit $i \in \mathbb{N}^{*}$, d\'esignons par $\mathcal{A}_{(X,x)}^{i}$ l'ensemble des morphismes de $\mathbb{K}$-alg\`ebre locale de $\mathcal{O}_{X,x} \rightarrow  \mathbb{K}\{t\}/(t)^{i+1}$, et par $\Pi^{i}: \mathcal{A}_{(X,x)} \rightarrow \mathcal{A}_{(X,x)}^{i}$ le morphisme induit par la surjection canonique de $\mathbb{K}\{t\} \rightarrow \mathbb{K}\{t\}/(t)^{i+1}$. Tout choix de syst\`eme de coordonn\'ees sur $\mathbb{K}^{n}$ induit une identification de $\mathcal{A}_{(\mathbb{K}^{n},0)}^{i}$ avec $ \mathbb{K}^{ni}$. Ordonnons ensuite, l'ensemble des suites d'entiers (resp. l'ensemble des suites de fonctions de $\mathbb{N}$ dans $\mathbb{N}$) par l'ordre lexicographique. De m\^eme l'ensemble $\mathcal{D}(n+1)$ des parties $N$ de $\mathbb{N}^{n+1}$ satisfaisant $N+\mathbb{N}^{n+1}=N$ \'etant muni de son ordre total naturel (cf. section 2), on ordonne l'ensemble des suites d'\'el\'ements de $\mathcal{D}(n+1)$ par l'ordre lexicographique. On a alors:

\begin{theo}
Soient $(X,x)$ et $i \in \mathbb{N}^{*}$. La fonction $\varphi^{i} \rightarrow \mathcal{N}^{i}_{X,\varphi^{i}}$ est  semi-continue sup\'erieurement pour la topologie de Zariski sur $\mathcal{A}_{(\mathbb{K}^{n},0)}^{i}$. C'est \`a dire, pour toute sous-vari\'et\'e alg\'ebrique irr\'eductible $V$ de $\mathcal{A}_{(\mathbb{K}^{n},0)}^{i}$, il existe $V'$  vari\'et\'e alg\'ebrique strictement incluse dans V telle que:
\begin{itemize}
\item[1)] $\forall \varphi^{i},\hspace{3pt}\theta^{i} \in V-V'$, $\mathcal{N}^{i}_{X,\varphi^{i}}=\mathcal{N}^{i}_{X,\theta^{i}}$;
\item[2)] $\forall \varphi^{i} \in V-V'$, $\forall \theta^{i} \in V'$,  $\mathcal{N}^{i}_{X,\varphi^{i}}<\mathcal{N}^{i}_{X,\theta^{i}}$ .
\end{itemize}
\end{theo}

\begin{cor}Soient $(X,x)$ et $i \in \mathbb{N}^{*}$. Il existe une partition:
$$\overline{\Pi^{i}(\mathcal{A}_{(X,x)})}=\bigcup_{1\leq j \leq l_{i}} S_{i,j}$$
o\`u $S_{i,j}=U_{i,j}-W_{i,j}$, $U_{i,j}$ et $W_{i,j}$ sont des sous-vari\'et\'es alg\'ebriques de $\mathcal{A}^{i}_{(\mathbb{K}^{n},0)}$ et $\mathcal{M}^{i}_{X,\varphi^{i}}$, $\mathcal{H}^{i}_{X,\varphi^{i}}$, $\mathcal{N}^{i}_{X,\varphi^{i}}$ sont constantes sur $S_{i,j}$. De plus pour $j\neq j'$ les valeurs de $\mathcal{N}^{i}_{X,\varphi^{i}}$ sur $S_{i,j}$ et $S_{i,j'}$ sont distinctes.
\end{cor}

Les constructibles $\Pi^{i}(\mathcal{A}_{(X,x)})$, dont l'\'etude fut initi\'ee par J. Nash (cf. [N]), portent donc une stratication naturelle, nous esp\'erons que celle-ci pourra s'av\'erer utile \`a une meilleure compr\'ehension de leurs g\'eom\'etries. \\
Comme premi\`ere illustration de l'utilit\'e de ces stratifications, nous d\'efinissons les notions  de strates principales et de partie principale d'ordre $i$ de $\mathcal{A}_{(X,x)}$ (c.f. section 5, def. 5.1). Nous montrons comme cons\'equence \'el\'ementaire du th\'eor\`eme de semi-continuit\'e comment le  volume motivique de $\mathcal{A}_{(X,x)}$ (c.f. [D-L1])  se r\'ealise comme \guillemotleft limite \guillemotright des parties principales (c.f. th. 5.2). Nous traitons par cette technique quelques exemples.

\vspace{10pt}

Notre expos\'e est organis\'e comme suit. La section 2 rappelle quelques r\'esultats et techniques que nous utiliserons par la suite. Ceux-ci sont issus de [B-M1], [B-M2], [T1], [T2]. Nous rappelons aussi quelques notations et d\'efinitions concernant la th\'eorie de l'int\'egration motivique. Le paragraphe 3 prouve le th\'eor\`eme 1.3. Notre sch\'ema de preuve suit d'assez pr\`es  celui de J.C. Tougeron dans le cas classique [T3]. La d\'emonstration du th\'eor\`eme 1.4 n\'ecessitant quelques \'el\'ements sur la suite des multiplicités de Nash, les résultats et preuves  concernant ces notions sont d'abord donn\'es dans la section 4.  Les techniques de base proviennent du travail de E. Bierstone et P.D. Milman [B-M1]. La section 5 pr\'esente les notions de \guillemotleft strates principales et de partie principale d'ordre  $i$ \guillemotright  de $\mathcal{A}_{(X,x)}$ et leur rapport avec le volume motivique tel que d\'efini dans [D-L1]. Quelques exemples viennent illustrer cette technique. Le th\'eor\`eme 1.4 est d\'emontr\'e au paragraphe 6. L'existence de d\'esingularisation peut servir ici de substitut aux \'el\'ements provenant de la section 4. On perd alors en \'el\'ementarit\'e mais aussi beaucoup d'informations concernant le diam\`etre des ouverts dont il est question  dans l'\'enonc\'e 1.4.

\section{Rappels et notations}
Nous ferons un usage r\'ep\'et\'e du r\'esultat suivant. Nous l'appelerons par commodit\'e th\'eor\`eme des fonctions implicites de Tougeron (cf. [T2]).

\begin{theo}[J.C. Tougeron]
Soient $x=(x_{1},\ldots,x_{n})$, $y=(y_{1},\ldots,y_{p})$, $f=(f_{1},\ldots,f_{q})\in \mathbb{K}\{x,y\}^{q}$ telle que $f(0,0)=0$. Soit $I$ l'id\'eal engendr\'e dans $ \mathbb{K}\{x\}$ par les mineurs d'ordre q de la matrice jacobienne $f'_{y}(x,0)=((\partial f_{i}/\partial y_{j}) (x,0))$, et soit $I'$ un id\'eal propre de $\mathbb{K}\{x\}$. Si $f(x,0)\in\oplus_{q} I^{2}.I'$ , il existe $y(x)\in\oplus_{p} I.I'$ tel que $f(x,y(x))=0$. 
\end{theo}
Le résultat est aussi valable en  formel et en $C^{\infty}$.

Soit $\alpha=(\alpha_{1},\ldots,\alpha_{n})\in\mathbb{N}^{n}$, on  notera: $\vert\alpha\vert=\alpha_{1}+\cdots+\alpha_{n}$. $\mathbb{N}^{n}$ est totalement ordonn\'e par l'ordre lexicographique sur les $n+1$ upplets  $(\vert\alpha\vert,\alpha_{1},\ldots,\alpha_{n})$. Si $A$ est un anneau commutatif unitaire intègre, et $f$ un élément non nul de $A\lbrack\lbrack X \rbrack\rbrack=A\lbrack\lbrack X_{1},\ldots,X_{n}\rbrack\rbrack$, nous noterons par $\nu (f)$ son exposant initial. C'est \`a dire si : 
$$ f=\sum_{\alpha\in\mathbb{N}^n}a_{\alpha}.X^{\alpha}, \hspace{4pt} \nu (f)=Min\{\alpha\in\mathbb{N}^n \vert a_{\alpha}\not=0\}.$$ 
Pour un id\'eal, $I \subset A[[X]]$, on notera $N_{I}$ le diagramme des exposants initiaux de $I$, i.e.: 
$$N_{I}=\{\alpha \in \mathbb{N}^{n} \vert \exists g \in I \textit{ tel que }  \nu(g)=\alpha\}.$$
On a : $N_{I} + \mathbb{N}^{n} =N_I$. Si $N \subset \mathbb{N}^{n}$ satisfait $N + \mathbb{N}^{n} =N$, le lemme de Dickson assure l'existence d'une partie finie unique de $\mathbb{N}^{n}$, $\{\alpha^{1},\ldots,\alpha^{p}\}$, telle que : 
$$N=\cup_{1\leq i\leq p}(\alpha^{i} +  \mathbb{N}^{n}) \textit { et } \alpha ^{j} \notin \cup_{i\neq j} (\alpha^{i} +  \mathbb{N}^{n})$$
Les $\alpha^{i}$ sont dits les sommets de $N$. L'ensemble $\mathcal{D}(n)$ des parties  de $\mathbb{N}^{n}$ stables par translations (i.e. satisfaisant $N +\mathbb{N}^{n}=N$) est totalement ordonn\'e comme suit. Soit $N_{1},N_{2}\in \mathcal{D}(n)$. Pour chaque $i=1,2$, soient $\beta ^{k}_{i}$,$k=1,\ldots,t_{i}$ les sommets de $N_{i}$ index\'es dans l'ordre croissant. Apr\`es avoir \'eventuellement permut\'e $N_{1}$ et $N_{2}$, il existe $t\in \mathbb{N}$ tel que: $\beta^{k}_{1}=\beta^{k}_{2}$, $1\leq k\leq t$ et (1) $t_{1}=t=t_{2}$,(2) $t_{1}>t=t_{2}$ ou bien (3) $t_{1},t_{2}>t$ et $\beta_{t+1}^{1}<\beta_{t+1}^{2}$. Dans le cas (1), $N_{1}=N_{2}$. Dans les cas (2) et (3), $N_{1}<N_{2}$. Il revient au m\^eme de dire que la suite $(\beta ^1_1,\ldots,\beta ^{t_1}_1,\infty,\ldots)$ est strictement plus petite que la suite $(\beta ^1_2,\ldots,\beta ^{t_2}_2,\infty,\ldots)$ pour l'ordre lexicographique, avec la convention $\beta <\infty $ pour tout $\beta \in \mathbb{N}^n$.\\

Si $N$ est un \'el\'ement de $\mathcal{D}(n)$, on notera $H_{N}$ la fonction de $\mathbb{N}$ dans $\mathbb{N}$ d\'efinie par: $$H_{N}(k)=card\{\alpha \in \mathbb{N}^{n}\vert \alpha \not \in N\textrm{ et } \vert \alpha \vert \leq k\}.$$

Soit I un id\'eal de $\mathbb{K}[[X]]$ et $g_{1},\ldots,g_{p}$  des \'el\'ements de $I$ tels que  $\nu (g_{1}),\ldots$ $,\nu (g_{p})$ sont les sommets de $N(I)$. Nous disons que $g_{1},\ldots,g_{p}$ est une base standard de $I$. Il existe par ailleurs une unique famille $g_{1},\ldots,g_{p}$ telle que : 
$$InMon(g_{i})=X^{\beta^{i}} $$ 
$$Supp(X^{\beta^{i}}-g_{i}) \subset \mathbb{N}^{n}-N(I), \hspace{5pt} \forall 1\leq i\leq p,$$ 
o\`u $\beta^{i},\ldots,\beta^{p}$ sont les sommets de $N(I)$, $Supp$ d\'esigne le support d'un \'el\'ement de $\mathbb{K}[[X]]$, $InMon(f)=a_{\nu (f)}X^{\nu (f)}$.
Une telle famille est appel\'ee la base standard distingu\'ee de $I$. Si:
$$H_{I}(k) = dim_{\mathbb{K}} \frac {\mathbb{K}[[X]]} {I+\mathcal{M}^{k+1}} $$
d\'esigne la fonction de Hilbert-Samuel de $I$, on a : 
$$H_{I}(k)=card\{\alpha \in  \mathbb{N}^{n}-N(I)\vert \hspace{3pt} \vert \alpha  \vert \leq k\}=H_{N(I)}(k).$$ 
Pour plus d'informations, nous renvoyons le lecteur \`a [B-M2]. Nous ferons par ailleurs un usage libre de la notion de transform\'e strict d'un id\'eal de  $\mathbb{K}[[X]]$. Rappelons simplement que si $g_{1},\ldots,g_{p}$  est une base standard de $I$, alors leurs transform\'ees strictes  $g'_{1},\ldots,g'_{p}$,  par l'\'eclatement de centre l'origine engendrent le transform\'e strict $I'$ de $I$.\\

Dans le texte de l'article, nous consid\`erons des anneaux de la forme $A[[t,X]]$, $X=(X_1,\ldots,X_n)$. Lorsque nous parlons de diagramme des exposants initiaux  $\subset \mathbb{N}^{n+1}$ d'un id\'eal $I$ de $A[[t,X]]$, nous consid\'erons toujours des multi-indices $(\alpha _0,\alpha _1,\ldots,\alpha _n)$ o\`u le premier indice $\alpha _0$  correspond toujours \`a la variable $t$.

 Si X est un espace analytique et $\varphi\in \mathcal{A}_{X}$, nous noterons par $\mathcal{B}_{X,i}(\varphi)$ la boule : $$\mathcal{B}_{X,i}(\varphi)=\{\varphi '\in \mathcal{A}_ {X}/\hspace{3pt} \mathcal{D}_ {X}(\varphi ,\varphi ')<e^{-i}\}=\{\varphi '\in \mathcal{A}_ {X}/\hspace{3pt} ord(\varphi -\varphi ')>i\}.$$ 
Lorsque le contexte le permettra, nous omettrons l'indice $X$.\\

Rappelons que l'on d\'esigne par $K_0(\mathcal{V}_{\mathbb{K}})$ l'anneau de Grothendieck des $\mathbb{K}$-vari\'et\'es (i.e. les sch\'emas r\'eduits s\'epar\'es de type fini sur $\mathbb{K}$). Par d\'efinition $K_0(\mathcal{V}_{\mathbb{K}})$ est l'anneau engendr\'e par les symboles $[V]$, pour une $\mathbb{K}$-vari\'et\'e $V$, modulo les relations:\\
- $[V]=[V']$, si $V$ et $V'$ sont $\mathbb{K}$-isomorphes\\
- $[V]=[V']+[V-V']$, si $V'$ est ferm\'e dans $[V]$\\
- $[V]\times [V']=[V\times _{\mathbb{K}}V']$.\\
La classe de la droite affine $\mathbb{A}^1_{\mathbb{K}}$ sera not\'ee $\mathbb{L}$. Suivant les notations de Looijenga [Lo], $\mathcal{M}_{\mathbb{K}}$ d\'esignera l'anneau $\K_0(\mathcal{V}_{\mathbb{K}})[\mathbb{L}^{-1}]$ (ce qui correspond \`a la notation $\mathcal{M}_{loc}$ de [D-L1]). Pour $m\in \mathbb{Z}$, on note $F^m\mathcal{M}_{\mathbb{K}}$ le sous-groupe engendr\'e par les symboles $[S]\mathbb{L}^{i}$ o\`u $i-dimS\geq m$. On notera $\widehat{\mathcal{M}}_{\mathbb{K}}$ le compl\'et\'e de $\mathcal{M}_{\mathbb{K}}$ par rapport \`a cette filtration et $\overline{\mathcal{M}}_{\mathbb{K}}$ l'image de  $\mathcal{M}_{\mathbb{K}}$ par l'application naturelle de $\mathcal{M}_{\mathbb{K}}$ dans $\widehat{\mathcal{M}}_{\mathbb{K}}$. Concernant les mesures motiviques de [D-L1], nous adoptons la convention de [Lo], i.e. les mesures de Denef-Loeser sont multipli\'ees par $\mathbb{L}^d$, $d$ d\'esignant la dimension de la vari\'et\'e ou du germe que l'on consid\`ere.

\section{In\'egalit\'es de Lojasiewicz fibr\'ees au niveau de l'espace des arcs}

Nous commençons d'abord par constater qu'il suffit d'obtenir le th\'eor\`eme 1.3, lorsque $X$ et $Y$ sont respectivement des ouverts $\Omega _{n}$ et $\Omega _{p}$ de $\mathbb{K}^n$ et $\mathbb{K}^p$.

\vspace{10pt}

 En effet, \'etant donn\'e un point $(x,y)\in X\times Y$, il suffit d'\'etablir l'existence de voisinages $V_{x}$ et $V_{y}$ et d'un $\alpha \geq 0$ tels que la propri\'et\'e de 1.3 soit valide pour tout $(\phi _{0}, \varphi) \in \mathcal{A}_ {X,V_{x}} \times \mathcal{A}_ {Y,V_{y}}$. On peut par cons\'equent supposer que $X$ et $Y$ sont des mod\`eles locaux, $X$ d\'efini par $f_{1},\ldots,f_{k}$, analytiques sur un ouvert $\Omega _{n}$, et que le morphisme $F$ est induit par des fonctions analytiques  $F_{1},\ldots,F_{p}$ sur $\Omega _{n}$. 
Consid\'erons alors le morphisme $G$: 

$$G: \Omega _{n} \rightarrow \Omega _{p} \times \mathbb{K}^{k}$$

$$ \hspace{45pt} x \rightarrow (F_{1}(x),\ldots,F_{p}(x),f_{1}(x),\ldots,f_{k}(x)).$$

Si le th\'eor\`eme 1.3  est valide pour le morphisme $\mathcal{G}$, il l'est pour $\mathcal{F}$. En effet,  si $\varphi \in \mathcal{A}_{Y} \subset \mathcal{A}_{\Omega _{p}}$, consid\'erons l'arc $\theta=(\varphi,0) \in \mathcal{A}_{\Omega_{p}\times \mathbb{K}^{k}}$. On a pour tout $\phi \in \mathcal{A}_{X} \subset \mathcal{A}_{\Omega _{n}}$: 
$$\mathcal{D}_{Y}(\mathcal{F}(\phi),\varphi)=\mathcal{D}_{\Omega _{n}\times \mathbb {K}^{k}}(\mathcal{G}(\phi),\theta)$$
$$\mathcal{D}_{X}(\phi,\mathcal{F}^{-1}(\varphi))=\mathcal{D}_{\Omega _{n}} (\phi,\mathcal{G}^{-1}(\theta)).$$
Nous \'etablissons maintenant $1.3$ dans la situation cit\'ee plus haut. N\^otre sch\'ema de d\'emonstration suit, tout en l'adaptant \`a notre cas,  celui de J.C. Tou\-ge\-ron [T3]. Pour cela d\'esignons par $\mathcal{O}$ le faisceau des germes de fonctions $\mathbb{K}$-analytiques sur $\Omega _{n} \times \Omega _{p}$ et $\mathcal{I} \subset \mathcal{O}$ un sous-faisceau coh\'erent d'id\'eaux. Nous utiliserons les notations suivantes :
\begin{itemize}
\item[-] $V(\mathcal{I})$ l'espace analytique de z\'eros de $\mathcal{I}$, i.e. $V(\mathcal{I})=(\vert X \vert , (\mathcal{O}/ \mathcal{I})\vert  \vert X \vert )$ o\`u  $\vert X\vert=Support(\mathcal{O}/\mathcal{I})$; 
\item[-] $\mathcal{V}(\mathcal{I})=\mathcal{A}_{V(\mathcal{I})}$;
\item[-] Soit $\varphi \in \mathcal{A}_{\Omega_{p}}$, $\mathcal{V}(\varphi ^{*} \mathcal{I})=\{\phi \in \mathcal{A}_{\Omega_{n}}\vert (\phi,\varphi) \in \mathcal{V}(\mathcal{I})\}$;
\item[-] Si $(x_{0},y_{0}) \in \Omega_{n} \times \Omega_{p}$, on d\'esignera par $\mathcal{I}_{(x_{0},y_{0})}$ la fibre de $\mathcal{I}$ en $(x_{0},y_{0})$ et par $f_{1}(x,y),\ldots,f_{k_{0}}(x,y)$  un syst\`eme de g\'en\'erateurs de $\mathcal{I}_{(x_{0},y_{0})}$;
\item[-] Si $\theta=(\theta_{1},\ldots,\theta_{s})\in \mathcal{O}_{1}^{s}$, on pose: $\vert \theta \vert =\vert \theta_{1} \vert+ \ldots +\vert \theta_{s} \vert$, o\`u $\vert \theta_{j}\vert=e^{-ord(\theta_{j})}$.
\end{itemize}

\vspace{10pt}

Le th\'eor\`eme 1.3 d\'ecoule alors du r\'esultat suivant, appliqu\'e au sous-faisceau de $\mathcal{O}$ engendr\'e par $F_{1}(x)-y_{1},\ldots, F_{p}(x)-y_{p}$, o\`u $F_{1},\ldots, F_{p}$ sont les composantes du morphisme $F:\Omega_{n} \rightarrow \Omega_{p}$.

\begin{theo}
Soient $\mathcal{I} \subset \mathcal{O}$ et $K,K'$ deux compacts comme pr\'ec\'edemment. Il existe alors $\alpha \geq 0$ tel que : 

$\forall \varphi \in \mathcal{A}_{(\Omega_{p},K')}$, $\forall \phi _{0} \in \mathcal{A}_{(\Omega_{n},K)}$, $\exists C_{(\phi_{0},\varphi)}>0$ et un voisinage $\mathcal{V}_{\phi_{0}}$ de $\phi_{0}$ satisfaisant:
$$\forall \phi \in \mathcal{V}_{\phi_{0}}, \sum_{i=1}^{k_{0}} \vert f_{i}(\phi,\varphi)\vert \geq  C_{(\phi_{0},\varphi)} \mathcal{D}(\phi,\mathcal{V}(\varphi^{*} \mathcal{I}))^{\alpha}$$
o\`u $f_{1}(x,y),\ldots,f_{k_{0}}(x,y)$ engendrent $\mathcal{I}_{(\phi_{0}(0),\varphi(0))}$.
\end{theo}
\textit{Preuve:}

Pour prouver le r\'esultat, il suffit de voir que pour un point $(x_{0},y_{0}) \in \Omega_{n} \times \Omega_{p}$, il existe des voisinages $V_{x_{0}}$, $W_{y_{0}}$ de $x_{0}$ et $y_{0}$, et un $\alpha \geq 0$ tels que la propri\'et\'e de 3.1 soit satisfaite pour $(\phi_{0}, \varphi) \in \mathcal {A}_{V_{x_{0}}} \times \mathcal {A}_{W_{y_{0}}}$. Si $(x_{0},y_{0}) \notin V(\mathcal{I})$, ceci est trivialement satisfait avec $\alpha = 0$. Nous supposerons donc  $(x_{0},y_{0}) = (0,0) \in V(\mathcal{I})$. Pour prouver l'assertion pr\'ec\'edente, il suffit de voir que : 

\vspace{10pt}

(3.2) il existe un voisinage $V_{0}$ de 0 dans $\mathbb{K}^{p}$ et un $\alpha \geq 0$ satisfaisant :
$\forall \varphi \in \mathcal{A}_{V_{0}}$ tel que $(0,\varphi) \in \mathcal{V}(\mathcal{I})$, $\exists C_{\varphi}>0$ et un voisinage $\mathcal{V}_{0}$ de $0$ dans $\mathcal{A}_{(\mathbb{K}^{n},0)}$ tel que: 
$$\forall \phi \in \mathcal{V}_{0}, \sum_{i=1}^{k_{0}} \vert f_{i}(\phi,\varphi) \vert \geq C_{\varphi}\mathcal{D}(\phi,\mathcal{V}(\varphi^{*}\mathcal{I}))^{\alpha}.$$ 
En effet supposons (3.2)  prouv\'ee pour tout faisceau coh\'erent d'id\'eaux. Partant d'un faisceau coh\'erent d'id\'eaux $\mathcal{I}$ engendr\'e par $f_{1}(x,y),\ldots,f_{k_{0}}(x,y)$ sur un voisinage de l'origine dans $\mathbb{K}^{n} \times \mathbb{K}^{p}$, on consid\`ere le faisceau d'id\'eaux $\mathcal{J}$ engendr\'e sur un voisinage de l'origine dans $\mathbb{K}^{n} \times \mathbb{K}^{n+p}$ par :  
$$G_{1}(x,z,y)=f_{1}(x+z,y),\ldots,G_{k_{0}}(x,z,y) =f_{k_{0}}(x+z,y).$$ 
L'application de (3.2) à $\mathcal{J}$ fournit un voisinage $U_{0}=V_{0} \times W_{0}$ de 0 dans $\mathbb{K}^{n+p} = \mathbb{K}^{n} \times \mathbb{K}^{p}$ et un $\alpha \geq 0$ satisfaisant: 

\vspace{10pt}

(3.3) $\forall (\phi_{0}, \varphi) \in \mathcal{A}_{U_{0}}= \mathcal{A}_{V_{0}} \times \mathcal{A}_{W_{0}}$ tel que $(0,\phi_{0}, \varphi) \in \mathcal{V}(\mathcal{J})$, 

$\exists \hspace{5pt} C_{(\phi_{0}, \varphi)}$  et un voisinage $\mathcal{V}_{0}$ de 0 dans $\mathcal{A}_{(\mathbb{K}^n,0)}$ tels que : 
$$\forall \hspace{5pt} \theta \in \mathcal{V}_{0}, \sum_{i=1}^{k_{0}} \vert G_{i}(\theta,\phi_{0},\varphi) \vert \geq  \mathcal{C}_{\phi_{0}, \varphi} \mathcal{D} (\theta,\mathcal{V}((\phi_{0}, \varphi )^{*} \mathcal{J})^{\alpha}.$$
 Soient alors $\varphi \in \mathcal{A}_{V_{0}}$ et $\phi_{0} \in \mathcal{A}_{W_{0}}$. Si $(\phi_{0}, \varphi) \notin \mathcal{V}(\mathcal{I})$, posons :
 $$s=Min(ord f_{i}(\phi_{0}, \varphi)) \hspace{5pt} \mathcal {C}_{(\phi_{0}, \varphi)}=e^{-s} \hspace{5pt} \textrm{et }  \mathcal{V}_{\phi_{0}}=\mathcal{B}(\phi_{0},e^{-s}).$$
 On a alors:  
$$ \forall \hspace{5pt} \phi \in \mathcal{V}_{\phi_{0}}, \sum_{j=1}^{k_{0}} \vert f_{j}(\phi,\varphi) \vert \geq \mathcal{C}_{(\phi_{0}, \varphi)} \geq \mathcal{C}_{(\phi_{0}, \varphi)} \mathcal{D}(\phi,\mathcal{V}(\varphi^{*}\mathcal{I}))^{\alpha}.$$ 
Si maintenant $(\phi_{0}, \varphi) \in \mathcal{V}(\mathcal{I})$, alors $(0,\phi_{0}, \varphi) \in \mathcal{V}(\mathcal{J})$. D'o\`u l'existence d'un voisinage $\mathcal{V}_{0}$ de $0$  dans $\mathcal{A}_{(\mathbb{K}^n,0)}$ et d'une constante $\mathcal{C}_{(\phi_{0}, \varphi)}$ comme dans (3.3). 
Posons: $\mathcal{V}_{\phi_{0}}= \phi_{0} + \mathcal{V}_{0}$. Puisque:    $$\sum_{i=1}^{k_{0}}\vert f_{i}(\phi,\varphi)\vert = \sum_{i=1}^{k_{0}}\vert G_{i}(\phi -\phi_{0},\phi_{0},\varphi)\vert ;$$
$$\mathcal{D}(\phi,\mathcal{V}(\varphi^{*}\mathcal{I}))= \mathcal{D}(\phi - \phi_{0},\mathcal{V}((\phi_{0},\varphi)^{*}\mathcal{J}).$$
On a:
$$\forall \phi \in \mathcal{V}_{\phi_{0}}, \sum_{i=1}^{k_{0}}\vert f_{i}(\phi,\varphi)\vert \geq C_{(\phi_{0},\varphi)}\mathcal{D}(\phi,\mathcal{V}(\varphi^{*}\mathcal{I}))^{\alpha}.$$
Ce qui est l'assertion d\'esir\'ee. Pour prouver (3.2) nous laissons le soin au lecteur de v\'erifier que celle-ci est \'equivalente \`a:  

\vspace{10pt}

(3.4) $\exists \hspace{5pt} V_{0}$ voisinage de 0 dans $\mathbb{K}^{p}$ et un $\alpha \geq 0$ satisfaisant :
 $\forall \hspace{5pt} \varphi \in \mathcal{A}_{V_{0}}$ tel que $(0,\varphi) \in \mathcal{V}(\mathcal{I})$, il existe $a_{\varphi},k_{\varphi} \in \mathbb{N}$ tels que l'implication suivante soit v\'erifi\'ee :  
$$ \forall \hspace{5pt} i \in \mathbb{N}, \forall \hspace{5pt} \phi \in \mathcal{A}_{(\mathbb{K}^{n},0)}, (ord(\phi) > k_{\varphi}\textit{ et } Min_{1 \leq j \leq k_{0}} ord(f_{j}(\phi,\varphi))> a_{\varphi}+ \alpha i)$$  

$$\Longrightarrow  (\exists \hspace{5pt} \phi ' \in \mathcal{A}_{(\mathbb{K}^{n},0)}, \phi ' \in \mathcal{V}(\varphi^{*} \mathcal{I}) \textit{ et } ord(\phi -\phi ') >i).$$ 

\vspace{10pt}

Nous allons \`a pr\'esent prouver (3.4) par r\'ecurrence descendante sur la hauteur de $\mathcal{I}_{(0,0)}$.

 Si $haut(\mathcal{I}_{(0,0)})=n+p$, quitte \`a nous restreindre \`a un voisinage assez petit de $(0,0)$ dans $\mathbb{K}^{n} \times \mathbb{K}^{p}$, on a $V(\mathcal{I})=(0,0)$  et $\exists m \in \mathbb{N}^{*}$ tel que:  $(x_{1},\ldots,x_{n},y_{1},\ldots, y_{p})^{m} \subset \mathcal{I}_{(0,0)}$. Par suite (3.4) est satisfaite avec $\alpha=m$, $a_{0}=k_{0}=0$. 

Soit maintenant $\mathcal{I}$ tel que $haut(\mathcal{I}_{(0,0)})=r<n+p$. Il existe un entier $m$ et $\mathcal{P}^{1}, \ldots,\mathcal{P}^{s}$, des sous-faisceaux  coh\'erents d'id\'eaux de $\mathcal{O}$ dans un voisinage de $(0,0)$ dans $\mathbb{K}^{n} \times \mathbb{K}^{p}$ tels que :       $$\forall j, \hspace{5pt} 1\leq j \leq s,\hspace{5pt} \mathcal{P}_{(0,0)}^{j} \textrm{ est premier de hauteur } \geq r;$$
$$\exists m \in \mathbb{N} \vert (\bigcap _{1\leq j \leq s} \mathcal{P}^{j})^{m} \subset \mathcal{I} \subset  \bigcap _{1\leq j \leq s} \mathcal{P}^{j}.$$
Supposons que pour tout $j$, $1\leq j \leq s$, (3.4) soit satisfaite pour $\mathcal{P}^{j}$ avec une constante $\alpha_{j}\geq 0$. Il est alors ais\'e de v\'erifier que (3.4) est satisfaite pour $\mathcal{I}$ avec $\alpha=m(\sum_{j=1}^{s}\alpha_{j})$. En effet, soit $(0,\varphi) \in \mathcal{V}(\mathcal{I})$. D\'esignons par $\Delta_{\varphi}$ l'ensemble (n\'ecessairement non vide):
$$\Delta_{\varphi}=\{j \in \{1,\ldots,s\}\vert (0,\varphi) \in \mathcal{V}(\mathcal{P}^{j})\}.$$
Consid\'erons aussi $\overline{\Delta}_{\varphi}=\{1,\ldots,s\}-\Delta_{\varphi}$. On note alors:
\begin{itemize}
\item[-] $a_{\varphi}^{j}=ord((0,\varphi)^{*}\mathcal{P}^{j})$, si $j \in \overline{\Delta}_{\varphi}$;
\item[-] $k_{\varphi}^{j}, a_{\varphi}^{j}$ les entiers associ\'es \`a $\mathcal{P}^{j}$ par (3.4), si $j \in \Delta_{\varphi}$.
\end{itemize}
Soient alors: $k_{\varphi}=Max( Max_{j \in \overline{\Delta}_\varphi }a_{\varphi }^j, Max_{j \in \Delta_{\varphi}} k_{\varphi}^{j})$ et $a_{\varphi}=m(\sum _{j=1}^{s}a_{\varphi}^{j})$. Si $\phi$ est tel que $ord(\phi)>k_{\varphi}$ et $Min_{1 \leq l \leq k_{0}}ord(f_{l}(\phi,\varphi))>a_{\varphi}+\alpha i$, il existe $j \in \Delta_{\varphi}$ tel que $ord((\phi,\varphi)^{*}\mathcal{P}^{j})>a_{\varphi}^{j}+\alpha_{j}i$. D'o\`u l'assertion d\'esir\'ee pour $\mathcal{I}$.

\vspace{10pt}

Par suite, nous supposerons que $\mathcal{I}$ est tel que $\mathcal{I}_{(0,0)}$ est premier de hauteur $r < n+p$. Soit maintenant $J$ l'id\'eal des \'el\'ements de $\mathbb{K}\{y_{1},\ldots,y_{p}\}$ nuls sur le germe en $(0,0)$ de $V(\mathcal{I})\cap (0\times \mathbb{K}^{p})$. On peut supposer $J \subset \mathcal{I}_{(0,0)}$. En effet, si ce n'\'etait pas le cas, consid\'erons $h \in J$, $h \not\in \mathcal{I}_{(0,0)}$. Soit alors $\mathcal{J}$ le sous-faisceau coh\'erent d'id\'eaux  de $\mathcal{O}$ engendr\'e sur un voisinage de $(0,0)$ par $h$ et $\mathcal{I}_{(0,0)}$. On a alors $haut(\mathcal{J}_{(0,0)})>r$.

Donc par hypoth\`ese de r\'ecurrence, il existe un voisinage $V_{0}$ de $0$ dans $\mathbb{K}^{p}$ tel que (3.4) soit v\'erifi\'ee pour $\mathcal{J}$. On peut toujours supposer que $h$ est nul sur $V(\mathcal{I})\cap (0\times V_{0})$. Mais alors pour $\varphi \in \mathcal{A}_{V_{0}}$:
\begin{itemize}
\item[-] $(0,\varphi) \in \mathcal{V}(\mathcal{I})  \Longleftrightarrow (0,\varphi) \in \mathcal{V}(\mathcal{J})$;
\item[-] $\mathcal{V}(\varphi^{*}\mathcal{I})=\mathcal{V}(\varphi^{*}\mathcal{J})$, pour $\varphi$ tel que $(0,\varphi) \in \mathcal{V}(\mathcal{I})$;
\item[-] $\forall \phi \in \mathcal{A}_{(\mathbb{K}^{n},0)}$, $\forall \varphi \in \mathcal{A}_{V_{0}}$ tel que $(0,\varphi) \in \mathcal{V}(\mathcal{I}),$
\end{itemize}
$$Min_{1\leq j \leq k_{0}}(ord(f_{j}(\phi,\varphi))=Min(ord(h(\varphi)),ord(f_{j}(\phi,\varphi)),1\leq j \leq k_{0}).$$
Il en r\'esulte que (3.4) est aussi v\'erifi\'ee pour $\mathcal{I}$ avec le m\^eme $\alpha \geq 0$ et le m\^eme $V_{0}$. Nous supposerons donc que $J \subset \mathcal{I}_{(0,0)}$ et que cette inclusion est stricte (sans quoi (3.4) est triviale). On a alors le lemme suivant:

\begin{lem}[ \textrm{[T3]}, lemme 2.4] Supposons $J \varsubsetneq \mathcal{I}_{(0,0)}$. Soit $r-l$ la hauteur de l'id\'eal premier $J$. Apr\`es une \'eventuelle permutation sur les coordonn\'ees $x_{1},\ldots,x_{n}$, il existe $u_{1}(x,y),\ldots,u_{l}(x,y) \in \mathcal{I}_{(0,0)}$ tels que:
$$u(x,y)=det(\frac{\partial u_{k}(x,y)}{\partial x_{j}})_{1 \leq k,j \leq l} \not\in \mathcal{I}_{(0,0)}.$$
\end{lem}
Puisque $haut(J)=r-l$ et $J$ est premier, il existe $u_{l+1}(y),\ldots,u_{r}(y) \in J$
tels qu' apr\`es une \'eventuelle permutation sur les coordonn\'ees $y_{1},\ldots,y_{p}$:
$$v(y)=det(\frac{\partial u_{l+k}(y)}{\partial y_{j}})_{1 \leq k,j \leq r-l} \not\in J.$$
Par suite:
$$v(y)u(x,y)=det(\frac{\partial u_{k}(x,y)}{\partial z_{j}})_{1\leq k,j \leq r} \not\in \mathcal{I}_{(0,0)},$$
o\`u $z_{j}=x_{j}, 1 \leq j \leq l$ et $z_{j}=y_{j-l}$,$\hspace{5pt}l+1 \leq j \leq r$.
Par le crit\`ere jacobien des points simples, le localis\'e $(\mathcal{O}_{(0,0)})_{\mathcal{I}_{(0,0)}}$ est r\'egulier de dimension $r$ et son maximal est engendr\'e par $u_{1},\ldots,u_{r}$. Par cons\'equent, il existe $w(x,y) \in \mathcal{O}_{(0,0)}) \setminus \mathcal{I}_{(0,0)}$ tel que : $w.\mathcal{I}_{(0,0)} \subset (u_{1},\ldots,u_{r}) \mathcal{O}_{(0,0)}$.

Soit $\Delta =u(x,y)v(y)w(x,y)$. Consid\'erons maintenant $\mathcal{J}$ le sous-faisceau coh\'e\-rent d'id\'eaux de $\mathcal{O}$ engendr\'e par $\mathcal{I}$ et $\Delta$ sur un voisinage de $0$ dans $\mathbb{K}^{n}\times \mathbb{K}^{p}$, suffisamment petit pour que : 
$$ \Delta . \mathcal{I} \subset (u_{1},\ldots,u_{r}) .\mathcal{O}.$$
Puisque $haut(\mathcal{J}_{(0,0}) \geq r+1$, l'hypoth\`ese de r\'ecurrence nous assure l'existence d'un voisinage $V_{0}$ de $0$ dans $\mathbb{K}^{p}$ et d'un $\alpha \geq 0$ tel que (3.4) soit satisfaite pour $\mathcal{J}$. Nous allons montrer que (3.4) est valide pour $\mathcal{I}$ avec le m\^eme voisinage $V_{0}$ et avec $\beta=Max(2\alpha,1)$.

\vspace{8pt}

Soit $\varphi \in \mathcal{A}_{V_{0}}$ tel que $(0,\varphi) \in \mathcal{V}(\mathcal{I})$. Deux cas peuvent se pr\'esenter : 

\vspace{8pt}

$\underline {1^{er} cas:} \hspace{10pt}$ $(0,\varphi) \notin \mathcal{V}(\mathcal{J})$. 

Soit alors $k_{\varphi}=ord( \Delta(0,\varphi))<+\infty$. Posons $a_{\varphi}=2k_{\varphi}$. Consid\'erons alors, $\phi \in \mathcal{A}_{(\mathbb{K}^{n},0)}$ tel que $ord(\phi)>k_{\varphi}$ et $Min_{1 \leq j \leq k_{0}} ord(f_{j}(\phi,\varphi))>a_{\varphi}+i$.
On a: 
\vspace{10pt}
\begin{itemize}
\item[-] $ord(\Delta(\phi,\varphi))= ord(\Delta(0,\varphi))=k_{\varphi}$, en particulier $ord( u(\phi,\varphi)) \leq k_{\varphi}$; 
\item[-] $Min_{1 \leq j \leq r} ord(u_{j}(\phi,\varphi)) \geq Min_{1 \leq j \leq k_{0}} (ord(f_{j}(\phi,\varphi)) > 2k_{\varphi}+ i$.   
\end{itemize}
\vspace{10pt} 
Consid\'erons alors les $l$ \'el\'ements $U_{1}(t,z_{1},\ldots,z_{l}),\ldots,U_{l}(t,z_{1},\ldots,z_{l})$
de $\mathbb{K}\{t,z_{1},\-\ldots,\-z_{l}\}$ suivants:
$$U_{j}(t,z_{1},\ldots,z_{l})=u_{j}(\phi_{1}(t)+z_{1},\ldots,\phi_{l}(t)+z_{l},\phi_{l+1}(t),\ldots,\phi_{n}(t),\varphi(t)), 1 \leq j \leq l.$$
$ \hspace{15pt}U_{j}(0,0)=0,1\leq j\leq l.$

\vspace{8pt}

Alors:
$$det(\frac{\partial U_{k} }{\partial z_{j}}(t,0))_{1 \leq k,j\leq l}=u(\phi(t),\varphi(t)).$$
Par suite:             
$$ U_{j}(t,0)=u_{j}(\phi(t),\varphi(t)) \in (det(\frac{\partial U_{k} }{\partial z_{j}}(t,0))^{2}.t^{i+1+2(k_{\varphi}-valu(\phi,\varphi))}).$$
Par cons\'equent le th\'eor\`eme des fonctions implicites de Tougeron donne l'existence de $\theta_{1},\ldots,\theta_{l} \in \mathbb{K}\{t\}$ tels que : 

$$ 
\theta_{j}(t) \in det(\frac{\partial U_{k} }{\partial z_{j}}(t,0)).(t)^{i+1+2(k_{\varphi}-valu(\phi,\varphi))}
\subset (t)^{i+1+k_{\varphi}+k_{\varphi}-valu(\phi,\varphi)} \subset (t)^{i+1+k_{\varphi}}$$
$$U_{j}(t,\theta_{1}(t),\ldots,\theta_{l}(t))=0, 1\leq j\leq l.
$$

Posons alors $\phi'=(\phi_{1}+\theta_{1},\ldots,\phi_{l}+\theta_{l},\phi_{l+1},\ldots,\phi_{n})$.

On a: 
$$ord(\phi' - \phi)>i+k_{\varphi}$$
$$u_{j}(\phi', \varphi)=0, \hspace{6pt} 1\leq j\leq r \textrm{ (pour }l+1 \leq j \leq r, u_j(\phi' ,\varphi )=u_j(\varphi )=0).$$
Il en d\'ecoule que $ord (\Delta(\phi', \varphi))=ord (\Delta(\phi, \varphi))=k_{\varphi}$, donc en particulier $\Delta(\phi', \varphi) \neq 0$. 

Mais sur $V_{0} \subset \mathbb{K}^{n} \times \mathbb{K}^{p}$, on a $\Delta.\mathcal{I} \subset (u_{1},\ldots,u_{r}).\mathcal{O}$.
Donc $(\phi',\varphi) \in \mathcal{V}( \mathcal{I})$ et $ord(\phi'- \phi)>i$, ce que nous d\'esirions.   
\vspace{10pt}

$\underline {2^{ieme} cas:}$ $\hspace{10pt} (0,\varphi) \in \mathcal{V}(\mathcal{J})$.

Il existe alors $k_{\varphi},a_{\varphi} \in \mathbb{N}$ tels que $\forall \phi \in \mathcal{A}_{(\mathbb{K}^{n},0)}$, on ait:
$$(ord(\phi)>k_{\varphi} \hspace{5pt} et \hspace{5pt} Min(ord(\Delta(\phi, \varphi)),Min_{1 \leq j \leq k_{0}} ord(f_{j}(\phi,\varphi)))>a_{\varphi}+\alpha i)$$
$$\Longrightarrow  (\exists \phi' \in \mathcal{A}_{(\mathbb{K}^{n},0)}\vert ord(\phi'-\phi)>i \hspace{5pt} et \hspace{5pt} (\phi',\varphi) \in \mathcal{V}(\mathcal{J}))$$
Soit alors $\phi \in \mathcal{A}_{(\mathbb{K}^{n},0)}$ tel que: $ord(\phi)>k_{\varphi}$ et $Min_{1 \leq j \leq k_{0}}(ord(f_{j}(\phi,\varphi)))>2a_{\varphi}+2\alpha i$.

\vspace{10pt}
   
Si $ord(\Delta(\phi,\varphi))>a_{\varphi}+\alpha i$, il existe $\phi' \in \mathcal{A}_{(\mathbb{K}^{n},0)}$ avec $ord(\phi' -\phi)>i$ et $(\phi',\varphi)\in \mathcal{V}(\mathcal{J}) \subset \mathcal{V}(\mathcal{I})$. Nous avons alors fini. On peut donc supposer que : 
$$ord(\Delta(\phi,\varphi)) \leq a_{\varphi} + \alpha i.$$  
Donc en particulier : $ord(u(\phi,\varphi)) \leq a_{\varphi} + \alpha i.$ 

On proc\`ede alors comme dans le premier cas, en consid\'erant le syst\`eme d'\'equations implicites : 
$$ U_{j}(t,z)=u_{j}(   \phi_{1}(t)+ z_{1},\ldots,\phi_{l}(t)+ z_{l},\phi_{l+1}(t),\ldots,\phi_{n}(t),\varphi(t)  )=0, \hspace{5pt} 1 \leq j \leq l. $$
On a : 
$$U_{j}(t,0) \in (det(\frac{\partial U_{k} }{\partial z_{j}}(t,0))^{2}.t^{ 2(a_{\varphi}+ \alpha i -ord (u(\phi,\varphi))) +1 }) .$$
D'o\`u l'existence de $\theta_{1},\ldots,\theta_{l} \in \mathbb{K}\{t\}$ tels que : 
$$ord(\theta_{s}) \geq a_{\varphi}+ \alpha i + 1 + (a_{\varphi}+ \alpha i -ord (u(\phi,\varphi))) \geq a_{\varphi}+ \alpha i + 1, \hspace{5pt}, 1 \leq s \leq l, $$ 
et $U_{s}(t,\theta_{1}(t),\ldots,\theta_{l}(t))=0, 1\leq s\leq l.$

Posons $\phi'=(\phi_{1}+\theta_{1},\ldots,\phi_{l}+\theta_{l},\phi_{l+1},\ldots,\phi_{n})$. On a : $ord(\Delta(\phi',\varphi))$ $= ord(\Delta(\phi,\varphi))$ $\leq a_{\varphi} + \alpha i$. On conclut alors comme dans le premier cas.  

\endproof

\section{Suite des multiplicit\'es de Nash d'un espace analytique le long d'un arc}
Les notions de suite des mutiplicit\'es de Nash (Resp. suite de Nash des fonctions de Hilbert-Samuel et suite de Nash des diagrammes des exposants initiaux) d'un espace analytique le long d'un arc ont \'et\'e introduite dans la d\'efinition 1.5. Pour les notations employ\'ees, nous prions le lecteur de se reporter \`a l'introduction.\\ 
 Le th\'eor\`eme 1.6 d\'ecoule directement du r\'esultat suivant. Celui-ci ne surprendra pas les "d\'esingularisateurs". On en trouve trace au moins dans le cas complexe, sous une forme un peu diff\'erente, dans [L-J-T] th. 2.13 p.242. Nous donnons ici une preuve valable en r\'eel et complexe, qui proc\`ede par r\'eduction au cas hypersurface. Cette preuve est nous semble-t-il plus effective que celle de [L-J-T].

 Pour \'enoncer ce r\'esultat, consid\'erons un germe d'espace analytique $(X,x)$ $\hookrightarrow (\mathbb{K}^{n},0)$, et $\varphi \in \mathcal{A}_{(\mathbb{K}^{n},0)}$ un arc non trivial. Nous notons par $(C,0)$ le germe de courbe image de $(\mathbb{K},0)$ par $\varphi$. Soit $H_{X,C}$ la valeur g\'en\'erique de la fonction de Hilbert-Samuel de $X$ le long de la courbe $C$. C'est donc la valeur commune des $H_{X,\varphi (t)}$ pour $t \neq 0$ assez petit.

\vspace{10pt}
Consid\'erons maintenant le diagramme commutatif suivant:

\begin{equation*}
\raisebox{1.6cm}
{\xymatrix{ 
\ar[r]^-\id         & (\K,\o) \ar[r]^-\id \ar[d]^{\varphi_i} & (\K,\o) \ar[r]^-\id \ar[d]^{\varphi_{i-1}} & \cdots \ar[r]^-\id     & (\K,\o) \ar[r]^-\id \ar[d]^{\varphi_1} & (\K,\o) \ar[d]^{\varphi_0} \\
\ar[r]^-{\alpha_{i+1}} & (B_i,\o_i) \ar[r]^-{\alpha_i}            & (B_{i-1},\o_{i-1}) \ar[r]^-{\alpha_{i-1}}    & \cdots \ar[r]^-{\alpha_2} & (B_1,\o_1) \ar[r]^-{\alpha_1}            & (B_0,\o_0)                \\
\ar[r]^-{\rho_{i+1}}   & (X_i,\o_i) \ar@{^{(}->}[u]_{k_i} \ar[r]^-{\rho_i} & (X_{i-1},\o_{i-1}) \ar@{^{(}->}[u]_{k_{i-1}} \ar[r]^-{\rho_{i-1}}      & \cdots \ar[r]^-{\rho_2}   & (X_1,\o_1) \ar@{^{(}->}[u]_{k_1} \ar[r]^-{\rho_1}              & (X_0,\o_0) \ar@{^{(}->}[u]_{k_0}
}}
\end{equation*}
o\`u:
\begin{itemize}
\item[-] $(B_{0},O_{0})=(\mathbb{K}^{n},0)$ et $(X_{0},O_{0})=(X,x)$;
\item[-] $\alpha _{i}$ est l'\'eclatement de $B_{i-1}$ de centre $O_{i-1}$;
\item[-] $\varphi _{i}$ rel\`eve $\varphi _{i-1}$ \`a travers $\alpha _{i}$, $O_{i}=\varphi_{i}(0)$;
\item[-] $X_{i}$ est le transform\'e strict de $X_{i-1}$ par $\alpha_{i}$, $\rho _{i}= \alpha _{i}/X_{i}$ et les $k_{i}$ sont les plongements induits.
\end{itemize}
\begin{theo}
Soit $H_{i}=H_{X_{i},O_{i}}$ la fonction de Hilbert-Samuel du germe $(X_{i},O_{i})$. La suite $H_{i}$ est d\'ecroissante et il existe $i_{0}$ tel que $H_{i}=H_{X,C}$,$\forall i \geq i_{0}$. En particulier la suite des multiplicit\'es $m_{i}$ successives est d\'ecroissante et se stabilise \`a la multiplicit\'e g\'en\'erique de $X$ le long de $C$. Par suite, pour $\varphi \in \mathcal{R}_{(X,0)}$, les germes $(X_{i},O_{i})$ sont lisses pour $i \geq i_{0}$.
\end{theo}
Avant de prouver 4.1, notons que celui-ci entraine 1.6. En effet, partant d'un germe $(X,x)$ et d'un arc $\varphi \in \mathcal{A}_{(\mathbb{K}^{n},0)}$, on applique 4.1 \`a $Z_{0}=(\mathbb{K}\hspace{3pt},0)\times (X,x)$ et \`a l'arc $\Gamma _{0}$. Ceci fournit le  point 2) de 1.6. Pour le point 1), il suffit de remarquer que la multiplicit\'e de $X$ en $y$ est la m\^eme que celle de $Z_{0}$ en $(t,y)$.

\vspace{10pt}

\noindent \textit {Preuve de 4.1}:

D\'emontrons d'abord le cas o\`u $(X,x)$ est un germe d'hypersurface \`a l'origine de $\mathbb{K}^{n}$. Dans ce cas, la fonction d'Hilbert-Samuel est enti\`erement d\'etermin\'ee par la multiplicit\'e. Il s'agit donc de voir que la suite des multiplicit\'es successives $m_{i}$ des $(X_{i}, O_{i})$ tend vers $m'_{0}$, o\`u $m'_{0}$ est la multiplicit\'e g\'en\'erique de $X $ le long de $(C,0)$. Soit $m'_{i}$ la multiplicit\'e g\'en\'erique de $(X_{i}, O_{i})$ le long de $(C_{i},O_{i})$, o\`u $(C_{i}, O_{i})=Im \varphi_{i}$. On a: 

\begin{itemize}
\item[i)] $m'_{i} \leq m_{i}$ par semi-continuit\'e sup\'erieure de la multiplicit\'e. 
\item[ii)] $m'_{i}=m'_{0}$ car $\alpha_{i}$ est un isomorphisme de $B_{i}-\alpha_{i}^{-1}(O_{i-1}) \rightarrow B_{i}-\{O_{i-1}\}$.  
\end{itemize} 

\vspace{8pt}

Donc $m'_{0} \leq m_{i}$ pour tout $i \in \mathbb{N}$. 

D'autre part, soit $f \in \mathbb{K} \{x_{1},\ldots,x_{n}\}=\mathcal{O}_{n}$, un germe d'\'equation pour notre hypersurface. 
$$f(x,y)= \sum_{\alpha \in \mathbb{N}^n} \frac{1}{\alpha !} \frac {\partial ^{\vert \alpha \vert  }f(x)}{\partial x^{\alpha} } y^{\alpha}. $$ 
Consid\'erons, pour $k \in \mathbb{N}$, $J_{k}$ l'id\'eal de $\mathcal{O}_{n}$ engendr\'e par les d\'eriv\'ees partielles $\frac {\partial ^{\vert \alpha \vert  }f(x)}{\partial x^{\alpha}}$ avec $\vert \alpha \vert \leq k$. Alors:
$$m'_{0}=Min\{k \in \mathbb{N}/\varphi^{*}J_{k} \neq 0 \}.$$ 
Posons $D=ord \varphi^{*}(J_{m'_{0}})<+\infty$. 
On a alors: $\forall \hspace{5pt} i \geq D,\hspace{3pt} m_{i}=m'_{0}$. 

En effet, si cela n'\'etait pas le cas, on aurait, compte-tenu de la d\'ecroissance de la suite des multiplicit\'es (cas particulier du th\'eor\`eme de Bennett [B]): 
$$(*) \hspace{10pt} m_{i} > m'_{0} \hspace{10pt} pour \hspace{10pt} 0 \leq i \leq D.$$
Ceci est absurde, comme il va r\'esulter du lemme suivant: 
\begin{lem}
Soit $i \in \mathbb{N}$, notons $f_{i} \in \mathcal{O}_{n}$, la i.\`eme transform\'ee stricte de $f$, et pour $k \in \mathbb{N}$, consid\'erons $J_{i,k}$ l'id\'eal de $\mathcal{O}_{n}$ engendr\'e par les $\frac {\partial ^{\vert \alpha \vert  }f_{i}(x)}{\partial x^{\alpha} }$ avec $\vert \alpha \vert \leq k$. Alors si $k<m_{i}$, on a: 
$$ord \hspace{5pt} (\varphi_{i}^{*} J_{i,k}) \geq 1 + ord \hspace{5pt} (\varphi_{i+1}^{*} J_{i+1,k})    .$$
\end{lem} 
Appliquant successivement ceci, pour $i$ variant de $0$ \`a $D$, avec $k=m'_{0}$, on aurait alors par $(*)$:
$$ ord \hspace{5pt} (\varphi^{*} J_{m'_{0}}) \geq (D+1) + ord \hspace{5pt} (\varphi_{D+1}^{*} J_{D+1,m'_{0}}) .   $$ 
Ce qui est absurde car $ord \hspace{5pt} (\varphi^{*} J_{m'_{0}})=D$. Il nous reste donc \`a prouver (4.2). 

\vspace{10pt}
\noindent $\textit{Preuve de 4.2}:$ 

On peut trouver un syst\`eme de coordonn\'ees et un isomorphisme local: $$u_{i+1}:\hspace{5pt} (B_{i+1},O_{i+1})  \rightarrow (\mathbb{K}^{n},0)$$
 tel que le germe d'\'eclatement 
$$\alpha_{i+1}:(B_{i+1},O_{i+1}) \simeq (\mathbb{K}^{n},0) \rightarrow (B_{i},O_{i}) \simeq (\mathbb{K}^{n},0)$$
 soit repr\'esent\'e par: 
$$(x_{1},x_{2},\ldots,x_{n}) \rightarrow (x_{1},x_{1}(a_{2}+x_{2}),\ldots,x_{1}(a_{n}+x_{n})),\hspace{3pt}a_{i} \in \mathbb{K}, \hspace{3pt}i \geq 2.$$
On a alors:
$$f_{i}(x_{1},x_{1}(a_{2}+x_{2}),\ldots,x_{1}(a_{n}+x_{n}))=x_{1}^{m_{i}}u(x_{1},\ldots,x_{n}) f_{i+1}(x_{1},\ldots,x_{n})  $$ 
o\`u $u$ est un inversible de $\mathcal{O}_{n}$. 

On v\'erifie par r\'ecurrence sur $k \leq m_{i}$, que pour $\vert \alpha \vert \leq k$, on a: 
$$\frac {\partial ^{\vert \alpha \vert  }f_{i}}{\partial x^{\alpha} } (x_{1},x_{1}(a_{2}+x_{2}),\ldots,x_{1}(a_{n}+x_{n})) \in x_{1}^{m_{i}-k} J_{i+1,k}.    $$
Composant ceci avec $\varphi_{i+1}$ (\'ecrite dans le syst\`eme de coordonn\'ees), on obtient: 
$$ord \hspace{5pt} (\varphi_{i}^{*} J_{i,k}) \geq (m_{i}-k)ord \hspace{5pt} (\varphi_{i}) + ord \hspace{5pt} (\varphi_{i+1}^{*} J_{i+1,k}) \geq 1 + ord \hspace{5pt} (\varphi_{i+1}^{*} J_{i+1,k})      $$    
Ceci ach\`eve la preuve du cas o\`u $(X,0)$ est un germe d'hypersurface. Notons dans ce cas que pour $i \geq D$, $X_{i}$ est \'equimultiple le long de $C_{i}=Im \varphi_{i}$.
\endproof

\vspace{8pt}

Prouvons \`a pr\'esent le cas g\'en\'eral. Soit $H_{X_{i},C_{i}}$ la valeur g\'en\'erique de la fonction de Hilbert-Samuel de $X_{i}$ le long de $C_{i}=Im \varphi_{i}$. 

Chaque $\alpha_{i}$ induisant  un isomorphisme de $B_{i}-\alpha_{i}^{-1}(O_{i})$ sur $B_{i-1}-\{O_{-1}\}$, on a: 
\begin{itemize}
\item[i)] la suite  $H_{X_{i},C_{i}}$ est constante et \'egale à $H_{X,C}$;\item[ii)] $H_{i+1}=H_{X_{i+1},O_{i+1}} \leq H_{i}=H_{X_{i},O_{i}}  $ par le th\'eor\`eme de Bennett [B];
\item[iii)] par semi-continuit\'e de la fonction de Hilbert-Samuel (cf. [L-J-T]): 
$$H_{X,C}=H_{X,C_{i}} \leq H_{i}=H_{X_{i},O_{i}}.$$
\end{itemize}   
Par cons\'equent,  montrer qu'il existe  $i_{0}$ tel que, $ \forall i \geq i_{0}$ on ait $H_{i}=H_{X,C}$, c'est montrer qu'il existe  $i_{0}$ tel que $ \forall i \geq i_{0}$ on ait $H_{i}=H_{X_{i},C_{i}}$. 
Ce qui signifie qu'il s'agit de voir que $(C_{i},O_{i})\subset (S_{i},O_{i})$ o\`u $S_{i}$ est la strate de Hilbert-Samuel de $X_{i}$ passant par $O_{i}$ i.e.:
$$ S_{i}=\{x \in X_{i}\vert H_{X_{i},x}=H_{X_{i},O_{i}} \}.$$ 
Plongeons $(X_{i},O_{i})$ dans $(\mathbb{K}^n,0)$ et choisissons un syst\`eme de coordonn\'ees quelconques, $y_{1},\ldots,y_{n}$, \`a l'origine de $\mathbb{K}^{n}$. Ecrivons:
$$\mathcal{O}_{X_{i},O_{i}}=\frac{\mathbb{K}\{y\}}{\mathfrak{a}_{i}}.$$Si $f_{i,1},\ldots,f_{i,p_{i}}$ d\'esigne la base standard distingu\'ee de $\mathfrak{a}_{i}$ dans ce syst\`eme de coordonn\'ees, on a:

$$ S_{i}= \cap _{j=1}^{p_{i}} S_{i,j}.$$
Avec $S_{i,j}=\{x \in \mathbb{K}^{n}\vert m_{x}(f_{i,j})=m_{O_{i}}(f_{i,j})\}$, o\`u $m_{x}$ d\'esigne la multiplicit\'e ou l'ordre au point $x$  (cf. [B-M2] Th. 5.3.1 p.821). 

Notre probl\`eme revient donc \`a voir qu'il existe $i_{0}$, tel que pour tout $i \geq i_{0}$, on puisse trouver un syst\`eme de coordonn\'ees \`a l'origine de $(\mathbb{K}^{n},0)$ tel que chaque \'el\'ement de la base standard distingu\'ee $f_{i,j}$, avec $1 \leq j \leq p_{i}$, relative \`a ce syst\`eme de coordonn\'ees, soit \'equimultiple le long de $C_{i}$. Nous allons prouver cette affirmation. Puisque pour tout $i$:
$$H_{i+1} \leq H_{i},   $$
par le th\'eor\`eme de stabilisation 5.2.2 p. 820 de [B-M2], on obtient l'existence d'un $i_{0}$ tel que $\forall i \geq i_{0}$: 
$$H_{i}=H_{i_{0}}.$$ 
Fixons un tel entier $i_{0}$, plongeons $(X_{i_{0}},O_{i_{0}})$ dans $(\mathbb{K}^{n},0)$, et choisissons un syst\`eme de coordonn\'ees (à l'origine de $\mathbb{K}^{n}$) $(w,z)$ avec $ w=(w_{1},\ldots,w_{n-r})$, $z=(z_{1},\ldots,z_{r})$ satisfaisant  les conclusions du lemme 7.2 p.828 de [B-M2] (des variables essentielles dans la terminologie de [B-M2]).

Consid\'erons maintenant $f_{i_{0},1},\ldots,f_{i_{0},p_{i_{0}}}$, la base standard distingu\'ee de $\mathfrak{a}_{i_{0}}$ relativement \`a ce syst\`eme de coordonn\'ees, o\`u $\mathcal{O}_{X_{i_{0}},O_{i_{0}}}= \frac {\mathbb{K} \{w,z\}}{\mathfrak{a}_{i_{0}}}$. D'apr\`es le th\'eor\`eme 7.3 p. 828 de [B-M2], il existe un syst\`eme de coordonn\'ees $(w',z')$ dans un voisinage de $O_{i_{0}+1}=0$ tel que: 
\begin{itemize}
\item[-] si $\mathfrak{a}_{i_{0}+1} \subset \mathbb{K}\{w',z'\}$ d\'esigne le transform\'e strict de $\mathfrak{a}_{i_{0}}$ (donc $\mathcal{O}_{X_{i_{0}+1},O_{i_{0}+1}}$ $\simeq \mathbb{K}\{w',z'\}/\mathfrak{a}_{i_{0}+1}$), la base standard distingu\'ee de $\mathfrak{a}_{i_{0}+1}$ est constitu\'ee par les transform\'ees strictes $f_{i_{0}+1,j}$ des $f_{i_{0},j}$, $1 \leq j \leq p_{i_{0}}$. 
\item[-] le syst\`eme de coordonn\'ees $(w',z')$ satisfait les conclusions du lemme 7.2 p.828 de [B-M2].
\end{itemize} 
On peut donc it\'erer ce qui pr\'ec\`ede. D'apr\`es le cas hypersurface, au bout d'un nombre fini d'\'etapes D, les $f_{i_{0}+k,j}$, $1 \leq j \leq p_{i_{0}}$, $k \geq D$ deviennent \'equimultiples le long de  $C_{i_{0}+k}$. Le th\'eor\`eme en r\'esulte comme nous l'avons vu plus haut. 
\endproof

\begin{rem}

Soit $\varphi  \in \mathcal{R}_{(X,0)}$. Supposons que $(X,0)$ soit r\'eduit et \'equi\-di\-men\-sion\-nel  de dimension $d$. Soit $I \subset \mathcal{O}_{n}$ tel que $\mathcal{O}_{X,0}= \mathcal{O}_{n}/I$. Si $J$ d\'esigne l'id\'eal engendr\'e par $I$ et tous les jacobiens d'ordre $n-d$ des $n-d$ upplets d'\'el\'ements de $I$, et si $D=ord( \varphi^{*}(J))$, on a $(X_{i}, O_{i})$ lisse pour $i \geq D$. Il suffit pour le constater de reprendre le calcul du lemme 4.2. 
\end{rem} 

Le th\'eor\`eme 4.1 et la remarque 4.3 nous permettent de d\'ecrire une base de voisinage d'un \'el\'ement de $\mathcal{R}_{(X,0)}$. Pour cela, plongeons $(X,0)$ dans $(\mathbb{K}^{n},0)$ et notons par d la dimension en $0$ de $(X,0)$. Soit $z=(z',z")$, $z'=(z_{1},\ldots,z_{d})$,   $z"=(z_{d+1},\ldots,z_{n})$ un syst\`eme de coordonn\'ees \`a l'origine de $\mathbb{K}^{n}$. Pour $\varphi  \in \mathcal{R}_{(X,0)}$, notons comme pr\'edemment:
$$\varphi(t)=(\varphi_{1}(t),\ldots,\varphi_{n}(t)) = \sum_{k=1}^{+\infty} A_{k}.t^{k}.$$
Soit $\xi \in B_{\mathbb{K}^{n},i}(\varphi)$. Posons: $\xi(t)=\sum_{k=0}^{i} A_{k}.t^{k} +t^{i}r_{i}(\xi)(t)$,
o\`u $r_{i}(\xi) \in \mathcal{A}_{(\mathbb{K}^{n},0)}$. 
 Notons: $r_i(\xi )= (r'_i(\xi ),r"_(\xi ))$, $r'_{i}(\xi)=(r_{i,1}(\xi), \ldots, r_{i,d}(\xi))\in \mathcal{A}_{(\mathbb{K}^d,0)}$, $r"_{i}(\xi)=(r_{i,d+1}(\xi), \ldots, r_{i,n}(\xi))\in \mathcal{A}_{(\mathbb{K}^{n-d},0)}$.

\begin{cor}
Soit $ \varphi(t)=\sum_{k=1}^{+\infty} A_{k}.t^{k}$, $\varphi  \in \mathcal{R}_{(X,0)}$. Il existe $i_{0} \in \mathbb{N}$, tel que $\forall i \geq i_{0}$, apr\`es une \'eventuelle permutation des coordonn\'ees dans $\mathbb{K}^{n}$, il existe des germes de fonctions analytiques nulles \`a l'origine dans $\mathbb{K}^{d+1}$, $u_{i,j}(z',t)$,  $d+1 \leq j \leq n$, tels que pour $\xi \in B_{(\mathbb{K}^{n},i)} (\varphi)$ on ait :   
$$\xi \in B_{X,i}(\varphi) \Longleftrightarrow r"_{i}(\xi)(t)=(u_{i,d+1}(r'_{i}(\xi)(t),t),\ldots,u_{n}(r'_{i}(\xi)(t),t))   $$
\end{cor}     

\noindent $\textit{Preuve}$:\\
Conservant les notations qui suivent $(D_{2})$, on a $\xi \in B_{X,i}(\varphi)$ si et seulement si l'arc $t \rightarrow (t,r_{i}(\xi)(t)) \in \mathcal{A}_{(Z'_{i},0)}$. Il ne nous reste plus qu'\`a d\'ecrire $\mathcal{A}_{(Z'_{i},0)}$. Fixons $i_{0}$ de telle sorte que pour $i \geq i_{0}$, $(Z'_{i},0)$ soit lisse. Comme $dim_{0} (Z'_{i},0)=d+1$, celui-ci peut \^etre d\'efini par:
$$(Z'_{i},0)=\{(t,z) \in (\mathbb{K}\hspace{3pt},0)\times (\mathbb{K}^{n},0)\vert f_{d+1}(t,z)=\ldots=f_{n}(t,z)=0\},$$
o\`u le rang en $0$ de la matrice jacobienne des $f_{j}$, $d+1 \leq j \leq n$, est $n-d$. Pour conclure, il nous suffit de voir, par le th\'eor\`eme des fonctions implicites ordinaires, qu'il existe $i_{d+1},\ldots,i_{n}$ $\in \{1,\ldots,n\}$ tels que:
$$det(\frac{\partial f_{k}}{\partial z_{i_{l}}})_{d+1 \leq k,l \leq n}(0) \neq 0.$$
Ceci est n\'ecessairement satisfait. En effet, si cela n'\'etait pas le cas, puisque le rang de la matrice jacobienne des $f_{j}$ est tout de m\^eme $n-d$, cela signifierait que l'on peut d\'efinir (apr\`es une \'eventuelle permutation des $z_{j}$) $(Z'_{i},0)$ comme:
$$(Z'_{i},0)=\{(t,z) \in (\mathbb{K}\hspace{3pt},0)\times (\mathbb{K}^{n},0)/t-u_{1}(z_{1},\ldots,z_{d+1})=z_{d+2}-u_{2}(z_{1},\ldots,z_{d+1})$$ $$=\ldots=z_{n}-u_{n-d}(z_{1},\ldots,z_{d+1})=0\}$$
avec $u_{j}(0)=0$, $1\leq j \leq n-d$ et $\partial u_{1}/ \partial z_{1}(0)=\partial u_{1}/ \partial z_{2}(0)=\ldots=\partial u_{1}/ \partial z_{d+1}(0)=0$. 

Mais ceci est absurde puisque l'arc $t \rightarrow (t,r_{i}(\varphi)(t))$ est \`a valeurs dans $(Z'_{i},0)$. Or celui-ci ne peut satisfaire l'\'equation:
$$t-u_{1}(r_{i,1}(\varphi)(t),\ldots,r_{i,d+1}(\varphi)(t))=0,$$
car $ord(u_{1}(r_{i}(\varphi))\geq 2$.
\endproof
\begin{rem}

1) Soit $\varphi \in \mathcal{R}_{(X,0)}$. Fixons $i_{0}$ satisfaisant les conditions de $4.4$. Soit $i \geq i_{0}$, alors l'application $\Theta _{i}$:
$$B_{X,i}(\varphi) \rightarrow \mathcal{A}_{(\mathbb{K}^{d},0)}$$
$$\hspace{10pt} \xi \rightarrow r'_{i}(\xi)$$
est un "isomorphisme analytique" puisque sa r\'eciproque $\Theta ^{-1}_{i}$ est donn\'ee par:
$$\mathcal{A}_{(\mathbb{K}^{d},0)} \rightarrow B_{X,i}(\varphi)$$
$$\hspace{10pt} \alpha \rightarrow \Theta ^{-1}_{i}(\alpha)$$
o\`u $\Theta ^{-1}_{i}(\alpha)(t)=(\alpha_{1}(t),\ldots,\alpha_{d}(t),u_{i,d+1}(\alpha (t),t),\ldots,u_{i,n}(\alpha (t),t))$. On peut d'au\-tre part consid\'erer que la notion naturelle de morphisme analytique entre $\mathcal{A}_{(X,0)}$ et $\mathcal{A}_{(Y,0)}$ est celle des morphismes induits par les applications analytiques:
$$(\mathbb{K}\hspace{3pt},0)\times (X,0) \rightarrow (\mathbb{K}\hspace{3pt},0)\times (Y,0)$$
$$\hspace{10pt} (t,x) \rightarrow (t,a(t,x)).$$
Pour ces raisons, il nous semble justifi\'e de dire que de ce point de vue $\mathcal{R}_{(X,x)}$ est l'ensemble des points r\'eguliers de $\mathcal{A}_{(X,x)}$.

\vspace{5pt}
2) Si $(X,0)$ est r\'eduit et \'equidimensionnel, l'entier $i_{0}$ est pr\'ecis\'e par 4.3.

\vspace{5pt}
3) Sp\'ecifier la d\'ependance des $u_{i,j}(z',t)$ par rapport aux $A_{k}$, $k\leq i$, est la notion de stabilit\'e de $\textrm{[D-L1]}$. Ceci d\'ecoulera du th\'eor\`eme 4.10. L'entier $i_{0}$ de 4.3 est plus pr\'ecis que celui fourni par la preuve du lemme 4.1 de [D-L1]. Il est utile de pr\'eciser cet entier $i_{0}$ si on cherche \`a obtenir des informations pr\'ecises sur la s\'erie de Poincar\'e consid\'er\'ee dans $\textrm{[D-L1]}$.
\end{rem}

En g\'en\'eral $\mathcal{R}_{(X,0)}$ n'est pas dense dans $\mathcal{A}_{(X,0)}$ comme le montre l'exemple suivant:
\begin{ex}
Consid\'erons le germe en $0$ du parapluie de Whitney:
$$(X,0)=\{(z_{1},z_{2},z_{3}) \in (\mathbb{C}^{3},0)\vert z_{1}^{2}-z_{2}z_{3}^{2}=0 \}$$
$(Sing(X),0)=\{(z_{1},z_{2},z_{3}) \in (\mathbb{C}^{3},0)\vert z_{1}=z_{3}=0\}$. Soit l'arc $\varphi$, $t \rightarrow (0,t,0)$. Il est facile de v\'erifier que $\forall i \geq 1$, $B_{X,i}(\varphi) \subset \mathcal{A}_{(Sing(X),0)}$. Par cons\'equent $B_{X,i}(\varphi) \cap \mathcal{R}_{(X,0)}=\varnothing$.
\end{ex}
 Cependant ce genre de pathologie peut toujours \^etre \'evit\'e par un changement de param\'etrisation. Plus pr\'ecis\'ement on a:
\begin{lem}
Soient $(X,0)$ un germe d'espace analytique complexe r\'eduit et $\varphi \in \mathcal{A}_{(X,0)}$. Alors il existe un entier $p \in \mathbb{N}^{*}$ tel que, d\'esignant par $\theta$ l'arc $t \rightarrow \varphi(t^{p})$, on ait, $\forall \hspace{5pt} (Y,0)$ sous-ensemble analytique propre de $(X,0)$:
$$\forall i \in \mathbb{N}, \hspace{3pt}B_{X,i}(\theta)\cap(\mathcal{A}_{(X,0)}-\mathcal{A}_{(Y,0)})\neq \varnothing.$$
\end{lem}
\noindent $\textit{Preuve:}$

Fixons un repr\'esentant de $(X,0)$ dans un voisinage de $0$ assez petit. Soit alors $\Pi:  Z \rightarrow X$, avec $Z$ lisse d\'enombrable \`a l'infini, $\Pi$ propre et surjectif (un tel $\Pi$ existe en vertu des r\'esultats sur la d\'esingularisation). Soit maintenant $C$ la courbe image de $\varphi$. Posons $A'=\Pi^{-1}(C)$. Consid\'erons la d\'ecomposition de $A'$ en composantes irr\'eductibles $A'=\cup_{\alpha \in \Lambda} A'_{\alpha}$ ($\Lambda$ est d\'enombrable). Pour chaque $\alpha \in \Lambda$, puisque $\Pi$ est propre, $\Pi(A'_{\alpha})$ est un sous-ensemble analytique irr\'eductible inclus dans $C$. C'est donc un point ou $C$ tout entier. Comme $\Lambda$ est d\'enombrable, il existe $\alpha$ tel que $\Pi(A'_{\alpha})=C$. Fixons un tel $\alpha$, disons $\alpha=1$. $\Pi^{-1}(0)\cap A'_{1}$ est un sous-ensemble analytique de $A'_{1}$ tel que pour tout $p \in A'_{1}$, $dim_{p} \Pi^{-1}(0)\cap A'_{1}<dim_{p}A'_{1}$. Soit $p \in \Pi^{-1}(0)\cap A'_{1}$, consid\'erons un arc $\rho: (\mathbb{C},0) \rightarrow (A'_{1},p)$ tel que $\rho \in \mathcal{A}_{(A'_{1},p)}-\mathcal{A}_{( \Pi^{-1}(0)\cap A'_{1},p)}$. On peut supposer que $\rho$ d\'efinit la normalisation de son image $(C_{0},p)$. On a alors un diagramme commutatif induit par $\Pi$:

\begin{equation*}
\begin{CD}
(C_{0},p)@>\Pi>>(C,0)\\
@AA{\rho} A  @AA{n}A\\
(\mathbb{C},0)@>\tilde{\Pi}>>(\mathbb{C},0) @<{a}<< (\mathbb{C},0)
\end{CD}
\end{equation*}
o\`u $\varphi=n \circ a$, $n$ est est la normalisation de $(C,0)$, $\tilde{\Pi}$ le rel\`evement de $\Pi$ \`a la normalisation. Par un choix convenable de coordonn\'ee dans le $(\mathbb{C},0)$ source de  $\tilde{\Pi}$, on peut \'ecrire $\tilde{\Pi}(t)=t^{p}$. Ainsi $\Pi(\rho (t))=n(t^{p})$. Ecrivons maintenant $a(t)=t^{k}u(t)$, $u(0)\neq 0$. Puis soit $\alpha \in \mathcal{O}_{1}$ tel que $\alpha (t)^{p}=u(t^{p})$. On a:
$$\varphi(t^{p})=n(t^{kp}u(t^{p}))=n((t^{k} \alpha(t))^{p})=\Pi(\rho(t\alpha(t))).$$
Si $\tilde{\rho}$ d\'esigne l'arc $t \rightarrow \rho(t\alpha (t))$ $\in \mathcal{A}_{(Z,p)}$, on a donc $\Pi(\tilde{\rho}(t))=\varphi(t^{p})=\theta (t)$. Soit maintenant $(Y,0)$ un sous-ensemble analytique propre de $(X,0)$. Notons $Y'=\Pi^{-1}(Y)$, $(Y',p)$ est un sous-ensemble analytique propre de $(Z,p)$. Soit $i \in \mathbb{N}$, il existe en vertu du lemme 4.9 ci-dessous $\tilde{\xi}\in B_{Z,i}(\tilde{\rho})$ tel que $\tilde{\xi} \not\in \mathcal{A}_{(Y',p)}$. Par suite l'arc $\xi=\Pi(\tilde{\xi})$ $\in B_{X,i}(\theta)$ et $\not\in \mathcal{A}_{(Y,0)}$. 
\endproof

\begin{rem}
Le lemme 4.7 ne sera pas utilis\'e par la suite. C'est le seul endroit o\`u nous utilisons l'existence de d\'esingularisation. Il a \'et\'e motiv\'e par une affirmation dans l'article de J. Becker r\'ef\'erenc\'e dans $\textrm{[I5]}$.
\end{rem}
Le lemme 4.9 est une version effective du second point du lemme 1.4 de [T4]. Cette effectivit\'e sera utilis\'ee \`a la section suivante.
\begin{lem}
Soient $\varphi \in \mathcal{A}_{(\mathbb{K}^{n},0)}$ et $i \in \mathbb{N}$. Pour $f \in \mathbb{K}[[X_{1},\ldots,X_{n}]]$, $f \neq 0$ on a:
$$m_{0}(f)\leq Min_{\theta \in B_{\mathbb{K}^{n},i}(\varphi)}ord(f\circ \theta)\leq m_{0}(f)(i+1)$$
o\`u $m_{0}(f)$ d\'esigne la multiplicit\'e ou l'ordre de $f$.
\end{lem}

\noindent $\textit{Preuve:}$

Soient $\varphi \in \mathcal{A}_{(\mathbb{K}^{n},0)}$ et $i \in \mathbb{N}$. Pour prouver 4.9, on peut supposer $ord(\varphi) \leq i$. En effet sinon, soit $(a_{1},\ldots,a_{n})$ ne figurant pas dans les z\'eros de la forme initiale de $f$. Alors l'arc $t \rightarrow \theta (t)=(a_{1}t^{i+1},\ldots,a_{n}t^{i+1})$ est dans $B_{\mathbb{K}^{n},i}(\varphi)$ et $ord(f\circ \theta)=m_{0}(f)(i+1)$. Nous supposerons donc $ord(\varphi)\leq i$. On peut par ailleurs supposer que $f \in \mathcal{O}_{n}$. Soit $(Z_{0},0)=(\mathbb{K}\hspace{3pt},0)\times (X,x)$, o\`u $(X,x)$ est le germe d\'efini par $f$. Posons: $f_{0} \in \mathbb{K}\{t,X\}$, $f_{0}=f$. Puis soit $f_{j}(t,x) \in \mathcal{O}_{n+1}=\mathbb{K}\{t,X\}$, la transform\'ee stricte de $f_{j-1}$ par $\Pi'_{j}$ (les notations sont celles de $(D_{2}))$. Un calcul \'el\'ementaire en coordonn\'ees locales montre que:
$$t^{m_{0}(f_{0})+m_{0}(f_{1})+\ldots+m_{0}(f_{i-1})}f_{i}(t,X)=f(\sum _{k=1}^{i}A_{k}t^{k} +t^{i}X)$$
pour $(t,X)$ voisin de $(0,0)$.

Soit alors un point $(1,a_{1},\ldots,a_{n})$ qui n'est pas dans les z\'eros de la forme initiale de $f_{i}$. On a alors:
$$t^{m_{0}(f_{0})+m_{0}(f_{1})+\ldots+m_{0}(f_{i-1})}f_{i}(t,a_{1}t,\ldots,a_{n}t)=f(\sum _{k=1}^{i}A_{k}t^{k} +t^{i+1}A)$$
avec $A=(a_{1},\ldots,a_{n})$. D'o\`u:
$$ord(f(\sum _{k=1}^{i}A_{k}t^{k} +t^{i+1}A))=m_{0}(f_{0})+m_{0}(f_{1})+\ldots+m_{0}(f_{i-1})+m_{0}(f_{i}).$$
Or c'est un cas extr\^emement simple (v\'erifiable \`a la "main") du th\'eor\`eme de Benett [B]  de voir que:
$$m_{0}(f_{i})\leq m_{0}(f_{i-1})\leq \ldots\leq m_{0}(f_{0}).$$
Par suite $ord(f(\sum _{k=1}^{i}A_{k}t^{k} +t^{i+1}A))\leq m_{0}(f)(i+1)$.
\endproof

\vspace{10pt}

Nous allons \`a pr\'esent prouver le th\'eor\`eme de semi-continuit\'e 1.7.\\

\noindent $\textit{Preuve de 1.7.:}$

Consid\'erons une sous-vari\'et\'e alg\'ebrique affine irr\'eductible $V \subset$ $\mathbb{K}^{l}$. Notons $\mathbb{K}[V]$ la $\mathbb{K}$-alg\`ebre de type fini qui la d\'efinit. Soient $F_{1}(t,X),\ldots,F_{s}(t,X)$ des \'el\'ements de $\mathbb{K}[V][[t,X]]$. Posons:
$$F_{j}(t,x)=\sum_{\alpha _{1} \in \mathbb{N},\alpha \in \mathbb{N}^{n}} a^{j}_{\alpha_{1},\alpha}t^{\alpha _{1}}X^{\alpha}, \hspace{3pt} a^{j}_{\alpha_{1},\alpha} \in \mathbb{K}[V].$$
Soient $R=\sum_{j=1}^{s} \mathbb{K}[V][[t,X]]F_{j}(t,X)$ et $N(R)$ le diagramme des exposants initiaux de $R$. Pour chaque $v \in V$, il y a un morphisme "d'\'evaluation des coefficients en $v$":
$$\mathbb{K}[V][[t,X]] \longrightarrow \mathbb{K}[[t,X]]$$
$$\hspace{10pt} F(t,X) \longrightarrow F(v,t,X)$$
o\`u  si $F(t,x)=\sum_{\alpha _{1} \in \mathbb{N},\alpha \in \mathbb{N}^{n}} a_{\alpha_{1},\alpha}t^{\alpha _{1}}X^{\alpha}$, $ a_{\alpha_{1},\alpha} \in \mathbb{K}[V]$, 
$$F(v,t,x)=\sum_{\alpha _{1} \in \mathbb{N},\alpha \in \mathbb{N}^{n}} a_{\alpha_{1},\alpha}(v)t^{\alpha _{1}}X^{\alpha}, \hspace{5pt} a_{\alpha_{1},\alpha}(v)\in \mathbb{K}.$$
On notera par $R_{v}$ l'id\'eal de $\mathbb{K}[[t,X]]$ engendr\'e par les $F_{j}(v,t,X)$ et par $N(R_{v})$ son diagramme des exposants initiaux. On a alors:
\begin{lem}
Sous les notations pr\'ec\'edentes:
\begin{itemize}
\item[1)] $\forall v \in V$,  $N(R) \leq N(R_{v})$;
\item[2)] Il existe $V'$ sous-vari\'et\'e alg\`ebrique propre de V telle que:
\item[-] $\forall v \in V-V'$, $N(R)=N(R_{v})$;
\item[-] pour tout sommet $\beta$ de $N(R)$, il existe $G \in R$ tel que $\nu(G)=\beta$ et $\forall v \in V-V'$, $\nu(G(v,t,X))=\nu(G)=\beta$.

\end{itemize}
\end{lem}
\noindent $\textit{Preuve:}$

Il suffit de r\'ep\'eter dans nos notations la preuve des lemmes 7.1 et 7.2 de E. Bierstone et P.D. Milman dans [B-M1].

\begin{lem}
Soient $i \in \mathbb{N}^{*}$  et $V$ une sous-vari\'et\'e alg\'ebrique irr\'eductible de $\mathcal{A}_{(\mathbb{K}^{n},0)}^{i}=\mathbb{K}^{ni}$, il existe une sous-vari\'et\'e propre $V'$ de $V$ telle que:
\begin{itemize}
\item[1)] $\forall \varphi^{i},\hspace{3pt}\theta^{i} \in V-V', \hspace{3pt} \mathcal{H}_{X,\varphi^{i}}^{i}=\mathcal{H}_{X,\theta^{i}}^{i}$;
\item[2)] $\forall \varphi^{i},\hspace{3pt}\theta^{i} \in V-V', \hspace{3pt} \mathcal{N}_{X,\varphi^{i}}^{i}=\mathcal{N}_{X,\theta^{i}}^{i}$.
\end{itemize}
\end{lem}

\noindent $\textit{Preuve:}$

Nous tentons d'abord d'all\'eger les notations. Soit $\varphi^{i}$ un \'el\'ement de $\mathcal{A}_{(\mathbb{K}^{n},0)}^{i}$, c'est donc un morphisme de $\mathbb{K}$-alg\`ebres locales de $\mathbb{K}\{X_{1},\ldots,X_{n}\}$$\rightarrow \mathbb{K}\{t\}/(t)^{i+1}$. La donn\'ee de $\varphi^{i}$ est donc \'equivalente \`a la donn\'ee d'un \'el\'ement:
$$\varphi^{i}=\sum_{k=1}^{i}A_{k}t^{k},\textrm{ avec } A_{k} \in \mathbb{K}^{n}.$$
Les notations \'etant celles de $(D_{1})$ et $(D_{2})$, nous noterons par $H_{A^{j}}=H_{A_{1},\ldots,A_{j}}$, $A^{j}=(A_{1},\ldots,A_{j})$, la fonction de Hilbert-Samuel de $(Z'_{j},0)$. C'est donc la fonction de Hilbert-Samuel de l'anneau local $\mathbb{K}\{t,X\}/I_{j,A^{j}}$. Dans ces notations, la suite $\mathcal{H}_{X,\varphi^{i}}^{i}$ est la suite $(H_{0},H_{A_{1}},\ldots,H_{A_{1},\ldots,A_{i}})$, o\`u $H_{0}$ est la fonction de Hilbert-Samuel de $(Z'_{0},0)=(\mathbb{K}\hspace{3pt},0)\times(X,0)$. Soit $j$, $1 \leq j \leq i$. D\'esignons par $V_{j}$ l'adh\'erence de Zariski de $P_{j}(V) \subset \mathbb{K}^{nj}$, o\`u $P_{j}:$ $\mathbb{K}^{ni}\rightarrow \mathbb{K}^{nj}$ d\'esigne la projection oublie des $n(i-j)$ derniers termes. $V_{j}$ est donc une sous-vari\'et\'e irr\'eductible de $\mathbb{K}^{nj}$. Consid\'erons maintenant une base standard, $f_{1}(X),\ldots,f_{p_{0}}(X)$ $\in \mathbb{K}[[X]]$ de $I_{0}$, o\`u $I_{0}$ est l'id\'eal d\'efinissant $\mathcal{O}_{X,0}$ comme $\mathbb{K}\{X\}/I_{0}$. Puis posons:
$$ F_{k,0}(t,X)=f_{k}(X) \in \mathbb{K}[[X]], \hspace{5pt} 1\leq k \leq p_{0}.$$
Soit alors:
$$F_{k,1}(t,X)=\frac{F_{k,0}(t,t(Y_{1}+X))}{t^{\vert \nu(F_{k,0})\vert}} \in \mathbb{K}[Y_{1}][[t,X]], 1 \leq k \leq p_{0}.$$
Prenant la classe des coefficients des $F_{k,1}$ dans $\mathbb{K}[V_{1}]$, on obtient des \'el\'ements $\overline{F}_{k,1}(t,X)$ $\in \mathbb{K}[V_{1}][[t,X]]$. On notera $R_{1}$ l'id\'eal de $\mathbb{K}[V_{1}][[t,X]]$ engendr\'e par les $\overline{F}_{k,1}(t,X)$ et soit $N(R_{1})$ son diagramme des exposants initiaux. D\'esignons par $\beta_{1,1},\ldots,\beta_{p_{1},1}$ les sommets de $N(R_{1})$. D'apr\`es le lemme pr\'e\-c\'e\-dent il existe une sous-vari\'et\'e propre $V'_{1}$ de $V_{1}$ et des \'el\'ements $G_{1,1},\ldots,G_{p_{1},1}$ $\in R_{1} \subset \mathbb{K}[V_{1}][[t,X]]$ tels que:
\begin{itemize}
\item[-] $\forall A^{1} \in V_{1}-V'_{1}$, $N(R_{1})=N(R_{1,A^{1}})$;
\item[-] $\forall A^{1} \in V_{1}-V'_{1}$, $\forall k$, $1 \leq k \leq p_{1}$, $\beta_{k,1}=\nu(G_{k,1})=\nu(G_{k,1}(A^{1},t,X))$.
\end{itemize}

\vspace{8pt}

Posons ensuite:
$$F_{k,2}(t,X)=\frac{G_{k,1}(t,t(Y_{2}+X))}{t^{\vert \nu(G_{k,1})\vert}} \in \mathbb{K}[V_{1}][Y_{2}][[t,X]], 1 \leq k \leq p_{1}.$$
Prenant la classe des coefficients des $F_{k,2}$ dans $\mathbb{K}[V_{2}]$, on obtient des \'el\'ements $\overline{F}_{k,2}(t,X)$ $\in \mathbb{K}[V_{2}][[t,X]]$. On notera $R_{2}$ l'id\'eal de $\mathbb{K}[V_{2}][[t,X]]$ engendr\'e par les $\overline{F}_{k,2}(t,X)$ et soit $N(R_{2})$ son diagramme des exposants initiaux. D\'esignons par $\beta_{1,2},\ldots,\beta_{p_{2},2}$ les sommets de $N(R_{2})$. A nouveau le lemme 4.10 assure l'existence d'une sous-vari\'et\'e propre $V'_{2}$ de $V_{2}$ et d'\'el\'ements $G_{1,2},\ldots,G_{p_{2},2}$ $\in R_{2} \subset \mathbb{K}[V_{2}][[t,X]]$  tels que:
\begin{itemize}
\item[-] $\forall A^{2}=(A_{1},A_{2}) \in V_{2}-V'_{2}$, $N(R_{2})=N(R_{2,A^{2}})$;
\item[-] $\forall A^{2}=(A_{1},A_{2}) \in V_{2}-V'_{2}$, $\forall k$, $1 \leq k \leq p_{2}$, $\beta_{k,2}=$ $\nu(G_{k,2})=$ $\nu(G_{k,2}(A^{2},t,X))$.
\end{itemize}

\vspace{8pt}

On continue le proc\'ed\'e par r\'ecurrence. Soit $j<i$. Ayant d\'etermin\'e une sous-vari\'et\'e propre $V'_{j}$ de $V_{j}$ et des \'el\'ements $\overline{F}_{k,j}(t,X)$ $\in \mathbb{K}[V_{j}][[t,X]]$, $1 \leq k \leq p_{j-1}$, engendrant un id\'eal $R_{j}$, il existe des \'el\'ements $G_{1,j},\ldots,G_{p_{j},j}$ $\in R_{j} \subset \mathbb{K}[V_{j}][[t,X]]$ tels que:
\begin{itemize}
\item[-] $\forall A^{j}=(A_{1},\ldots,A_{j}) \in V_{j}-V'_{j}$, $N(R_{j})=N(R_{j,A^{j}})$;
\item[-] $\forall A^{j}=(A_{1},\ldots,A_{j}) \in V_{j}-V'_{j}$, $\forall k$, $1 \leq k \leq p_{j}$, $\beta_{k,j}$ $=$ $\nu(G_{k,j})$ $=$ $\nu(G_{k,j}(A^{j},t,X))$,
\end{itemize}
o\`u les $\beta_{k,j}$ sont les sommets de $N(R_{j})$. 

On pose alors:
$$F_{k,j+1}(t,X)=\frac{G_{k,j}(t,t(Y_{j+1}+X))}{t^{\vert \nu(G_{k,j})\vert}} \in \mathbb{K}[V_{j}][Y_{j+1}][[t,X]], 1 \leq k \leq p_{j}.$$
On note $\overline{F}_{k,j+1}(t,X)$ $1 \leq k \leq p_{j}$, les \'el\'ements de $\mathbb{K}[V_{j+1}][[t,X]]$ obtenus en prenant les classes des coefficients. Soient $R_{j+1}$ l'id\'eal de $\mathbb{K}[V_{j+1}][[t,X]]$ engendr\'e par les $\overline{F}_{k,j+1}(t,X)$ et $N(R_{j+1})$ son diagramme des exposants initiaux. Le lemme 4.10 assure l'existence d'\'el\'ements $G_{k,j+1} \in R_{j+1}$, $1 \leq k \leq p_{j+1}$ et d'une sous-vari\'et\'e prorpre $V'_{j}$ de $V_{j}$ tels que:
\begin{itemize}
\item[-] $\forall A^{j+1}=(A_{1},\ldots,A_{j+1}) \in V_{j+1}-V'_{j+1}$, $N(R_{j+1})=N(R_{j+1,A^{j+1}})$;
\item[-] $\forall A^{j+1}=(A_{1},\ldots,A_{j+1}) \in V_{j+1}-V'_{j+1}$, $\forall k$, $1 \leq k \leq p_{j+1}$, $\beta_{k,j+1}$ $=\nu(G_{k,j+1})$ $=\nu(G_{k,j+1}(A^{j+1},t,X))$,
\end{itemize}
o\`u les $\beta_{k,j+1}$ sont les sommets de $N(R_{j+1})$.

On pose alors:
$$V'=\bigcup_{1 \leq j \leq i}((V'_{j}\times \mathbb{K}^{n(i-j)})\cap V)).$$
Alors $V'$ est une sous-vari\'et\'e propre de $V$. Soient alors $A^{i}=(A_{1},\ldots,A_{i})$ (resp.
$B^{i}=(B_{1},\ldots,B_{i})$) $\in V-V'$. Puisque le diagramme des exposants initiaux d\'etermine la fonction de Hilbert-Samuel, il nous suffit pour conclure de voir que:
$$N(I_{j,A^{j}})=N(I_{j,B^{j}}),\hspace{5pt} 1 \leq j \leq i.$$
Pour cela du fait que les $F_{k,0}$ sont une base standard de $I_{0}$, leurs transform\'ees strictes engendrent $I_{1,A^{1}}$ pour tout $A^{1} \in V_{1}$. Par cons\'equent pour tout $A^{1} \in V_{1}$, $I_{1,A^{1}}$ est engendr\'e par les $\overline{F}_{k,1}(A^{1},t,X)$, $1 \leq k \leq p_{0}$. Maintenant si $A^{1} \in V_{1}-V'_{1}$, par construction m\^eme les $G_{k,1}(A^{1},t,X)$ sont une base standard de $I_{1,A^{1}}=R_{1,A^{1}}$. Par suite $N(I_{1,A^{1}})=N(I_{1,B^{1}})=N(R_{1})$ et donc $H_{A^{1}}=H_{B^{1}}$. De plus, les transform\'ees strictes des $G_{k,1}(A^{1},t,X)$ engendrent $I_{2,A_{1},A_{2}}$, $\forall A_{1} \in V_{1}-V'_{1}, A_{2} \in \mathbb{K}^{n}$. Par cons\'equent les $\overline{F}_{k,2}(A_{1},A_{2},t,X)$, $1 \leq k \leq p_{1}$, engendrent $I_{2,A_{1},A_{2}}$,  $\forall A^{1} \in V_{1}-V'_{1}, A^{2} \in \mathbb{K}^{n}$. Par suite, $\forall A_{1} \in V_{1}-V'_{1}, A_{2} \in \mathbb{K}^{n}$, $R_{2,A_{1},A_{2}}=I_{2,A_{1},A_{2}}$. Par cons\'equent:
$$N(I_{2,A_{1},A_{2}})=N(I_{2,B_{1},B_{2}})=N(R_{2}), \textrm{ et } H_{A_{1},A_{2}}=H_{B_{1},B_{2}}.$$
Par construction, $\forall (A_{1},A_{2}) \in V_{2}-V'_{2}$, les $G_{k,2}(A_{1},A_{2},t,X)$ sont une base standard de $I_{2,A_{1},A_{2}}$. Leurs transform\'ees strictes engendrent donc $I_{3,A^{3}}$, $\forall A^{3} \in (V_{2}-V'_{2})\times \mathbb{K}^{n}$. Par suite, les $F_{k,3}(A^{3},t,X)$ engendrent $I_{3,A^{3}}$, $\forall A^{3} \in (V_{2}-V'_{2})\times \mathbb{K}^{n}$. Donc $\forall A^{3} \in V_{3}-V'_{3}$, $N(I_{3,A^{3}})=N(R_{3})=N(R_{3,A^{3}})$. D'o\`u:
$$H_{A^{3}}=H_{B^{3}}.$$
On continue ainsi par r\'ecurrence jusqu'\`a l'ordre $i$. Notons  que pour tout $\varphi^{i} \in V-V'$,
$$\mathcal{N}^{i}_{X,\varphi ^{i}}=(N(I_{1,A_{1}}),\ldots,N(I_{i,A^{i}}))=(N(R_{1}),\ldots,N(R_{i}))$$
\endproof

\begin{lem}
Soient $i \in \mathbb{N}^{*}$ et $V$ une sous-vari\'et\'e alg\'ebrique irr\'eductible de $\mathcal{A}_{(\mathbb{K}^{n},0)}^{i}=\mathbb{K}^{ni}$. Soit $V'$ comme dans le lemme 4.11. On a:
 $$\forall \hspace{3pt} \varphi ^{i} \in V-V',\hspace{5pt} \forall \theta^{i} \in V',\hspace{5pt} \mathcal{N}_{X,\varphi^{i}}^{i}\leq \mathcal{N}_{X,\theta^{i}}^{i}.$$
\end{lem}

\noindent \textit{Preuve:}

On proc\`ede par r\'ecurrence sur $i$. Pour $i=1$, la propri\'et\'e  r\'esulte du lemme 4.10 puisque  pour tout $A_{1}\in V$, $I_{1,A_{1}}=R_{1,A_{1}}$.  Nous supposerons donc $i>1$ et la propri\'et\'e \'etablie pour $i-1$. Pour \'etablir la propri\'et\'e au rang i, on raisonne par r\'ecurrence sur $m$ la dimension de $V$. Si $m=0$, le r\'esultat est \'evident. On peut donc supposer $m>0$ et la propri\'et\'e \'etablie pour toute sous-vari\'et\'e irr\'eductible de $\mathbb{K}^{ni}$ de dimension $\leq m-1$. Soient $V$, $V'$ comme dans 4.11. Consid\'erons la d\'ecomposition en composantes irr\'eductibles de $V'$:
$$V'=\bigcup _{1\leq l\leq t}V'^{l}.$$
Le lemme 4.11 fournit pour chaque $l$ une sous-vari\'et\'e propre $V"^{l}$ telle que:
$$\forall \hspace{3pt} \varphi'^{i},\hspace{3pt} \forall \hspace{3pt} \theta'^{i}\hspace{3pt} \in V'^{l}-V"^{l},\hspace{3pt} \mathcal{N}_{X,\varphi'^{i}}^{i}=\mathcal{N}_{X,\theta'^{i}}^{i}.$$ 
Notre hypoth\`ese de r\'ecurrence sur la dimension nous dit que:
$$\forall \hspace{3pt} \varphi'^{i}\hspace{3pt} \in V'^{l}-V"^{l},\hspace{3pt}\forall\hspace{3pt} \theta'^{i}\hspace{3pt} \in V"^{l},\hspace{3pt} \mathcal{N}_{X,\varphi'^{i}}^{i}\leq \mathcal{N}_{X,\theta'^{i}}^{i}.$$
Par cons\'equent pour \'etablir la propri\'t\'e au rang $i$, on peut supposer que $V'$ est irr\'eductible et qu'il existe une sous-vari\'et\'e propre $V"$ de $V'$ comme dans 4.11. Il s'agit alors de voir que:
$$\forall\hspace{3pt}\varphi^{i} \in V-V',\hspace{3pt}\forall\hspace{3pt}\theta^{i} \in V'-V",\hspace{3pt} \mathcal{N}_{X,\varphi^{i}}^{i}\leq \mathcal{N}_{X,\theta^{i}}^{i}.$$
Maintenant par hypoth\`ese de r\'ecurrence sur $i$, notant $\varphi^{i-1}$ et $\theta^{i-1}$ les tronqu\'es \`a l'ordre $i-1$ de $\varphi^{i}$ et $\theta^{i}$, on a:
$$\mathcal{N}_{X,\varphi^{i-1}}^{i-1}\leq \mathcal{N}_{X,\theta^{i-1}}^{i-1}.$$
Si l' in\'egalit\'e ci-dessus est stricte, par definition de l'ordre lexicographique on a alors aussi une in\'egalit\'e stricte au rang $i$. On peut donc supposer que:
$$\mathcal{N}_{X,\varphi^{i-1}}^{i-1}= \mathcal{N}_{X,\theta^{i-1}}^{i-1}.$$

La preuve du lemme 4.11, appliqu\'ee \`a $V'$, nous dit qu'il existe des sous-vari\'et\'es propres $V"_{j}$ de $V'_{j}$, $1\leq j \leq i$, des \'el\'ements $\overline{F}'_{k,j}$ $\in \mathbb{K}[V'_{j}][[t,X]]$, $1\leq k \leq p'_{j-1}$, engendrant un id\'eal $R'_{j}$ dont le diagramme des exposants initiaux est not\'e $N(R'_{j})$, et des \'el\'ements $G'_{k,j} \in R'_{j}$ tels que:
\begin{itemize}
\item[-] $\forall B^{j}=(B_{1},\ldots,B_{j})\in V'_{j}-V"_{j}$, $N(R'_{j})=N(R'_{j,B^{j}})$;
\item[-] $\forall B^{j}=(B_{1},\ldots,B_{j})\in V'_{j}-V"_{j}$, $\nu(G'_{k,j})=\nu(G'_{k,j}(B^{j},t,X))$, $1 \leq k \leq p'_{j},$ o\`u les $\beta'_{k,j}$ sont les sommets de $N(R'_{j})$.
\end{itemize}
On a de plus:
$$F'_{k,j+1}(t,X)=\frac{G'_{k,j}(t,t(Y_{j+1}+X))}{t^{ \vert \nu(G'_{k,j})\vert}} \in \mathbb{K}[V'_{j}][Y_{j+1}][[t,X]].$$
Dont les classes donnent les $\overline{F}'_{k,j+1} \in \mathbb{K}[V'_{j+1}][[t,X]]$. Par hypoth\`ese, on a $p'_{j-1}=p_{j-1}$ et $\beta'_{k,j}=\beta_{k,j}$, $j\leq i-1$.
 Pour conclure il suffit d'\'etablir la propri\'et\'e  suivante:

\begin{slem}\textrm{ }\\
 Il existe des \'el\'ements $\tilde{F}'_{k,i}$ $\in \mathbb{K}[V][[t,X]]$, $1 \leq k \leq p_{i-1}$, tels que:
\begin{itemize}
\item[1)]  L'id\'eal de $\mathbb{K}[[t,X]]$ engendr\'e par les $\tilde{F}'_{k,i}(B_{1},\ldots,B_{i},t,X)$ est le m\^eme que celui engendr\'e par les $\overline{F}'_{k,i}(B_{1},\ldots,B_{i},t,X)$, $\forall B^{i}=(B_{1},\ldots,B_{i})$ $\in$  $V'-W"$,  o\`u $W"$ est une sous-vari\'et\'e propre de $V'$;
\item[2)]$\tilde{F}'_{k,i}(A_{1},\ldots,A_{i-1},A_{i},t,X) \in I_{i,A_{1},\ldots,A_{i-1},A_{i}}$, $\forall \hspace{3pt} A^{i}=(A_{1},\ldots,A_{i}) \in \hspace{3pt} V$.
\end{itemize}
\end{slem}

En effet supposons $(4.13)$ \'etablie. Notons $\tilde{R}_{i}$ l'id\'eal de $\mathbb{K}[V][[t,X]]$ engendr\'e par les $\tilde{F}'_{k,i}$ et par $N(\tilde{R}_{i})$ son diagramme des exposants initiaux. Par le second point de $(4.13)$, on a pour tout $A^{i} \in V$, $\tilde{R}_{i,A^{i}} \subset I_{i,A^{i}}$. Donc $N(I_{i,A^{i}})$ $\leq$ $ N(\tilde{R}_{i,A^{i}})$. Maintenant par le lemme 4.10, il existe une sous vari\'et\'e propre $W'$ de $V$ telle que $\forall A^{i}$  $\in$ $V-W'$, $N(\tilde{R}_{i,A^{i}})=N(\tilde{R}_{i})$. Par irr\'eductibilit\'e de $V$, $V'\cup W' \varsubsetneq V$. Consid\'erant un point $A^{i} \in V-(V'\cup W')$ on a donc:
$$N(R_{i})=N(I_{i,A^{i}})\leq N(\tilde{R}_{i,A^{i}})=N(\tilde{R}_{i})$$
Maintenant par le premier point de $(4.13)$, consid\'erons un point $B^{i}$ de $V'-(V"\cup W")$ (qui est n\'ecessairement non vide par irr\'eductibilit\'e de $V'$), on a:
$$N(R'_{i})=N(I_{i,B^{i}})=N(\tilde{R}_{i,B^{i}}).$$
Par suite par le lemme 4.10, $N(\tilde{R}_{i})\leq N(R'_{i})$. D'o\'u $N(R_{i})\leq N(R'_{i})$, ce que nous cherchions \`a d\'emontrer  . Il nous reste donc \`a prouver $(4.13)$.

\vspace{10pt}
Nous allons \`a nouveau proc\'eder par r\'ecurrence sur $i$. Les notations \'etant celles de la preuve du lemme  4.11 on a:
$$F_{k,1}(t,X)=\frac{F_{k,0}(t,t(Y_{1}+X))}{t^{\vert \nu(F_{k,0})\vert}} \in \mathbb{K}[Y_{1}][[t,X]],1\leq k \leq p_{0}.$$
Les $\overline{F}'_{k,1}$ sont les classes dans $\mathbb{K}[V'_{1}][[t,X]]$ des $F_{k,1}$. Posons: $\tilde{F}'_{k,1}(t,X)=\overline{F}_{k,1}(t,X) \in \mathbb{K}[V_{1}][[t,X]]$. $(4.13)$ est alors claire pour $i=1$. D'autre part, il existe des $\overline{U}^{1}_{k,l}$ $\in$  $\mathbb{K}[V'_{1}][[t,X]]$ tels que:
$$G'_{k,1}=\sum_{l=1}^{p_{0}} \overline{U}^{1}_{k,l}\overline{F}'_{l,1} \in \mathbb{K}[V'_{1}][[t,X]].$$
Posons:
$$G"_{k,1}=\sum_{l=1}^{p_{0}} U^{1}_{k,l}\tilde{F}'_{l,1} \in \mathbb{K}[V_{1}][[t,X]]$$
o\`u les $U^{1}_{k,l}$ repr\'esentent les $\overline{U}^{1}_{k,l}$. Nous allons construire des $\tilde{G}'_{k,1}$ $\in$  $\mathbb{K}[V_{1}][[t,X]]$, $1 \leq k \leq p_{1}$, tels que:
\begin{itemize}
\item[-] $\tilde{G}'_{k,1}$ $\in$ $R_{1}$;
\item[-] $\nu (\tilde{G}'_{k,1})=\beta_{k,1}$;
\item[-] $\forall \hspace{3pt} B_{1} \hspace{3pt} \in V'_{1}-V"_{1}$, $\nu (\tilde{G}'_{k,1}(B_{1},t,X))=\nu (G'_{k,1}(B_{1},t,X))$ $= \beta_{k,1}$;
\item[-] $\forall \hspace{3pt} B_{1} \hspace{3pt} \in V'_{1}-V"_{1}$, $\tilde{G}'_{k,1}(B_{1},t,X)$ $=G'_{k,1}(B_{1},t,X)$ modulo un \'el\'ement de $\mathbb{K}^{*}$.
\end{itemize}
On proc\`ede par r\'ecurrence sur $k$. On a d'abord:
$$\nu(G"_{1,1})\leq \nu(G'_{1,1})=\beta_{1,1}.$$
Mais $\nu(G"_{1,1})$ $\in$ $N(R_{1})$. Donc $\nu(G"_{1,1}) \geq \beta_{1,1}$. Par suite $\nu(G"_{1,1})= \nu(G'_{1,1})$. Posons $\tilde{G}'_{1,1}=G"_{1,1}$ $\in$ $\mathbb{K}[V_{1}][[t,X]]$. Soit $k$, $1\leq k <p_{1}$, supposons avoir construit $\tilde{G}'_{1,1},\ldots,\tilde{G}'_{k,1}$ $\in$ $\mathbb{K}[V_{1}][[t,X]]$ satisfaisant les propri\'et\'es ci-dessus. On a: $\nu(G"_{k+1,1})$ $\leq \nu(G'_{k+1,1})=\beta_{k+1,1}$. Mais $\nu(G"_{k+1,1})$ $\in$ $N(R_{1})$. Donc:
$$\nu(G"_{k+1,1}) \hspace{3pt}\in \cup_{k+1 \leq s \leq p_{1}}(\beta_{s,1}+ \mathbb{N}^{n+1}),$$
ou bien
$$\nu(G"_{k+1,1}) \hspace{3pt} \in\cup_{1 \leq s \leq k}(\beta_{s,1}+ \mathbb{N}^{n+1}).$$
Dans le premier cas on a $\nu(G"_{k+1,1})\geq \beta_{k+1,1}$ et donc $\nu(G"_{k+1,1}) = \beta_{k+1,1}$. On pose alors $\tilde{G}'_{k+1,1}=G"_{k+1,1}$. Dans le second cas si $\nu(G"_{k+1,1})<\nu(G'_{k+1,1})$, on \'ecrit:
$$\nu(G"_{k+1,1})=\beta_{l,1}+\gamma,\textrm{ avec } 1\leq l \leq k \textrm{ et } \gamma \in \mathbb{N}^{n+1}.$$
On consid\`ere alors:
$$G'''_{k+1,1}(t,X)=\tilde{g}'_{\beta_{l,1},l,1}(Y_{1})G"_{k+1,1}(t,X)-g"_{\beta_{l,1}+\gamma,k+1,1}(Y_{1})t^{\gamma_{1}}X^{\gamma'}\tilde{G}'_{l,1}(t,X).$$
O\`u $\gamma=(\gamma_{1},\gamma')$ et $\tilde{g}'_{\beta_{l,1},l,1}(Y_{1})$ est le coefficient dans le d\'eveloppement de $\tilde{G}'_{l,1}$ correspondant au mon\^ome $t^{\beta_{l,1,1}}X^{\beta'_{l,1}}$ avec $\beta_{l,1}=(\beta_{l,1,1},\beta'_{l,1})$ et $g"_{\beta_{l,1}+\gamma,k+1,1}(Y_{1})$ est celui correspondant au mon\^ome $t^{\beta_{l,1,1}+\gamma_{1}}X^{\beta'_{l,1}+\gamma'}$ dans le developpement de $G"_{k+1,1}$. On a alors:
$$(\mathbf{**})\hspace{20pt}\nu(G"_{k+1,1})<\nu(G'''_{k+1,1})\leq \nu(G'_{k+1,1})=\beta_{k+1,1}.$$
De plus:
\begin{itemize}
\item[-] $\forall$ $B_{1}$ $\in$ $V'_{1}-V"_{1}$, $G'''_{k+1,1}(B_{1},t,X)$ $=$ $G'_{k+1,1}(B_{1},t,X)$ \`a la multiplication pr\`es par un \'el\'ement de $\mathbb{K}^{*}$;
\item[-] $\forall$ $B_{1}$ $\in$ $V'_{1}-V"_{1}$, $\nu(G'''_{k+1,1}(B_{1},t,X))=\beta_{k+1,1}$;
\item[-] $G'''_{k+1,1}$ $\in$ $R_{1}$.
\end{itemize}
On recommence alors le m\^eme raisonnement avec $G'''_{k+1,1}$ au lieu de $G'_{k,1}$. Puisque la premi\`ere in\'egalit\'e de $(\mathbf{**})$ est stricte, au bout d'un nombre fini d'\'etapes, on aura construit un \'el\'ement $\tilde{G}'_{k+1,1}$ $\in$ $\mathbb{K}[V_{1}][[t,X]]$ poss\'edant les propri\'et\'es d\'esir\'ees.

Maintenant ayant construit $\tilde{G}'_{1,1},\ldots,\tilde{G}'_{p_{1},1}$, on pose:
$$\tilde{F}'_{k,2}(t,X)=\frac{\tilde{G}'_{k,1}(t,t(Y_{2}+X))}{t^{\vert \beta_{k,1} \vert}} \in \mathbb{K}[V_{1}][Y_{2}][[t,X]], 1 \leq k \leq p_{1},$$
et on note encore $\tilde{F}'_{k,2}$ les \'el\'ements de $\mathbb{K}[V_{2}][[t,X]]$ obtenus en prenant la classe des coefficients dans $\mathbb{K}[V_{2}]$. Par construction m\^eme $\tilde{F}'_{k,2}(B_{1},B_{2},t,X)$ et $\overline{F}'_{k,2}(B_{1},B_{2},t,X)$ diff\`erent au plus multiplicativement d'une constante non nulle pour tout $(B_{1},B_{2})$ $\in$ $V'_{2}-V"_{2}$. Le premier point de $(4.13)$ en d\'ecoule pour $i=2$. Concernant le second point, pour tout $A_{1}$ $\in$ $V_{1}$, $\tilde{G}'_{k,1}(A_{1},t,X) \in R_{1,A_{1}}=I_{1,A_{1}}$. Donc pour tout $A_{2}$ $\in$ $\mathbb{K}^{n}$, leurs transform\'ees strictes:
$$\frac{\tilde{G}'_{k,1}(A_{1},t,t(A_{2}+X))}{t^{\vert  \nu(\tilde{G}'_{k,1}(A_{1},t,X)) \vert}}$$
sont des \'el\'ements de $I_{2,A_{1},A_{2}}$. Mais:
$$\vert \nu(\tilde{G}'_{k,1}(A_{1},t,X)) \vert \geq \vert \nu(\tilde{G}'_{k,1}(t,X)) \vert= \vert \beta_{k,1} \vert .$$
Donc les $\tilde{F}'_{k,2}(A_{1},A_{2},t,X)$ sont des multiples de ces transform\'ees strictes et sont par suite des \'el\'ements de $I_{2,A_{1},A_{2}}$. Donc $(4.13)$ est \'etablie pour $i=2$.

Maintenant pour passer du rang $i-1$ au rang $i$, on raisonne exactement comme nous venons de le faire pour le passage du rang $1$ au rang $2$. On consid\`ere  l'id\'eal $\tilde{R}_{i-1}$ de $\mathbb{K}[V_{i-1}][[t,X]]$ engendr\'e par les $\tilde{F}'_{k,i-1}$, ainsi que l'id\'eal $\tilde{R}'_{i-1}$ de $\mathbb{K}[V'_{i-1}][[t,X]]$ engendr\'e par les classes des $\tilde{F}'_{k,i-1}$. On applique le lemme 4.10 \`a $\tilde{R}'_{i-1}$ , ce qui fournit des \'el\'ements $\hat{G}'_{k,i-1}$ analogues des $G'_{k,1}$. Ceux-ci permettent de construire par le m\^eme proc\'ed\'e que dans le passage du rang $1$ au rang $2$ des $\tilde{G}'_{k,i-1}$ $\in$ $\mathbb{K}[V_{i-1}][[t,X]]$, $1 \leq k \leq p_{i-1}$, tels que:
\begin{itemize}
\item[-] $\nu(\tilde{G}'_{k,i-1})=\beta_{k,i-1}$;
\item[-] $\nu(\tilde{G}'_{k,i-1}(B_{1},\ldots,B_{i-1},t,X))=\beta_{k,i-1}$ pour tout $(B_{1},\ldots,B_{i-1})$ dans un ouvert de Zariski non vide de $V'_{i-1}$;
\item[-] $ \tilde{G}'_{k,i-1}(B_{1},\ldots,B_{i-1},t,X)=\hat{G}'_{k,i-1}(B_{1},\ldots,B_{i-1},t,X)$ \`a la multiplication pr\`es par un \'el\'ement de $\mathbb{K}^{*}$, pour tout $(B_{1},\ldots,B_{i-1})$ dans un ouvert de Zariski non vide de $V'_{i-1}$;
\item[-] $\tilde{G}'_{k,i-1}$ $\in$ $\tilde{R}_{i-1}$.
\end{itemize}
 Ceux-ci permettent de d\'efinir les $\tilde{F}'_{k,i}$ par la formule:
$$\tilde{F}'_{k,i}(t,X)=\frac{\tilde{G}'_{k,i-1}(t,t(Y_{i}+X))}{t^{\vert \beta_{k,i-1} \vert}} \in \mathbb{K}[V_{i}][[t,X]], 1 \leq k \leq p_{i-1}.$$
On v\'erifie alors les propri\'et\'es de $(4.13)$ comme plus haut.
\endproof

\vspace{5pt}

Pour finir la preuve du th\'eor\`eme de semi-continuit\'e, il nous suffit de dire que celui-ci r\'esulte des lemmes 4.11 et 4.12.
\endproof

\vspace{5pt}

\noindent \textit{Preuve de 1.8:}\\

La semi-continuit\'e pour la topologie de Zariski d'une fonction  $\mathcal{F}^{i}$ de $\mathcal{A}^{i}_{(\mathbb{K}^{n},0)}$ dans un ensemble ordonn\'e $\Lambda$ entraine que pour tout $\lambda \in \Lambda$ l'ensemble:
$$S^{i}_{\lambda}=\{\varphi^{i}\hspace{3pt}\in \hspace{3pt}\mathcal{A}^{i}_{(\mathbb{K}^{n},0)}\vert \mathcal{F}^{i}(\varphi^{i})\geq \lambda\}$$
est une sous-vari\'et\'e alg\'ebrique de $\mathcal{A}^{i}_{(\mathbb{K}^{n},0)}$. D'autre part d'apr\`es le th\'eor\`eme de semi-continuit\'e $\varphi^{i} \rightarrow \mathcal{N}_{X,\varphi^{i}}^{i}$ ne prend qu'un nombre fini de valeurs $\mathcal{N}_{i,j}$, $1\leq j \leq l_{i}$, sur $\overline{\Pi^{i}(\mathcal{A}_{(X,x)})}$. On pose alors:
$$U_{i,j}=\{\varphi^{i} \in \overline{\Pi^{i}(\mathcal{A}_{(X,x)})}\vert \mathcal{N}_{X,\varphi^{i}}\geq \mathcal{N}_{i,j}\} \textrm{ et } W_{i,j}=\cup_{j'\neq j}U_{i,j}.$$
Les propri\'et\'es de 1.8 sont alors claires puisque le diagramme des exposants initiaux d\'etermine la fonction de Hilbert-Samuel ainsi que la multiplicit\'e.
\endproof
\begin{rem}
1) Il est naturel de se demander aussi si les fonctions $\varphi^{i} \rightarrow \mathcal{M}_{X,\varphi^{i}}^{i}$ et $\varphi^{i} \rightarrow \mathcal{H}_{X,\varphi^{i}}^{i}$ sont semi-continues sup\'erieurement pour la topologie de Zariski sur $\mathcal{A}_{(\mathbb{K}^{n},0)}^{i}$. Ceci est au moins le cas pour les hypersurfaces. Cependant l'auteur n'a pas r\'eussi pour le moment \`a le prouver dans le cas g\'en\'eral.\\

2) On peut raffiner la partition $S_{i,j}$ ci-dessus en \'ecrivant $S_{i,j}=\cup_{1 \leq k \leq t_{i,j}}U_{i,j,k}-W_{i,j,k}$, o\`u $U_{i,j,k}$ et $W_{i,j,k}$ sont des vari\'et\'es alg\'ebriques. De telle sorte que pour chaque $k$, il existe des $F_{k,l}(t,X)$ $\in$ $R(U_{i,j,k},W_{i,j,k})[[t,X]]$, $R(U_{i,j,k},W_{i,j,k})$ d\'esignant les fonctions r\'eguli\`eres sur $ U_{i,j,k}$ \`a p\^oles dans $W_{i,j,k}$, telles que $\forall A^{i}$ $\in$ $U_{i,j,k}-W_{i,j,k}$ les $F_{k,l}(A^{i},t,X)$ soient la base standard distingu\'ee de $I_{i,A^{i}}$.
\end{rem}
\section{Strates Principales et volume motivique}
Soient $(X,x)$ un germe de $\mathbb{K}$-espace analytique et $i\geq 1$ un entier. D'apr\`es les r\'esultats de la section pr\'ec\'edente, on dispose d'une partition:
$$ \mathbb{K}^{ni}=\bigcup_{1\leq l\leq k_i}S_{i,l}$$
o\`u les $S_{i,l}$ sont des sous-ensembles constructibles de $\mathbb{K}^{ni}$ tels que:\\
- sur chaque $S_{i,l}$, $\mathcal{N}^i_{X,\varphi ^i}=(N_{0,\varphi ^i},\ldots,N_{i,\varphi ^i})$ est constante \'egale \`a $\mathcal{N}^i_{X,l}=(N_0^l,\ldots,N_i^l)$\\
- Si $l\neq l'$, $\mathcal{N}^i_{X,l}\neq\mathcal{N}^i_{X,l'}$\\
Par voie de cons\'equences $\mathcal{M}^i_{X,\varphi ^i}$ et $\mathcal{H}^i_{X,\varphi ^i}$ sont constantes sur chaque $S_{i,l}$. On notera respectivement $\mathcal{M}^{i}_{X,l}=(m_0^l,\ldots,m_i^l)$ et $\mathcal{H}^{i}_{X,l}=(H_0^l,\ldots,H_i^l)$ la valeur de ces constantes.

\begin{defi} Soient $(X,x)$ un germe de $\mathbb{K}$-espace analytique et $i\geq 1$ un entier donn\'e. On dit qu'une strate $S_{i,l}$ de la statification de $\mathbb{K}^{ni}$ par les suites de Nash des diagrammes des exposants initiaux est principale si et seulement si:
\begin{itemize}
\item[1)] $\forall m \in \mathbb{N}$, $(m,0,\ldots,0 )\not \in N_i^l$
\item[2)] $m_i^l=1$
\end{itemize}
On appelle $i^{ieme}$ partie principale de la stratification de $\mathcal{A}_{(X,x)}$ par les suites de Nash des diagrammes initiaux et on note $P^i_{(X,x)}$ la r\'eunion (disjointe) des $S_{i,l}$ qui sont des strates principales.
\end{defi}  
Soient $m \in \mathbb{N}$, $\alpha =(\alpha _0,\alpha _1,\ldots,\alpha _n)\in \mathbb{N}^{n+1}$ tels que $\vert \alpha \vert =m$. On a $\alpha \leq (m,0,\ldots,0)$. Ainsi pour $f \in \mathbb{K}[[t,X]]$, dire que son exposant initial $\nu (f)$ est \'egal \`a $(m,0,\ldots,0)$ \'equivaut \`a dire que sa forme initiale $In(f)$ vaut $a_mt^m$ avec $a_m\neq 0$. Par cons\'equent dire qu'un multi-indice $(m,0\ldots,0)$ figure dans le diagramme des exposants initiaux de $I_{i,A^i}$ pour un point $A^i$ de $\mathbb{K}^{ni}$ \'equivaut \`a dire que $t^m$ appartient \`a l'id\'eal initial $In(I_{i,A^i})$ de $I_{i,A^i}$. Si $\mathbb{K}=\mathbb{C}$, l'existence d'un tel $m$ \'equivaut au fait que le r\'eduit du c\^one tangent de $(Z'_i,0)$, not\'e $(C(Z'_i,0))_{red}$, est inclus dans $\{t=0\}$ (cette derni\`ere \'equation est celle du diviseur exeptionnel par la suite d\'eclatements d\'efinie par $A^i$). La seconde condition correspond \'evidemment au fait que $(Z'_i,0)$ est lisse. Ces deux remarques montrent en particulier que lorsque $\mathbb{K}=\mathbb{C}$, $P^i_{(X,x)}$ est ind\'ependant du syst\`eme de coordon\'ees choisi sur $\mathbb{K}^n$.\\
\begin{theo}\textrm{ }\\
Soit $(X,x)$ un germe de $\mathbb{K}$-espace analytique $\mathbb{K}$-r\'eduit et equidimensionel de dimension $d$.
\begin{itemize}
\item[1)] La suite $[P^i_{(X,x)}]\mathbb{L}^{-di}$  converge vers une limite non nulle dans $\widehat{\mathcal{M}}_{\mathbb{K}}$.
\item[2)] La suite $[\Pi ^i(\mathcal{A}_{(X,x)})]\mathbb{L}^{-di}-[P^i_{(X,x)}]\mathbb{L}^{-di}$  converge vers zero dans $\widehat{\mathcal{M}}_{\mathbb{K}}$. Ainsi:
$$Lim_{ i\rightarrow +\infty}[P^i_{(X,x)}]\mathbb{L}^{-di}= Lim_{ i\rightarrow +\infty}[\Pi ^i(\mathcal{A}_{(X,x)})]\mathbb{L}^{-di}$$
\end{itemize}
Par cons\'equent $$Lim_{ i\rightarrow +\infty}[P^i_{(X,x)}]\mathbb{L}^{-di}=\mu (\mathcal{A}_{(X,x)})$$
o\`u $\mu $ d\'esigne la mesure motivique.
\end{theo}

\noindent \textit{Preuve: }\\
Nous commen\c{c}ons par constater que $P^i_{(X,x)}\subset \Pi ^i(\mathcal{A}_{(X,x)})$. En effet, soit $A^i \in P^i_{(X,x)}$. Par le second point de 5.1 le germe correspondant $(Z'_i,0) \subset (\mathbb{K}^{n+1},0)$ est lisse de dimension $d+1$. Par cons\'equent $N(I_{i,A^i})=\cup _{1\leq k\leq n-d}(\alpha_ {i_k}+\mathbb{N}^{n+1})$ avec $\vert \alpha _{i_k}\vert=1$. Par la premi\`ere condition de 5.1 $\alpha _{i_{k}}\neq (1,0,\ldots,0)$, $1\leq k\leq n-d$. Ceci implique que la base standard distingu\'ee de $I_{i,A^i}$ est de la forme $X_{i_{k}}-U_{i_{k}}(t,X")$, $1\leq k\leq n-d$, o\`u $X"$ d\'esigne les $d$ coordonn\'ees parmi $X_1,\ldots,X_n$, autres que $X_{i_1},\ldots,X_{i_{n-d}}$. Par definition du diagramme $(D_2)$, ceci implique l'inclusion souhait\'ee. Montrons maintenant que la suite des $[P^i_{(X,x)}]\mathbb{L}^{-di}$ est de Cauchy dans $\widehat{\mathcal{M}}_{\mathbb{K}}$. Soient $i,k\in \mathbb{N}^*$, par le second point de la remarque 4.14 et la description ci-dessus, l'application:
$$P^{i+k}_{(X,x)}\cap \Pi _{i+k,i}^{-1}(P^i_{(X,x)})\longrightarrow P^i_{(X,x)}$$
d\'eduite de la projection oublie des \guillemotleft $k$ derniers termes\guillemotright de $\mathbb{K}^{n(i+k)}$ dans $\mathbb{K}^{ni}$ est une fibration localement triviale de fibtre $\mathbb{L}^{dk}$. Ainsi $[P^{i+k}_{(X,x)}\cap \Pi _{i+k,i}^{-1}(P^i_{(X,x)}]=[P^i_{(X,x)}]\mathbb{L}^{dk}$ et donc:
\begin{eqnarray*}
[P^{i+k}]\mathbb{L}^{-d(k+i)}-[P^i]\mathbb{L}^{-di}& = &[P^{i+k}]\mathbb{L}^{-d(k+i)}-[P^{i+k}\cap \Pi _{i+k,i}^{-1}(P^i)]\mathbb{L}^{-d(i+k)}\\
 & = &[P^{i+k}-P^{i+k}\cap \Pi _{i+k,i}^{-1}(P^i)]\mathbb{L}^{-d(i+k)}
\end{eqnarray*}
A un point $A^{i+k}=(A_1,\ldots,A_i,\ldots,A_{i+k})$ de $P^{i+k}-P^{i+k}\cap \Pi _{i+k,i}^{-1}(P^i)$ correspond une suite des multiplicit\'es de Nash $m_0,\ldots,m_i,\ldots,m_{i+k}$ avec $m_i\geq 2$ et $m_{i+k}=1$. Soit $p=n-d$ et $J_{p,0}$ (resp. $J_{p,A^i}$) l'id\'eal de $\mathbb{K}[[t,X]]$ engendr\'e par $I_0$ et tous les d\'eterminants mineurs d'ordre p des jacobiennes de $p$ \'el\'ements de $I_0$ (resp. engendr\'e par $I_{i,A^i}$ et les d\'eterminants mineurs d'ordre $p$ des jacobiennes de $p$ \'el\'ements de $I_{i,A^i}$). Soit $\varphi \in \mathcal{A}_{(\mathbb{K}^n,0)}$ tel que $\varphi ^i=A^i$. Puisque la suite des multiplicit\'es de Nash est d\'ecroissante, on a $m_l\leq 2$, $0\leq l\leq i$. Un calcul \'el\'ementaire montre alors que:
$$ord_t(\Gamma _0^{*}(J_{p,0})) \geq i+ord_t(\Gamma _i^{*}(J_{p,A^i}))\geq i+1.$$
Soit $\beta '$ la fonction d'Artin-Greenberg de $\mathbb{K}[[t,X]]/J_{p,0}$ et $c\geq 1$ tel que $\beta '(i)\leq ci$. Notons $j=[i/c]$ la partie enti\`ere de $i/c$. Par d\'efinition de $\beta '$, on a donc $A^j \in \Pi^j(\mathcal{A}_{(sing(X),x)})$. Il en d\'ecoule que:
$$d_{i+k}=dim (P^{i+k}-P^{i+k}\cap \Pi _{i+k}^{-1}(P^i))\leq (d-1)j+d(i+k -j)=d(i+k)-j$$
Par suite $Lim_{i\rightarrow  +\infty}d(i+k)-d_{i+k}=+\infty$. Ce qui assure que la suite $[P^i]\mathbb{L}^{-di}$ est de Cauchy dans $\widehat{\mathcal{M}}_{\mathbb{K}}$ et donc convergente. \\
Soit maintenant $A^i\in \Pi^i(\mathcal{A}_{(X,x)})-P^i_{(X,x)}$. Soit $\varphi(t)=\sum_{k=1}^iA_kt^k+t^i\varphi _i(t)$ avec $\varphi _i(t)=\sum_{k=1}^{+\infty}A_{i+k}t^k$ un \'el\'ement de $\mathcal{A}_{(X,x)}$ qui rel\`eve $A^i$. Par construction du diagramme $(D_2)$, l'arc $\Gamma _i(t)=(t,\varphi _i(t))$ est trac\'e dans $(Z'_i,0)$. Ainsi le point $(1,A_{i+1})$ est dans l'ensemble des z\'eros de l'id\'eal initial $In(I_{i,A^i})$ de $I_{i,A^i}$. Ceci implique que pour tout $m\in \mathbb{N}$, $(m,0,\ldots,0)\not \in N_{i,A^i}$. Par cons\'equent puisque $A^i\not \in P^i$, la suite des multiplicit\'es de Nash $(m_0,\ldots,m_i)$ de $A^i$ est telle que $m_i\geq 2$. Proc\'edant comme plus haut en posant $j=[i/c]$, ceci implique que:
$$ord_t(\phi  ^{*}J_{p,0}) \geq i+1,\hspace{3pt}\forall \phi \in \mathcal{A}_{(\mathbb{K}^n,0)}\vert \phi ^i=A^i$$
Par cons\'equent $A^j \in \Pi ^j(\mathcal{A}_{(Sing(X),x)})$. Il en r\'esulte que:
$$l_i=dim(\Pi ^i(\mathcal{A}_{(X,x)}-P^i_{(X,x)})\leq (d-1)j+d(i-j)=di-j$$
Par suite $Lim_{i\longrightarrow \infty}di-l_i=+\infty$  et $[\Pi^i(\mathcal{A}_{(X,x)})-P^i_{(X,x)}]\mathbb{L}^{-di}$ converge vers z\'ero dans $\widehat{\mathcal{M}}_{\mathbb{K}}$ . Il en d\'ecoule la convergence des deux suites consid\'er\'ees et la coincidence de leurs limites. Le fait que cette limite soit \'egale \`a $\mu (\mathcal{A}_{(X,x)})$ est cons\'equence de 7.1 p.229 dans [D-L1] (ou bien peut \^etre consid\'erer comme une d\'efinition du volume motivique c.f. [D-L2]).
\endproof 
 
R\'ecemment S. Ishii et J.  Koll\'ar  ont consid\'er\'e l'hypersurface de $\mathbb{C}^5$ d\'efinie par $X_1^3+X_2^3+X_3^5+X_4^5+ X_5^6=0$ et montr\'e qu'elle fournissait un contre exemple en dimension 4 au \guillemotleft probl\'eme de Nash\guillemotright (c.f. [N], [I.K], [R1], [R2] par exemple, pour la formulation de celui-ci). Nous nous proposons ici, pour illustrer ce qui pr\'ec\`ede et voir le r\^ole jouer par les suites de Nash, de calculer par les m\'ethodes ci-dessus le volume motivique de telles hypersurfaces. On peut bien entendu faire le calcul de mani\`ere ind\'ependante en d\'eterminant  une r\'esolution des singularit\'es et utiliser \guillemotleft les formules de changement de variables \guillemotright de l'int\'egration motivique de [D-L1] ou [S].

\begin{ex} $\mathbb{K}=\mathbb{C}$. Soient $n\geq 3$, $k\geq 2$  et $(H_{n,k},0)$ le germe d'hypersurface \`a l'origine de $\mathbb{K}^{n+1}$ d\'efini par:
$$X_1^k+\ldots+X_n^k+Y^{2k}=0 ,\textrm{ (Le contre-exemple d'Ishii-Koll\'ar est le cas n=4, k=3)}$$
Posons pour $p\geq 1$, $V_{p,k}=\{(X_1,\ldots,X_p)\in \mathbb{K}^p/1+X_1^k+\ldots+X_p^k=0\}$.\\
On a alors:
\begin{eqnarray*}
\mu (\mathcal{A}_{(H_{n,k},0)})&  = &  ([V_{n-1,k}]+[V_{n-2,k}]+\ldots+k)\mathbb{L}\frac{\mathbb{L}-1}{\mathbb{L}^{n}-1}+[V_{n,k}]\frac{\mathbb{L}-1}{\mathbb{L}^{2n-1}-1}\\
 & + & ([V_{n-1,k}]+[V_{n-2,k}]+\ldots+k)\frac{(\mathbb{L}-1)^2}{(\mathbb{L}^{n-1}-1)(\mathbb{L}^{2n-1}-1)}
\end{eqnarray*}
\end{ex}
Nous indiquons bri\`evement les principales \'etapes des calculs. Pour cela, d\'esignons par $C_k$ et $W_k$ les constructibles:
\begin{eqnarray*}
C_k=\{X\in \mathbb{K}^n/X\neq 0\textrm{ et }X_1^k+\ldots+X_n^k=0\}\\
W_k=\{(X,Y)\in \mathbb{K}^{n+1}/Y\neq 0\textrm{ et }X_1^k+\ldots+X_n^k+Y^{2k}=0\}.
\end{eqnarray*}
Il est \'el\'ementaire de v\'erifier que dans $K_0(\mathcal{V}_{\mathbb{K}})$ on a:
$$[C_k]=([V_{n-1,k}]+[V_{n-2,k}]+\ldots+k)(\mathbb{L}-1)\textrm{ et }[W_k]=[V_{n,k}](\mathbb{L}-1).$$
 Soit maintenant $\varphi ^i=A^i=(A_1,\ldots,A_i)\in \mathbb{K}^{(n+1)i}$, $A_l=(B_l,Y_l)\in \mathbb{K}^{n+1}$ avec $B_l=(b_{l,1},\ldots,b_{l,n})$. Posons $v=Min\{l/A_l\neq 0\}$. Si $A^i\in P^i$ alors on a $v\leq i$, sans quoi $m_i=k$. L'\'equation locale de $(Z'_v,0)$ est:
$$(X_1+b_{v,1})^k+\ldots+(X_n+b_{v,n})^k+t^{kv}(Y+Y_v)^{2k}=0.$$
Deux \'eventualit\'es peuvent se pr\'esenter:\\
$\bullet $ $B_v\neq $ et alors $m_v>0$ si et seulement si $B_v\in C_k$. On alors $m_v=1$. Ainsi les \'el\'ments $A^i$ de $P^i$ au dessus de tels points sont une fibration localement triviale de fibre $\mathbb{L}^{n(i-v)}$ et de base $C_k \times \mathbb{L}$. Ainsi leur contribution \`a $[P^i]\mathbb{L}^{-ni}$ est $[C_k]\mathbb{L}.\mathbb{L}^{-nv}$. Faisant la somme pour $v$ variant de $1$ \`a $i$, et prenant la limite dans $\widehat{\mathcal{M}}_{\mathbb{K}}$ lorsque $i$ tend vers $+\infty$ on obtient:
$$[C_k]\mathbb{L}\times \frac{1}{\mathbb{L}^n-1}=([V_{n-1,k}]+[V_{n-2,k}]+\ldots+k)\mathbb{L}\frac{\mathbb{L}-1}{\mathbb{L}^{n}-1}$$
Ceci est la contribution des strates principales correspondant au suites de Nash $(k,\ldots,k,1,\ldots,1)$, le premier $1$ en position $v$ et $ord_t(\varphi ^i)=v$. Le d\'ecoupage de $[C_k]$ en $[V_{n-1,k}]+[V_{n-2,k}]+\ldots+k$ correspond aux diff\'erentes possibilit\'es pour $N^v=\alpha +\mathbb{N}^{n+1}$, avec $\alpha =(0,\beta )$ et $\vert\beta\vert =1$.\\
$\bullet $ $B_v=0$ et donc $Y_v\neq 0$. On consid\`ere le premier $l$ avec $v+l\leq i$ tel que $B_{v+l}\neq 0$. Si $A^i \in P^i$, on a certainement $l\leq v$. Sinon l'\'equation locale de $(Z'_{2v},0)$ est:
$$X_1^k+\ldots+X_n^k+(t^vY+t^vY_{2v}+\ldots+Y_v)^{2k}=0$$
et l'on a $m_{2v}=0$. On a alors deux sous-cas.\\
   $-\hspace{2pt} a)$ $l=v$. L'\'equation locale de $(Z'_{2v},0)$ est alors:
$$(X_1+b_{2v,1})^k+\ldots+(X_1+b_{2v,1})^k+(t^vY+t^vY_{2v}+\ldots+Y_v)^{2k}=0.$$
On a $m_{2v}>0$ si et seulement si $(B_{2v},Y_v)\in W_k$ et alors $m_{2v}=1$. Les \'el\'ements de $P^i$ au dessus de tels points sont une fibration localement triviale de fibre $\mathbb{L}^{n(i-2v)}$ et de base $W_k \times \mathbb{L}^{v}$. Leur contribution a $[P^i]\mathbb{L}^{-ni}$ est donc $[W_k]\mathbb{L}^{-(2n-1)v}$. Prenant la somme pour $v$ variant de $1$ \`a $[i/2]$ et passant \`a la limite on obtient:
$$[W_k]\times \frac{1}{\mathbb{L}^{2n-1}-1}=[V_{n,k}]\frac{\mathbb{L}-1}{\mathbb{L}^{2n-1}-1}$$
Ces points correspondent aux suites de Nash $(k,\ldots,k, 1,\ldots,1)$, le premier 1 en position $2v$ et $ord_t(\varphi ^i)=v$.\\
 $-\hspace{2pt} b)$ $1\leq l<v$ et $l+v\leq i$. On a en particulier $v\geq 2$. L'\'equation locale de $(Z'_{v+l},0)$ est 
$$(X_1+b_{1,v+l})^k+\ldots+(X_n+b_{n,v+l})^k+t^{k(v-l)}(t^lY+t^lY_{l+v}+\ldots+Y_v)^{2k}=0$$
On  alors $m_{v+l}>0$ si et seulement si $B_{v+l}\in C_k$ et alors $m_{v+l}=1$. Ainsi les \'el\'ements de $P^i$ au dessus de tels points sont une fibration localement triviale de fibre $\mathbb{L}^{n(i-v-l)}$ et de base $(\mathbb{L}-1)[C_k]\mathbb{L}^l$. Ainsi la contribution \`a $[P^i]\mathbb{L}^{-ni}$ est $(\mathbb{L}-1)[C_k]\mathbb{L}^{-n(v+l)}\mathbb{L}^l$. Faisant varier $v$ de $2$ \`a $i$ et $l$ de $1$ \`a $i_v=Min(v-1,i-v)$, on obtient:
$$[C_k](\mathbb{L}-1)\sum_{v=2}^i\mathbb{L}^{-nv}\sum_{l=1}^{i_v}\mathbb{L}^{-(n-1)l}=[C_k](\mathbb{L}-1)\sum_{v=2}^i\mathbb{L}^{-nv}\frac{1-\mathbb{L}^{-(n-1)i_v}}{\mathbb{L}^{n-1}-1}$$
Posons $j=[(i+1)/2]$. Nous obtenons:
$$[C_k](\mathbb{L}-1)\frac{1}{\mathbb{L}^{n-1}-1}(\sum_{v=2}^i\mathbb{L}^{-nv}-\mathbb{L}^{n-1}\sum_{v=2}^{j} \mathbb{L}^{-(2n-1)v}-\sum_{v=j+1}^i\mathbb{L}^{-(n-1)i}\mathbb{L}^{-v})$$
La troisi\`eme somme est de limite nulle. La limite de la diff\'erence des deux premi\`eres sommes vaut:
$$\frac{\mathbb{L}^{n-1}-1}{(\mathbb{L}^{n}-1)(\mathbb{L}^{2n-1}-1)}$$
Il vient donc (apr\`es simplification par $\mathbb{L}^{n-1}-1$):
$$[C_k]\frac{\mathbb{L}-1}{(\mathbb{L}^{n}-1)(\mathbb{L}^{2n-1}-1)}=([V_{n-1,k}]+[V_{n-2,k}]+\ldots+k)\frac{(\mathbb{L}-1)^2}{(\mathbb{L}^{n}-1)(\mathbb{L}^{2n-1}-1)}$$

\section{Condition de rang de Gabrielov et espace des arcs}
Dans cette section, $\mathbb{K}= \mathbb{C}$. Soit $F_{x}:(X,x) \rightarrow(Y,y)$ un germe de morphisme entre deux germes irr\'eductibles. Notons $\mathcal{F}_{x}:\mathcal{A}_{(X,x)} \rightarrow \mathcal{A}_{(Y,y)}$ le morphisme induit. Rappelons que nous d\'esignons par rang g\'en\'erique de $F_{x}$ le nombre : 
$$r_{1}(F_{x})=grk(F_{x})=Inf_{U } (Sup_{x'} Rang(F_{x'}))$$
 o\`u $x' \in Reg(X \cap U)$
et $U$ parcourt une base de voisinage $x$. Nous cherchons \`a \'etablir:
\begin{theo}
Les propri\'et\'es suivantes sont \'equivalentes:
\begin{itemize}
\item[1)] $grk(F_{x})=dim_{y}Y$;
\item[2)] Pour tout ouvert $\mathcal{V}$ de $\mathcal{R}_{(X,x)}$, $\mathcal{F}_{x}(\mathcal{V})$ contient un ouvert de  $\mathcal{R}_{(Y,y)}$;
\item[3)] $\mathcal{F}(\mathcal{R}_{(X,x)})$ contient un ouvert de $\mathcal{R}_{(Y,y)}$.
\end{itemize}
\end{theo} 
$\textit{Preuve:}$ 

$1) \Rightarrow 2)$

Soient $\varphi \in \mathcal{R}_{(X,x)}$ et $ i \in \mathbb{N}$, nous allons prouver que: $\mathcal{F}_{x}(B_{i}(\varphi))$ contient un ouvert de $\mathcal{R}_{(Y,y)}$. 
Nous commençons d'abord par nous ramener au cas o\`u $(X,x)$ est lisse. 

 \underline{R\'eduction au cas o\`u $(X,x)$ est lisse:} 

\vspace{10pt}

Consid\'erons un plongement $(X,x) \rightarrow (\mathbb{K}^{m},0)$. Soient $I \subset \mathbb{K}\{x_{1},\ldots,x_{m}\}$ l'id\'eal d\'efinissant $(X,x)$, $l=haut(I)$. Notons $J_{l}$ l'id\'eal engendr\'e par $I$ et tous les d\'eterminants jacobiens d'ordre $l$ des $l$-upplets d'\'el\'ements de $I$. Puisque $\varphi \in \mathcal{R}_{(X,x)}, ord(\varphi^{*}J_{l})=u<+\infty$. 
Par le proc\'ed\'e utilis\'e au th\'eor\`eme 4.1, on peut trouver un germe $(Z,z)$ lisse, un germe de morphisme: $\Pi_{z}:(Z,z) \rightarrow (X,x)$ obtenu par une succession de $u$ \'eclatements ponctuels et $\phi \in \mathcal{A}_{(Z,z)}$ qui rel\`eve $\varphi$ \`a travers $\Pi_{z}$. Soit alors $G_{z}:(Z,z) \underrightarrow{\Pi_{z}} (X,x) \underrightarrow{F_{x}} (Y,y) $. Il est alors \'el\'ementaire de v\'erifier que $grk(G_{z})=grk(F_{x})$.  Si l'assertion 2) est prouv\'ee pour $\mathcal{G}_{z}$, elle le sera pour $\mathcal{F}_{x}$. En effet, soit $i\in \mathbb{N}$, s'il existe $j\in \mathbb{N}$ et $\theta \in \mathcal{R}_{(Y,y)}$ telle que:
$$B_{j}(\theta) \subset \mathcal{G}_{z}(B_{i}(\phi)),$$ 
on aura:
$$B_{j}(\theta) \subset \mathcal{G}_{z}(B_{i}(\phi))=\mathcal{F}_{x}(\Pi_{z}(B_{i}(\phi)) \subset \mathcal{F}_{x}(B_{i}(\varphi)).$$  
Nous supposerons donc d\'esormais que $(X,x)=(\mathbb{K}^{m},0)$. Soit alors $d=dim_{y}Y$. Plongeons $(Y,y)$ dans $(\mathbb{K}^{n},0)$. Nous noterons $F_{1}(x_{1},\ldots,x_{m}),\ldots,$  $F_{n}(x_{1},\ldots,x_{m})$ les composantes de $F_{x}$, et par $F$ un repr\'esentant de $F_{x}$. Pour un choix g\'en\'erique de coordonn\'ees lin\'eaires \`a l'origine de $(\mathbb{K}^{n},0)$, on peut supposer \^etre dans la situation qui suit. 

Soit $P_{0}:(Y,0) \rightarrow (\mathbb{K}^{d},0)$, le germe de morphisme induit par la projection : 
$$\mathbb{K}^{n} \rightarrow \mathbb{K}^{d}$$
$$\hspace{20pt} (z_{1},\ldots,z_{n})\rightarrow (z_{1},\ldots,z_{d})$$

Le morphisme: 
$$P_{0}^{*}:\hspace{5pt}\mathcal{O}_{d}=\mathbb{K} \{z_{1},\ldots,z_{d}\} \longrightarrow   \mathcal{O}_{Y,0}=\mathcal{O}_{n}/J$$
est injectif et fini. 

Si $k$ est la dimension du corps de fractions de $\mathcal{O}_{Y,0}$ sur le corps des fractions de $\mathcal{O}_{d}$, le polyn\^ome minimal de $\overline{z}_{d+1}$ est un polyn\^ome distingu\'e, $P(Z)=Z^{k}+ \sum_{l=1}^{k}a_{l} Z^{k-l} $ \`a coefficients dans $\mathcal{O}_{d}$.  

Le discriminant $\Delta(z_{1},\ldots,z_{d})$ de ce polyn\^ome est non nul et: 
$$\overline{\Delta}. \mathcal{O}_{Y,0} \subset \mathcal{O}_{d} +\mathcal{O}_{d}.\overline{z}_{d+1}+ \ldots +\mathcal{O}_{d}.\overline{z}_{d+1}^{k-1}.$$
En particulier, $z_{d+1}^{k}+ \sum_{l=1}^{k}a_{l} z_{d+1}^{k-l} \in J$, et pour $2 \leq j \leq n-d$, il existe des $a_{l,j} \in \mathcal{O}_{d}$ tels que: 
$$ \Delta z_{d+j} - \sum_{l=1}^{k-1}a_{l,j} z_{d+1}^{l} \in J.$$
D'autre part, le choix de la projection lin\'eaire $P_{0}$ \'etant g\'en\'erique et $grk(F_{x})=d$, on peut supposer quitte à renum\'eroter les coordonn\'ees dans $\mathbb{K}^{m}$ que:
$$J_{1}(x_{1},\ldots,x_{m})=det(\frac{\partial F_{l}} {\partial x_{j}}(x_{1},\ldots,x_{m}))_{1 \leq l \leq d, 1 \leq j \leq d   } $$ 
est non identiquement nul.  

Notons $(D,0) \subset (\mathbb{K}^{d},0)$ le germe en $0$ des z\'eros de $\Delta$. On a $\Delta \circ F_{x} \neq 0$, sans quoi $F_{x} ((\mathbb{K}^{m},0)) \subset P_{0}^{-1}(D,0)$ d'o\`u $grk(F_{x}) \leq dim_{0}P_{0}^{-1}(D,0) = d-1 <d$. 

De la m\^eme façon $(V,0)=F_{x}^{-1}(Sing(Y),0) \subset (\mathbb{K}^{m},0)$ est un sous-ensemble analytique propre de $(\mathbb{K}^{m},0)$. Soit $I_{V}$ l'id\'eal r\'eduit d\'efinissant $(V,0)$ et $h \in I_{V}$, $h \neq 0$. Posons $k_{1}=m_{0}(h)$, $k_{2}=m_{0}(J_{1})$, $k_{3}=m_{0}(\Delta \circ  F_{x})$, o\`u $m_{0}(\hspace{3pt} )$ d\'esigne la multiplicit\'e ou l'ordre \`a l'origine. 

Soient $\phi \in \mathcal{A}_{(\mathbb{K}^{m},0)}$ et $i \in \mathbb{N}$, en vertu du lemme 4.9 , il existe $\alpha \in B_{i}(\phi)$ telle que:
$$ord(h(\alpha(t)). J_{1}(\alpha(t)). \Delta(F(\alpha(t))) \leq (k_{1} + k_{2} + k_{3}) (i+1)=s.$$
Posons $\theta(t)=F(\alpha(t))$. Nous allons voir que:

\vspace{10pt}

\begin{itemize}
\item[-] $B_{2s+1}(\theta) \subset \mathcal{F}_{x}(B_{i}(\phi))=\mathcal{F}_{x}(B_{i}(\alpha))$;
\item[-] $B_{2s+1}(\theta) \subset \mathcal{R}_{(Y,0)}$.
\end{itemize}  

\vspace{10pt}

Soit $\theta'=( \theta '_{1}, \ldots,\theta '_{d}, \theta '_{d+1}, \ldots, \theta '_{n}) \in B_{2s+1}(\theta) \subset \mathcal{A}_{(Y,0)}$. Consid\'erons le syst\`eme d'\'equations implicites: 
$$ G_{1}(t,y_{1},\ldots,y_{m})= F_{1}(\alpha_{1}(t)+y_{1},\ldots,\alpha_{m}(t)+y_{m}) - \theta '_{1}(t)=0$$
$$\vdots$$
$$ G_{d}(t,y_{1},\ldots,y_{m})= F_{d}(\alpha_{1}(t)+y_{1},\ldots,\alpha_{m}(t)+y_{m}) - \theta '_{d}(t)=0.$$ 
Soit $I$ l'id\'eal de $\mathcal{O}_{1}$ engendr\'e par les mineurs d'ordre $d$ de la jacobienne $G'_{y}(t,0)$. On a $ J_{1}(\alpha(t))\in I$. Soit $l=ord(I) \leq ord(J_{1}(\alpha(t)) \leq s$. 

Par hypoth\`ese,  $G(t,0) \in (t)^{2s+1}. \mathcal{O}_{1}^{d}$. Donc $G(t,0) \in \oplus _{d} (t)^{2(s-l)+1}.I^{2}$. Par le th\'eor\`eme des fonctions implicites de J.C. Tougeron, il existe $y_{1}(t), \ldots,y_{m}(t)$  $ \in (t)^{2(s-l)+1}. (t)^{l} \subset (t)^{s+1}$ telles que: 
$$ F_{j}(\alpha (t) + y(t)) = \theta'_{j}(t),\hspace{3pt} 1 \leq j \leq d.$$ 
Posons $\rho (t) = \alpha (t) + y(t)$. On a donc $\rho \in B_{s}(\alpha) \subset B_{i}(\alpha)$. 

Pour conclure, il nous reste \`a voir que pour $d+1 \leq j \leq n$, on a: $F_{j}(\rho (t))= \theta '_{j}(t)$. 

\vspace{10pt}

Pour cela remarquons que, puisque $\theta '\in B_{2s+1}(\theta) \subset \mathcal{A}_{(Y,0)}$, $\Delta (\theta ')$ et $\Delta (\theta)$ coinc\"{i}dent jusqu'\`a l'ordre $2s$. Mais $ord(\Delta (\theta)) \leq s$, donc $ord(\Delta (\theta ')) \leq s$. 

Posons: $\tilde{\theta} '(t)=(\theta '_{1}(t), \ldots,\theta '_{d}(t))$. Maintenant, puisque $F$ est \`a valeur dans $(Y,0)$ et $(F_{1}(\rho (t)), \ldots,F_{d}(\rho (t))) = \tilde{\theta} '(t)$. 

On a:    
$$F_{d+1}(\rho(t))^{k}+\sum_{l=1}^{k}a_{l}(\tilde{\theta}'(t))F_{d+1}^{k-l}(\rho (t))=0$$
$$\theta '_{d+1}(t)^{k}+\sum_{l=1}^{k}a_{l}(\tilde{\theta}'(t))(\theta'_{d+1})^{k-l}(t)=0$$
Par suite $F_{d+1}(\rho(t))$ et $\theta'_{d+1}(t)$ sont donc deux racines du polyn\^ome $P(U)=U^{k}+\sum_{l=1}^{k}a_{l}(\tilde{\theta}'(t))U^{k-l}$. D'apr\`es le th\'eor\`eme de Puiseux, il existe, $u_{1}(t^{1/k!})$, $\ldots$, $ u_{k}(t^{1/k!})$ analytiques complexes en $t^{1/k!}$  solutions de $P(U)=0$. Supposons que $F_{d+1}(\rho(t))$ $\neq \theta_{d+1} '(t)$, alors:
$$(F_{d+1}(\rho(t^{k!}))-\theta'_{d+1}(t^{k!}))^{2} \textrm{ divise }\prod _{l<j}(u_{l}(t)-u_{j}(t))^{2}=\Delta(\tilde{\theta}'(t^{k!})).$$
Par suite:
$$ord(F_{d+1}(\rho(t))-\theta'_{d+1}(t))k!\leq ord(\Delta(\tilde{\theta}'(t))k!\leq sk!.$$
Or ceci est absurde car:
$$ord(F_{d+1}(\rho(t))-F_{d+1}(\alpha(t))\geq s+1 \textrm{ et } ord(\theta'_{d+1}(t)-F_{d+1}(\alpha(t)))\geq 2s+1.$$
Par suite $F_{d+1}(\rho(t))=\theta '_{d+1}(t)$. Ensuite, pour $j \geq 2$, on a: 
$$ \Delta(\tilde{\theta}'(t)) F_{d+j} (\rho (t)) = \sum_{l=0}^{k-1}a_{l,j}(\tilde{\theta}'(t))F_{d+1}^{l}(\rho (t))$$
$$ = \sum_{l=0}^{k-1}a_{l,j}(\tilde{\theta}'(t))(\theta '_{d+1})^{l}(t)= \Delta(\tilde{\theta}'(t)) \theta '_{d+j}(t).$$
Par suite, $F_{d+j}(\rho (t))=\theta '_{d+j}(t)$, $j \geq 2$. Donc $\rho$ rel\`eve $\theta '$ par $F$.  Enfin notons que $h(\rho(t)) \neq 0$ car $ord(h(\rho(t))- h(\alpha(t))) \geq s+1$ et $ord(h(\alpha(t)) \leq s$.
Par cons\'equent, $\rho(t)$ n'est pas trac\'ee dans $F^{-1}(Sing(Y),0)$. Donc \`a fortiori $F(\rho(t))= \theta '(t)$ n'est pas trac\'ee dans $(Sing(Y),0)$.

\vspace{10pt}    

Nous prouvons \`a pr\'esent l'implication $3) \Rightarrow 1)$. 

\vspace{5pt}

Consid\'erons $F_{x}: (X,0) \rightarrow (Y,0)$ tel que $\mathcal{F}(\mathcal{R}_{(X,0)})$ contienne un ouvert de $\mathcal{R}_{(Y,0)}$. Posons $d=dim_{0}Y$ et $P_{0}: (Y,0) \rightarrow (\mathbb{K}^{d},0)$ une projection lin\'eaire g\'en\'erique. D'apr\`es l'implication $1) \Rightarrow 2)$,  puisque $grk(P_{0})=d$, pour tout ouvert $\mathcal{V}$ de $\mathcal{R}_{(Y,0)}$, $\mathcal{P}_{0}(\mathcal{V})$ contient un ouvert de $\mathcal{A}_{(\mathbb{K}^{d},0)}$. Par suite, posant $F'=P_{0} \circ F$, $\mathcal{F}'(\mathcal{R}_{(X,0)})$ contient un ouvert de $\mathcal{A}_{(\mathbb{K}^{d},0)}$. Mais $grk(F') \leq grk(F)$. L'implication $3) \Rightarrow 1)$ sera donc prouv\'ee si l'on montre que $grk(F')=d$. Autrement dit, on peut supposer que notre germe d'application $F: (X,0) \rightarrow (\mathbb{K}^{d},0)$ est telle  que $\mathcal{F}(\mathcal{R}_{(X,0)})$ contient un ouvert de $\mathcal{A}_{(\mathbb{K}^{d},0)}$.

Consid\'erons une boule de rayon $r$ suffisamment petit, $B(0,r) \subset \mathbb{K}^{m}$ pour que:
$$X \cap B(0,r) = \{z \in B(0,r) / f_{1}(z)= \ldots=f_{p}(z)=0 \},$$
avec $f_{1},\ldots,f_{p}$ analytiques sur $B(0,r)$.
$$Sing(X) \cap B(0,r) =\{z \in B(0,r) / f_{1}(z)= \ldots=f_{p}(z)=0, $$ 
 $$ det(\frac{\partial{f_{i_{l}}}} {\partial{z_{j_{k}}}})_{1 \leq l \leq s, 1 \leq k \leq s }(z)=0,\forall \underline{i}=(i_{1},\ldots,i_{s}), \underline{j}=(j_{1},\ldots,j_{s})\}$$
o\`u $m-s = dim(X,0)$.

Soit $r_{n}$ une suite croissante telle que $0<r_{n}<r$ et $r_{n}\rightarrow r$. Posons:
$$R_{n,l}=\{z \in \mathbb{K}^{m}/\vert z \vert\leq r_{n},\hspace{1pt} z\in X\cap B(0,r),\sum_{\underline{i},\underline{j}}\vert det(\frac{\partial f_{i_{j}}}{\partial z_{j_{k}}})_{1 \leq j,k \leq s}(z)\vert \geq \frac{1}{l}\}.$$          
On a \'evidemment: $Reg(X)\cap B(0,r)=\cup _{l,n \in \mathbb{N}^{*}}R_{n,l}$. Par suite $F(Reg(X)\cap B(0,r))=\cup_{n,l}F(R_{n,l})$ est r\'eunion d\'enombrable de compacts de $\mathbb{K}^{d}$. Donc l'ensemble $F(Reg(X)\cap B(0,r))$ est mesurable pour la mesure de Lebesgue sur $\mathbb{K}^{d}$. Supposons $grk(F)<d$. Ecrivons $Reg(X)\cap B(0,r)=\cup_{n\in \mathbb{N}}D_{n}$ o\`u chaque $D_{n}$ est un domaine de carte. Posons $F_{n}=F/D_{n}$. En tout point de $D_{n}$ le rang de la jacobienne complexe de $F_{n}$ est strictement inf\'erieur \`a $d$. Donc le rang de la jacobienne r\'eelle de $F_{n}$ est strictement inf\'erieur \`a $2d$. Tout point de $F_{n}(D_{n})$ est donc valeur critique r\'eelle de $F_{n}$. Par le th\'eor\`eme de Sard, $F_{n}(D_{n})$ est donc de mesure nulle. Il en est donc de m\^eme de $F(Reg(X)\cap B(0,r))=\cup F_{n}(D_{n})$. Nous allons voir que ceci est contradictoire avec le fait que $\mathcal{F}(\mathcal{R}_{(X,0)})$ contienne un ouvert de $\mathcal{A}_{(\mathbb{K}^{d},0)}$.

\vspace{7pt}

Soient $\varphi \in \mathcal{A}_{(\mathbb{K}^{d},0)}$ et $i \in \mathbb{N}$ tel que $B_{i}(\varphi)
\subset \mathcal{F}(\mathcal{R}_{(X,0)})$. Posons $\varphi (t)=\sum_{k=1}^{+ \infty}A_{k}t^{k}$, $A_{k}=(a_{k}^{1},\ldots,a_{k}^{d})\in \mathbb{K}^{d}$. Consid\'erons alors l'application:
$$A:\hspace{5pt} \mathbb{K}\times\mathbb{K}^{d}\longrightarrow \mathbb{K}^{d}$$
$$\hspace{20pt} (t,x)\longrightarrow A(t,x)=\sum_{k=1}^{i}A_{k}t^{k}+x.t^{i+1}.$$
$F(Reg(X)\cap B(0,r))$ est r\'eunion d\'enombrable de compacts, donc de ferm\'es. Il en est donc de m\^eme de $V=A^{-1}(F(Reg(X)\cap B(0,r)))$. Soient $l\geq 1$ et $C_{l}=\{z \in \mathbb{K}/\vert z \vert =1/l \}$ le cercle de centre l'origine et de rayon $1/l$. Consid\'erons:
$$V_{l}=V\cap(C_{l}\times \mathbb{K}^{d}).$$
$V_{l}$ est donc r\'eunion d\'enombrable de compacts. Soit $X_{l}=P(V_{l})$, o\`u $P: \mathbb{K}\times \mathbb{K}^{d}\rightarrow \mathbb{K}^{d}$ est la projection canonique. $X_{l}$ est donc mesurable. Maintenant, le fait que $B_{i}(\varphi)
\subset \mathcal{F}(\mathcal{R}_{(X,0)})$ implique que:
$$\mathbb{K}^{d}=\cup_{l\geq 1} X_{l}\hspace{20pt}(\mathbf{*})$$
 En effet, soit $x \in \mathbb{K}^{d}$, l'arc $t\rightarrow \alpha(t)=A(t,x)$ est un \'el\'ement de $B_{i}(\varphi)$. Il se rel\`eve donc par $F$ en un arc $\rho \in \mathcal{R}_{(X,0)}$. Donc pour $\vert t \vert =1/l$ assez petit, $\alpha(t)=A(t,x)$ $\in F(Reg(X\cap B(0,r))$ et donc $x \in X_{l}$. Maintenant par $(\mathbf{*})$, l'un des $X_{l}$, disons $X_{l_{0}}$ est de mesure non nulle. Il en d\'ecoule que l'ensemble:
$$Z_{l_{0}}=\{z \in \mathbb{K}^{d}/ z=\sum_{k=1}^{i}A_{k}(\frac{1}{l_{0}})^{k}+u(\frac{1}{l_{0}})^{i+1}, \hspace{3pt} u\in X_{l_{0}} \}$$
est de mesure non nulle. Mais par d\'efinition, $Z_{l_{0}}\subset F(Reg(X)\cap B(0,r))$. Ce dernier ne peut donc \^etre de mesure nulle.
\endproof

\vspace{10pt}

Nous terminons cette section par quelques remarques et commentaires. Pour cela si $(X,x)$ est un germe d'espace analytique, nous noterons par $\widehat{\mathcal{O}} _{X,x}$ le compl\'et\'e $\mathcal{M}_{X,x}$ addique de $\mathcal{O}_{X,x}$. De la m\^eme façon si $F_{x}: (X,x) \rightarrow (Y,y)$ est un germe de morphisme, on notera $\widehat{F}_{x}^{*}: \widehat{\mathcal{O}}_{Y,y} \rightarrow \widehat{\mathcal{O}} _{X,x}$ le morphisme induit par $F_{x}$. On a alors le corollaire suivant qui est un cas particulier des th\'eor\`emes de A.M. Gabrielov et S. Izumi (cf. [I4]).  

\begin{cor}
Soit $(X,x)$ un germe d'espace r\'eduit et $F_{x}: (X,x) \rightarrow (\mathbb{C}^{d},0)$ tel que $grk(F_{x})=d$. Alors: 
\begin{itemize}
\item[1)] Si $S \in \widehat{\mathcal{O}} _{d}$ est tel que $\widehat{F_{x}}^{*}(S) \in \mathcal{O}_{X,x}$, alors $S \in \mathcal{O}_{d}$;
\item[2)] $\exists a,b \in \mathbb{R}^{+} / \forall S \in \widehat{\mathcal{O}} _{d}, m(\widehat{F_{x}}^{*}(S)) \leq a m(S) + b$ o\`u $"m"$ d\'esigne la multiplicit\'e ou l'ordre.
\end{itemize}
\end{cor}

\noindent \textit{Preuve}:

Puisque le rang g\'en\'erique de $F_{x}$ est le maximum des rangs g\'en\'eriques des restrictions de $F_{x}$ aux composantes irr\'eductibles de $(X,x)$, on peut  supposer que celui-ci est irr\'eductible. Mais alors d'apr\`es le th\'eor\`eme 1.4, il existe $\varphi \in \mathcal{A}_{(\mathbb{C}^{d},, 0)}$ et $i \in \mathbb{N}$ tels que $B_{i}(\varphi) \subset \mathcal{F}_{x}(\mathcal{R}_{(X,x)}) \subset \mathcal{F}_{x} (\mathcal{A}_{(X,x)})$. Soit alors $S \in \widehat{\mathcal{O}} _{d}$, tel que $\widehat{F_{x}}^{*}(S) \in \mathcal{O}_{(X,x)}$. Par suite, pour tout $\theta \in \mathcal{A}_{(X,x)}$, $\widehat{F_{x}}^{*}(S)(\theta) \in \mathcal{O}_{1}$, c'est \`a dire $S(\varphi ')  \in \mathcal{O}_{1}$, $\forall \varphi ' \in \mathcal{F}_{x} (\mathcal{A}_{(X,x)})$. Donc $S(\varphi ')  \in \mathcal{O}_{1}$, $\forall \varphi ' \in B_ {i}(\varphi)$. Or ceci implique que $S \in \mathcal{O}_{d}$ par le lemme 1.4 p.121 de [T4]. De la m\^eme façon, soit $k=m(\widehat{F_{x}}^{*}(S))$. Pour tout $\theta \in \mathcal{A}_{(X,x)}$, on a: $ord(\widehat{F_{x}}^{*}(S)(\theta)) \geq k$. Par suite, $Min_{\varphi ' \in B_ {i}(\varphi)} ord(S(\varphi ')) \geq k$. Mais d'apr\`es le lemme 4.9 de la section pr\'ec\'edente, ce minimum est inf\'erieur ou \'egal \`a $m(S)(i+1)$. Donc $m(\widehat{F_{x}}^{*}(S)) \leq \frac{1}{i+1} m(S)$. On notera donc que les constantes $a,b$ de $2)$ sont d\'etermin\'ees par le diam\`etre d'une boule incluse dans $\mathcal{F}_{x} (\mathcal{A}_{(X,x)})$. La preuve de $1) \Rightarrow 2)$ dans 1.4 fournit un moyen effectif de d\'eterminer un tel diam\`etre.
\endproof

\begin{rem}

Les th\'eor\`emes de A.M. Gabrielov et S.Izumi, qui sont respectivement les assertions de $1)$ et $2)$ dans le cas o\`u $(Y,y)$ est un irr\'eductible quelconque, sont beaucoup plus profonds. L'assertion $1)$ dans le cadre g\'en\'eral f\^ut d'abord prouv\'ee par A.M. Gabrielov. Puis S. Izumi donna une preuve de $1)$ notablement simplifi\'ee dans [I4]. L'assertion $2)$ dans le cas g\'en\'eral est due \`a S. Izumi [I1][I2][I3]. Le cas que nous traitons dans 5.2.1 a \'et\'e obtenu ind\'ependamment des façons pr\'ec\'edentes et de mani\`ere disjointes par P.M. Eakin-J.Harris, B. Malgrange, R. Moussu-J.C. Tougeron. On se ram\`ene, en g\'en\'eral, au cas o\`u $(X,x)$ est lisse en utilisant la r\'esolution des singularit\'es. Nous renvoyons le lecteur \`a [I4] pour une vision g\'en\'erale historique.
\end{rem}

Soit maintenant $(X,x)$ un germe irr\'eductible, $\varphi \in \mathcal{R}_{(X,x)}$ et $i\in \mathbb{N}$. Condid\'erons les deux assertions suivantes :
$$(T) \hspace{40pt}\forall S \in \widehat{\mathcal{O}}_{X,x}, (\forall \varphi ' \in B_{X,i}(\varphi) ,\hspace{5pt} S(\varphi ') \in \mathcal{O}_{1} \Rightarrow S \in \mathcal{O}_{X,x})     $$ 
$$ (I) \hspace{5pt} \exists a,b \in \mathbb{R}^{+}, \forall S \in \widehat{\mathcal{O}}_{X,x}, Min_{\varphi ' \in B_ {X,i}(\varphi)} (ord(S(\varphi '))) \leq a m(S) + b.$$
\footnotetext[1]{$(T)$ f\^ut conjectur\'ee par J.C. Tougeron dans [T4] avant qu'il n'ait connaissance du th\'eor\`eme de A.M. Gabrielov.}
\footnotetext[2]{$(I)$ est du \`a S. Izumi [I5] utilisant la r\'esolution des singularit\'es et son th\'eor\`eme mentionn\'e plus haut.} 
Les assertions $(T)$ et $(I)$ sont des cons\'equences tr\`es faciles de l'existence de la r\'esolution des singularit\'es et des th\'eor\`emes de A.M. Gabrielov et S. Izumi. Elles signifient grossi\`erement que le comportement d'un \'el\'ement de $\widehat{\mathcal{O}}_{X,x}$ sur un ouvert de $\mathcal{R}_{(X,x)}$ d\'ecrit son comportement g\'en\'eral. C'est pourquoi on peut consid\'erer que de tels ouverts sont repr\'esentatifs de $\mathcal{A}_{(X,x)}$. D'autre part, une preuve directe de ces assertions $(T)$ et $(I)$  impliquerait les th\'eor\`emes de A.M. Gabrielov et S. Izumi via 1.4 comme en 5.2. Pour prouver $(T)$ et $(I)$ directement, il suffirait d'\'etablir les propri\'et\'es $1)$ et $2)$ de 5.2 pour $\Pi: (Z,p) \rightarrow (X,x)$ l'\'eclatement de $(X,x)$ de centre $x$. Donc les th\'eor\`emes de A.M. Gabrielov et S. Izumi se r\'eduisent au cas apparemment simple ci-dessus. Contrairement aux apparences, ce cas est sans doute aussi compliqu\'e que le cas g\'en\'eral. Cependant, une d\'etermination "effective" de $a$ et $b$ comme dans 5.2 dans le cas de l'\'eclatement $\Pi: (Z,p) \rightarrow (X,x)$ de centre $x$ conduirait \`a des "estimations de Chevalley lin\'eaires uniformes" comme demand\'ees dans [B-M3] (question 1.28 p. 743).

\vspace{20pt}

\noindent \boldmath{\Large  Bibliographie}
\vspace{10pt}

\noindent [A1] M. Artin, \textit{Algebraic approximation of structures over complete local rings}, Pub. Math. I.H.E.S. $\mathbf{36}$ (1969), 23-58.

\noindent [B] B.M. Bennett, \textit{On the characteristic function of a local ring}, Ann. of Math. (2)  $\mathbf{91}$ (1970), 25-87.

\noindent [B-M1] E. Bierstone and P.D. Milman, \textit{Relations among analytic functions I}, Ann. Inst. Fourier (Grenoble), (1) $\mathbf{37}$ (1987), 187-239.

\noindent [B-M2] E. Bierstone and P.D. Milman, \textit{Uniformization of analytic spaces}, Jour. A.M.S. $\mathbf{2}$ (1989), 801-836.

\noindent [B-M3] E. Bierstone and P.D. Milman, \textit{Geometric and differential properties of subanalytic sets}, Ann. of Math. $\mathbf{147}$ (1998), 731-785.

\noindent [D-L1] J. Denef and F. Loeser, \textit{Germs of arcs on singular algebraic varieties  and motivic integration}, Invent. Math. $\mathbf{135}$ (1999), 201-232.

\noindent [D-L2] J. Denef and F. Loeser, \textit{Geometry on arc spaces of algebraic varieties}, in European Congress of Mathematics, vol I (Barcelona, 2000), Prog. Math., vol. 201, Birkh\"{a}user, 2001, 327-348.
 
\noindent [Ga] A.M. Gabrielov, \textit{Formal relations between analytic functions}, Izv. Akad. Nauk. SSSR. $\mathbf{37}$ (1973), 1056-1088.

\noindent [Gr] M. Greenberg, \textit{Rational points in henselian discrete valuation rings}, Pub. Math. I.H.E.S. $\mathbf{31}$ (1966), 59-64.

\noindent [G-S-L-J] G. Gonzalez-Sprinberg and M. Lejeune-Jalabert, \textit{Sur l'espace des cour\-bes trac\'ees sur une singularit\'e}, Prog. Math. $\mathbf{134}$ (1996), 8-32.

\noindent [H1] M. Hickel, \textit{Fonction de Artin et germes de courbes trac\'ees sur un germe d'espace analytique}, Amer. J. Math. $\mathbf{115}$ (1993), 1299-1334.

\noindent [H2] M. Hickel, \textit{Calcul de la fonction d'Artin-Greenberg d'une branche plane}, Pacific Journal of math. $\mathbf{213}$ (2004), 37-47, et Preprint Université Bordeaux I N°145, Avril 2002.

\noindent [I-K] S. Ishii and J. Koll\'ar, \textit{The Nash problem on arc families of singularities}, Duke Math. J. $\mathbf{120}$ (2003), 601-620.

\noindent [I1] S. Izumi, \textit{Linear complementary inequalities for orders of germs of analytic functions}, Invent. Math. $\mathbf{65}$ (1982), 459-471.

\noindent [I2] S. Izumi, \textit{A measure of integrity for local analytic algebras}, Pub. R.I.M.S, Kyoto Univ. $\mathbf{21}$ (1985), 719-735.

\noindent [I3] S. Izumi, \textit{Gabrielov's Rang Condition is equivalent to an inequality of reduced orders}, Math. Ann. $\mathbf276$ (1986), 81-89.

\noindent [I4] S. Izumi, \textit{The rank condition and convergence of formal  functions}, Duke Math. J. $\mathbf{59}$ (1989), 241-264.

\noindent [I5] S. Izumi, \textit{Note on convergence and orders of vanishing along curves}, Manuscripta math. $\mathbf66$ (1990), 261-275.

\noindent [L-J] M. Lejeune-Jalabert, \textit{Courbes trac\'ees sur un germe d'hypersurface}, Amer. J. Math. $\mathbf{112}$ (1990), 525-568.

\noindent [L-J-T] M. Lejeune-Jalabert and B. Teissier, \textit{Contribution \`a l'\'etude des singularit\'es du point de vue du polygone de Newton}, th\`ese universit\`e Paris VII, 1973.

\noindent [Lo] Looijenga E., Motivic Measures, \textit{S\'eminaire Bourbaki 1999/2000} $\mathbf{874}$, Ast\'e\-risque $\mathbf{276}$ (2002), 267-297.

\noindent [N] J.F. Nash Jr, \textit{Arc structure of singularities}, Duke Math. Journ. $\mathbf{81}$ (1995), 31-38.

\noindent [P] W. Pawlucki, \textit{On Gabrielov's regularity condition for analytic mappings}, Duke Math. J. $\mathbf{65}$ (1992), 299-311.

\noindent [R1] A.J. Reguera, \textit{Families of arcs on rational surface singularities}, Manuscripta Math. $\mathbf{88}$ (1995), 321-333.

\noindent [R2] A.J. Reguera, \textit{Image of the Nash map in terms of wedges}, C.R. Acad. Sci. Paris, Ser. I $\mathbf{338}$ (2004), 385-390.

\noindent [S] J. Sebag, \textit{Int\'egration motivique sur les sch\'emas formels}, \`a para\^{\i}tre au Bull. S.M.F. (2004).

\noindent [T1] J.C. Tougeron, \textit{Id\'eaux de fonctions diff\'erentiables }, Ann. Inst. Fourier $\mathbf{18}$ (1968), 177-240.

\noindent [T2] J.C. Tougeron, \textit{Id\'eaux de fonctions diff\'erentiables }, Ergeb. der Mathematik $\mathbf{71}$, Springer Verlag 1972.

\noindent [T3] J.C. Tougeron, \textit{An extension of Whithney's spectral theorem}, Pub. Math. I.H.E.S. $\mathbf{40}$ (1971), 139-148.

\noindent [T4] J.C. Tougeron, \textit{Courbes analytiques sur un germe d'espace analytique et applications}, Ann. Inst. Fourier $\mathbf{26}$ (1976), 117-131.\\

\vspace{5pt}
\noindent M. Hickel, Universit\'e Bordeaux 1, Laboratoire d'analyse et g\'eom\'etrie\\
U.M.R. 5467 du C.N.R.S., 351 cours de la lib\'eration\\
F-33405 TALENCE Cedex\\
\noindent e-mail address: hickel@math.u-bordeaux.fr
\end{document}